\documentclass[a4paper,12pt]{article}
\usepackage{amsmath,amsfonts,amssymb,amsthm}
\usepackage{latexsym,graphicx}
\usepackage{xypic}
\xyoption{all}
\def\Label{\label}
\begin{document}


\newtheorem{proposition}{\sc Proposition}[section]
\newtheorem{lemma}[proposition]{\sc Lemma}
\newtheorem{corollary}[proposition]{\sc Corollary}
\newtheorem{theorem}[proposition]{\sc Theorem}

\theoremstyle{definition}
\newtheorem{definition}[proposition]{\sc Definition}
\newtheorem{example}[proposition]{\sc Example}

\theoremstyle{remark}
\newtheorem{remark}[proposition]{\sc Remark}

\def\proof{{\sc Proof.~~}}
\def\endproof{\hbox{$\sqcup$}\llap{\hbox{$\sqcap$}}\vskip 4pt plus2pt}

\newcommand{\Section}{\setcounter{definition}{0}\section}
\renewcommand{\theequation}{\thesection.\thesubsection.\arabic{equation}}
\newcounter{c}
\renewcommand{\[}{\setcounter{c}{1}$$}
\newcommand{\etyk}[1]{\vspace{-7.4mm}$$\begin{equation}\Label{#1}
\addtocounter{c}{1}}
\renewcommand{\]}{\ifnum \value{c}=1 $$\else \end{equation}\fi}
\setcounter{tocdepth}{3}

\def\Rhom#1#2#3{{{\rm Hom}\sp{#1}(#2,#3)}} 
\def\Lhom#1#2#3{{{}\sp{#1}{\rm Hom}(#2,#3)}} 
\def\khom#1#2{{{\rm Hom}(#1,#2)}} 

\newcommand{\cC}{\mathcal{C}}
\newcommand{\can}{\operatorname{\it can}}
\newcommand{\alg}{\operatorname{Alg}}
\newcommand{\id}{\operatorname{id}}
\newcommand{\Me}{\mathcal{M}_A^C(\psi)}
\newcommand{\Mp}{\mathcal{M}_P^C(\psi)}
\newcommand{\Es}{(A,C, {\psi})}
\newcommand{\Ep}{(P,C,\psi)}
\newcommand{\Op}{\Omega^1 P}
\newcommand{\Ob}{\Omega^1 B}
\newcommand{\M}{\mathcal{M}}
\newcommand{\Ker}{\operatorname{Ker}}
\newcommand{\Ima}{\operatorname{Im}}
\newcommand{\Hom}{\operatorname{Hom}}
\newcommand{\tr}{\operatorname{Tr}}
\def\sw#1{{\sb{(#1)}}}
\def\sco#1{{\sp{(#1)}}} 
\def\su#1{{\sp{[#1]}}} 
\def\eps{\varepsilon}
\def\DC{\Delta_{\cC}}
\def\eC{\eps_\cC}
\def\DH{\Delta_{\cH}}
\def\eH{\eps_\cH}
\def\ut{\otimes}
\def\ov{\overline}
\def\act{\triangleleft}
\def\dl{{}_P\Delta}
\def\csop{\mbox{$(C^*)^{op}$}}
\def\lfa{{\large\mbox{$\forall\;$}}}
\def\fa{\forall\,\,}
\def\Rhom#1#2#3{{{\rm Hom}\sp{#1}(#2,#3)}}
\def\ten#1{\underset{#1}{\otimes}}
\def\coten#1{\underset{#1}{\square}}


\renewcommand{\labelenumi}{\theenumi)}
\newcommand{\nc}[2]{\newcommand{#1}{#2}}
\newcommand{\rnc}[2]{\renewcommand{#1}{#2}}

\rnc{\theequation}{\thesubsection.\arabic{equation}}

\def\note#1{}
\def\wegdamit#1{}
\def\CM{{\cal M}{\rm od}}
\def\st{such that}
\def\eps{{\epsilon}}
\def\vt{\triangleright}
\nc{\rp}{\mbox{$\Delta_{P}$}}
\nc{\beq}{\begin{equation}}
\nc{\eeq}{\end{equation}}
\rnc{\[}{\beq}
\rnc{\]}{\eeq}
\nc{\qpb}{quantum principal bundle}
\nc{\ff}{faithfully flat}
\nc{\ga}{\mbox{$\mathop{\mbox{\rm HGA}}^H_B(P)$}}
\nc{\red}{\mbox{$\mathop{\mbox{\rm Red}}^{H/J}_B(P)$}}
\nc{\bpr}{\begin{proposition}}
\nc{\bth}{\begin{theorem}}
\nc{\ble}{\begin{lemma}}
\nc{\bco}{\begin{corollary}}
\nc{\bre}{\begin{remark}}
\nc{\bex}{\begin{example}}
\nc{\bde}{\begin{definition}}
\nc{\ede}{\end{definition}}
\nc{\epr}{\end{proposition}}
\nc{\ethe}{\end{theorem}}
\nc{\phe}{principal homogenous extension}
\nc{\ele}{\end{lemma}}
\nc{\eco}{\end{corollary}}
\nc{\ere}{\hfill\mbox{$\Diamond$}\end{remark}}
\nc{\eex}{\hfill\mbox{$\Diamond$}\end{example}}
\nc{\epf}{\endproof}
\nc{\ot}{\otimes}
\nc{\bsb}{\begin{Sb}}
\nc{\esb}{\end{Sb}}
\nc{\ct}{\mbox{${\cal T}$}}
\nc{\ctb}{\mbox{${\cal T}\sb B$}}
\nc{\ba}{\begin{array}}
\nc{\ea}{\end{array}}
\nc{\bea}{\begin{eqnarray}}
\nc{\beas}{\begin{eqnarray*}}
\nc{\eeas}{\end{eqnarray*}}
\nc{\eea}{\end{eqnarray}}
\nc{\be}{\begin{enumerate}}
\nc{\ee}{\end{enumerate}}
\nc{\bcd}{\beq\begin{CD}}
\nc{\ecd}{\end{CD}\eeq}
\nc{\bi}{\begin{itemize}}
\nc{\ei}{\end{itemize}}
\nc{\kr}{\mbox{Ker}}
\nc{\te}{\!\ot\!}
\nc{\pf}{\mbox{$P\!\sb F$}}
\nc{\pn}{\mbox{$P\!\sb\nu$}}
\nc{\bmlp}{\mbox{\boldmath$\left(\right.$}}
\nc{\bmrp}{\mbox{\boldmath$\left.\right)$}}
\rnc{\phi}{\varphi}
\nc{\LAblp}{\mbox{\LARGE\boldmath$($}}
\nc{\LAbrp}{\mbox{\LARGE\boldmath$)$}}
\nc{\Lblp}{\mbox{\Large\boldmath$($}}
\nc{\Lbrp}{\mbox{\Large\boldmath$)$}}
\nc{\lblp}{\mbox{\large\boldmath$($}}
\nc{\lbrp}{\mbox{\large\boldmath$)$}}
\nc{\blp}{\mbox{\boldmath$($}}
\nc{\brp}{\mbox{\boldmath$)$}}
\nc{\LAlp}{\mbox{\LARGE $($}}
\nc{\LArp}{\mbox{\LARGE $)$}}
\nc{\Llp}{\mbox{\Large $($}}
\nc{\Lrp}{\mbox{\Large $)$}}
\nc{\llp}{\mbox{\large $($}}
\nc{\lrp}{\mbox{\large $)$}}
\nc{\lbc}{\mbox{\Large\boldmath$,$}}
\nc{\lc}{\mbox{\Large$,$}}
\nc{\Lall}{\mbox{\Large$\forall\;$}}
\nc{\bc}{\mbox{\boldmath$,$}}
\rnc{\epsilon}{\varepsilon}
\nc{\ra}{\rightarrow}
\nc{\ci}{\circ}
\nc{\cc}{\!\ci\!}
\nc{\T}{\mbox{\sf T}}
\nc{\cnl}{$\mbox{\sf T}\!\sb R$}
\nc{\lra}{\longrightarrow}
\rnc{\to}{\mapsto}
\rnc{\breve}{\widetilde}
\nc{\imp}{\Rightarrow}
\rnc{\iff}{\Leftrightarrow}
\nc{\ob}{\mbox{$\Omega\sp{1}\! (\! B)$}}
\nc{\op}{\mbox{$\Omega\sp{1}\! (\! P)$}}
\nc{\oa}{\mbox{$\Omega\sp{1}\! (\! A)$}}
\nc{\inc}{\mbox{$\,\subseteq\;$}}
\rnc{\subset}{\inc}
\nc{\spp}{\mbox{${\cal S}{\cal P}(P)$}}
\nc{\dr}{\mbox{$\Delta_{P}$}}
\nc{\dsr}{\mbox{$\Delta_{\Omega^1P}$}}
\nc{\ad}{\mbox{$\mathop{\mbox{\rm Ad}}_R$}}
\nc{\m}{\mbox{m}}
\nc{\hsp}{\hspace*}
\nc{\nin}{\mbox{$n\in\{ 0\}\!\cup\!{\Bbb N}$}}
\nc{\al}{\mbox{$\alpha$}}
\nc{\ha}{\mbox{$\alpha$}}
\nc{\hb}{\mbox{$\beta$}}
\nc{\hg}{\mbox{$\gamma$}}
\nc{\hd}{\mbox{$\delta$}}
\nc{\he}{\mbox{$\varepsilon$}}
\nc{\hz}{\mbox{$\zeta$}}
\nc{\hs}{\mbox{$\sigma$}}
\nc{\hk}{\mbox{$\kappa$}}
\nc{\hm}{\mbox{$\mu$}}
\nc{\hn}{\mbox{$\nu$}}
\nc{\la}{\mbox{$\lambda$}}
\nc{\hl}{\mbox{$\lambda$}}
\nc{\hG}{\mbox{$\Gamma$}}
\nc{\hD}{\mbox{$\Delta$}}
\nc{\Th}{\mbox{$\Theta$}}
\nc{\ho}{\mbox{$\omega$}}
\nc{\hO}{\mbox{$\Omega$}}
\nc{\hp}{\mbox{$\pi$}}
\nc{\hP}{\mbox{$\Pi$}}
\nc{\bpf}{\proof}
\nc{\as}{\mbox{$A(S^3\sb s)$}}
\nc{\bs}{\mbox{$A(S^2\sb s)$}}
\nc{\slq}{\mbox{$A(SL\sb q(2))$}}
\nc{\fr}{\mbox{$Fr\llp A(SL(2,{\Bbb C}))\lrp$}}
\nc{\slc}{\mbox{$A(SL(2,{\Bbb C}))$}}
\nc{\af}{\mbox{$A(F)$}}
\nc{\suq}{\mbox{$A(SU_q(2))$}}
\nc{\asq}{\mbox{$A(S_q^2)$}}
\nc{\tasq}{\mbox{$\widetilde{A}(S_q^2)$}}
\nc{\tc}{\widetilde{c}}
\rnc{\tilde}{\widetilde}

\def\esl{{\mbox{$E\sb{\frak s\frak l (2,{\Bbb C})}$}}}
\def\esu{{\mbox{$E\sb{\frak s\frak u(2)}$}}}
\def\zf{{\mbox{${\Bbb Z}\sb 4$}}}
\def\zt{{\mbox{$2{\Bbb Z}\sb 2$}}}
\def\ox{{\mbox{$\Omega\sp 1\sb{\frak M}X$}}}
\def\oxh{{\mbox{$\Omega\sp 1\sb{\frak M-hor}X$}}}
\def\oxs{{\mbox{$\Omega\sp 1\sb{\frak M-shor}X$}}}
\def\Fr{\mbox{Fr}}
\def\gal{-Galois extension}
\def\hge{Hopf-Galois $H$-extension}
\def\cge{coalgebra-Galois $C$-extension}
\def\cg{coalgebra-Galois}
\def\hgal{Hopf-Galois}
\def\age{algebra-Galois $A$-extension}
\def\agc{algebra-Galois $A$-coextension}
\def\lco{\!\!\!\phantom{I}^{co H/J}\!H}
\def\ta{\tilde a}
\def\tb{\tilde b}
\def\td{\tilde d}
\def\st{\stackrel}
\def\<{\langle}
\def\>{\rangle}
\def\d{\mbox{$\mathop{\mbox{\rm d}}$}}
\def\id{\mbox{$\mathop{\mbox{\rm id}}$}}
\def\im{\mbox{$\mathop{\mbox{\rm Im$\,$}}$}}
\def\ker{\mbox{$\mathop{\mbox{\rm Ker$\,$}}$}}
\def\coker{\mbox{$\mathop{\mbox{\rm Coker$\,$}}$}}
\def\hom{\operatorname{Hom}}
\def\ind{\operatorname{ind}}
\def\hor{\operatorname{hor}}
\def\endo{\operatorname{End}}
\def\map{\operatorname{Map}}
\def\o{\sp{[1]}}
\def\t{\sp{[2]}}
\def\mo{\sp{[-1]}}
\def\z{\sp{[0]}}
\def\te{\otimes}
\def\tens{\otimes}
\def\de{\hD}
\def\ses{short exact sequence}
\def\csa{$C^*$-algebra}
\def\lch{locally compact Hausdorff}
\def\lc{locally compact}
\def\cO{{\cal O}}
\newcommand{\Boxneu}{\square}
\def\C{{\Bbb C}}
\def\N{{\Bbb N}}
\def\R{{\Bbb R}}
\def\Z{{\Bbb Z}}
\def\T{{\Bbb T}}
\def\Q{{\Bbb Q}}
\def\cT{{\cal T}}
\def\cK{{\cal K}}
\def\cH{{\cal H}}
\def\cF{{\cal F}}
\def\im{{\rm Im}}
\def\id{{\rm id}}
\def\tr{{\rm tr}}
\def\Tr{{\rm Tr}}

\newcounter{zlist}
\newenvironment{zlist}{\begin{list}{(\arabic{zlist})}{
\usecounter{zlist}\leftmargin2.5em\labelwidth2em\labelsep0.5em
\topsep0.6ex
\parsep0.3ex plus0.2ex minus0.1ex}}{\end{list}}

\newcounter{blist}
\newenvironment{blist}{\begin{list}{(\alph{blist})}{
\usecounter{blist}\leftmargin2.5em\labelwidth2em\labelsep0.5em
\topsep0.6ex
\parsep0.3ex plus0.2ex minus0.1ex}}{\end{list}}

\newcounter{rlist}
\newenvironment{rlist}{\begin{list}{(\roman{rlist})}{
\usecounter{rlist}\leftmargin2.5em\labelwidth2em\labelsep0.5em
\topsep0.6ex
\parsep0.3ex plus0.2ex minus0.1ex}}{\end{list}}


\title{\Large\bf \vspace*{-15mm}
NONCOMMUTATIVE GEOMETRY APPROACH TO PRINCIPAL AND ASSOCIATED BUNDLES\\ 
\vspace*{5mm}}

\author{
\vspace*{-1mm}\large\sc
Paul F.~Baum   \\
\vspace*{-1mm}\normalsize
Mathematics Department,
McAllister Building  \\
\vspace*{-1mm}\normalsize
The Pennsylvania State University,
University Park, PA  16802, USA   \\
\vspace*{-0mm}\normalsize\sl e-mail: baum@math.psu.edu  \\
\vspace*{0mm}\and
\vspace*{-1mm}\large\sc
Piotr M.~Hajac\\
\vspace*{-1mm}\normalsize
Instytut Matematyczny, Polska Akademia Nauk\\
\vspace*{-1mm}\normalsize
ul.\ \'Sniadeckich 8, Warszawa, 00-956 Poland\\
\vspace*{-1mm}\normalsize
and\\
\vspace*{-1mm}\normalsize
Katedra Metod Matematycznych Fizyki, Uniwersytet Warszawski\\
\vspace*{-1mm}\normalsize
ul.\ Ho\.za 74, Warszawa, 00-682 Poland 
\vspace{0mm}\\
\normalsize\sl
http://www.impan.gov.pl/$\!\widetilde{\phantom{m}}\!$pmh\\
\vspace*{0mm}\and
\vspace*{-1mm}\large\sc
Rainer Matthes\\
\vspace*{-1mm}\normalsize
Katedra Metod Matematycznych Fizyki, Uniwersytet Warszawski\\
\vspace*{-1mm}\normalsize
ul.\ Ho\.za 74, Warszawa, 00-682 Poland 
\vspace{-0mm}\\
\normalsize\sl
e-mail: matthes@fuw.edu.pl\\
\vspace*{0mm}\and
\large\sc\vspace*{-1mm}
Wojciech Szyma\'nski\\
\vspace*{-1mm}\normalsize
School of Mathematical and Physical Sciences, University of Newcastle\\
\vspace*{-1mm}\normalsize
Callaghan, NSW 2308, Australia\\
\normalsize\sl
e-mail: Wojciech.Szymanski@newcastle.edu.au
}
\date{\normalsize }

\maketitle
\vspace*{-0mm}
\footnotesize\noindent
{\bf Abstract:}
We recast basic topological concepts underlying differential
geometry using the language and tools of noncommutative geometry.
This way we characterize principal (free and proper)
actions by a density condition
in (multiplier) $C^*$-algebras. We introduce the concept of piecewise
triviality to adapt the standard notion of local triviality to
fibre products of $C^*$-algebras. In the context of principal actions,
we study in detail an example 
of a non-proper free action with  continuous translation map,
and  examples of  compact principal bundles which are piecewise
trivial but not locally 
trivial,
 and neither piecewise trivial nor locally trivial, respectively. We show that the module of continuous sections of a
vector bundle associated to a compact principal bundle is a
 cotensor product of the algebra of functions defined on the total space 
(that are continuous along the base and polynomial along the fibres)
with the vector space of the representation. On the algebraic side,
we  
review the formalism of connections for the universal differential
algebras. 
\note{
Then, in an analytic vein,
we prove that even Fredholm modules can be considered
as appropriate pairs of bounded $*$-representations.
}
 In the differential
geometry framework,
we consider smooth connections on principal bundles as equivariant splittings of the
cotangent bundle, as 1-form-valued derivations of the algebra of smooth
functions on the structure group, and as axiomatically given covariant
differentiations of functions defined on the total space. Finally, we use the Dirac
monopole connection to compute the  pairing of the line bundles
associated to the Hopf fibration with the cyclic cocycle of integration
over~$S^2$. 

\small
\tableofcontents

\section*{Introduction}

This paper  is   mostly a survey. New results where derived to bridge
 gaps between the standard and noncommutative geometry approach
to principal and associated bundles. It is intended to serve as an introduction to the noncommutative geometry of
principal extensions of algebras. Since they play the role of quantum principal
bundles, herein we analyse in detail the topology of classical principal actions.
Following Henri Cartan, we define them as actions that are free and proper.
Next, we adopt the point of view that principal bundles are monoidal functors
from the category of finite-dimensional group representations to the category
of modules that are finitely generated projective over the 
algebra of functions defined on the base space. The functor is given by associating a vector bundle
to a principal bundle via a group representation, and then 
 taking the module of 
global sections of  the associated vector
bundle. This leads us to an algebraic formulation of the module
of continuous or smooth sections of a vector bundle. We focus in this context
on compact groups and spaces, so as to take an advantage of Peter-Weyl theory
and the Serre-Swan Theorem,
and thus arrive at the aforementioned general algebraic formalism without losing track of the
topological or smooth nature of vector bundles. 

Our exposition of connections on principal and vector bundles is done in a way that provides a passage
to the noncommutative setting. By standard Chern-Weil theory,
 connections are then used to define
the Chern character. As an example of this general theory, we study the
Dirac monopole connection on the Hopf fibration, and
prove that the integrated Chern classes of the associated line bundles
coincide with minus 
the winding numbers of representations that define them.
We interpret these topological invariants as the result of the pairing of 
$K$-theory and cyclic cohomology~\cite{c-a85}.

To make the paper reasonably self-contained, we enclose selected 
elements of general
noncommutative
algebra. In the special commutative case, they form an algebraic 
backbone of the
differential-geometric section. On the other hand, we refrain from 
discussing
the metric aspects of differential geometry. These are embodied in the 
theory
of spectral triples on the noncommutative side, and that is beyond the 
scope of this
paper.

\setcounter{equation}{0}
\section{Topological aspects of principal and  associated\\ bundles}
\label{topo}

The  category of
\lch\ spaces with all continuous maps as morphisms 
is equivalent to the (appropriately defined)
opposite category of commutative \csa s \cite{w-sl79}.
The goal of this section is to describe principal (free and
proper) actions of locally compact groups on locally compact
Hausdorff spaces purely in terms of their $C^*$-algebras. This
description will be the starting point for defining principal
actions of locally compact quantum groups on general
$C^*$-algebras. Throughout the paper, by a (locally) compact group
we mean a topological group which, as a set, is (locally) compact
 Hausdorff. 

A key technical concept appearing while considering principal actions and fibre bundles
is that of a fibred product. For any three sets $X,Y,Z$, and any two maps $f:X\ra Z$,
$g:Y\ra Z$,
we can define the {\em fibred product} of $X$ and $Y$ over $Z$ as the equaliser of
the maps $F(x,y)=(x,f(x),y)$, $G(x,y)=(x,g(y),y)$, i.e.
\[\label{fibred}
\xymatrix@+1pc{
*!<0ex, 0.8ex>{\{(x,y)\in X\times Y\;|\;f(x)=g(y)\}=:X\underset{Z}{\times}Y} \ar[r]^-{\subseteq} 
& X\times Y \ar@<1ex>[r]^-{F} \ar@<-1ex>[r]_-{G} & X\times Z\times Y\\
}
\]
In the noncommutative setting, the fibred product will appear reincarnated as
the cotensor product over a coalgebra
and the tensor product over an algebra.

Recall that a continuous map $F\in C(X,Y)$ is proper iff the
map
\[
F\times\id: X\times Z\lra Y\times Z
\]
is closed for
any topological space Z (see \cite[Definition 1, p. 97]{b-n71}).  
A continuous group action $X\times G\ra
X$ is called {\em proper} iff the (principal) map
\[\Label{pmap}
F^G:X\times G\ni (x,g)\longmapsto (x,xg)\in X\times X
\]
is proper. Now, if $X$ is Hausdorff and $Y$ is \lch, then a
continuous map $F\in C(X,Y)$ is proper if and only if
  $F^{-1}(K)$ is compact for any compact subset~$K$
  \cite[Proposition~7, p.104]{b-n71}.
  In particular, if $G$ is a locally compact group acting on a
  \lch\ space $X$, we can say that this action is proper iff
  $(F^G)^{-1}(K)$ is compact for any compact subset~$K$. (By a
\lc\ group we mean a topological group which is \lc\ and
Hausdorff.) A beautiful discussion and comparison study of
different definitions of proper actions can be found
in~\cite{b-h04}.

 An action is called {\em free} iff  $F^G$ is
injective, and {\em principal}
iff this map is both injective and proper. In other words, principal means free and
proper. It follows from Lemma~\ref{closed} that, if the action is  principal,
$F^G$ is a closed
injection, so that it yields a homeomorphism $F^G_X$
from $X\times G$ onto the fibred product
$X\times_{X/G}X$, which is the image of $F^G$. (Note here that $X/G$ is again
a \lch\ space \cite[Proposition 3, p. 253, Proposition 9, p. 257]{b-n71}.)
However, the opposite implication is not true: even when $F^G_X$
is a homeomorphism the action of $G$ on $X$ might still fail to be proper
(see Example~\ref{rysiu}).
Notice also that an action can be proper without being free (e.g., the trivial
action of a non-trivial compact group) and free without being proper (e.g.,
an ergodic action of $\R$ on a torus). Principal actions lead to the concept
of a topological principal bundle, and free but not proper actions give rise
to foliations.

 If the action is free, then  there exists a (not necessarily continuous) map
\[\label{classtrans}
\xymatrix@+1.5pc{
*!<0ex, 1.3ex>{\check{\tau}\::\: X\underset{X/G}{\times}X} \ar[r]^-{(F^G_X)^{-1}} 
& X\times G \ar[r]^-{\phantom{w}\pi_G\phantom{w}} & G
}
\]
where $\pi_G$ is the canonical surjection onto~$G$. Since it determines which element
of $G$ translates one point of $X$ to the other, it
is called the {\em translation map} \cite[Chapter~4, Definition 2.1]{h-d94}.
The following properties of $\check{\tau}$ are immediate from the definition:
\bea\label{tran}
\check{\tau}(x,xg)&=&g,\\
x\check{\tau}(x,y)&=&y,\\
\check{\tau}(xg,yh)&=&g^{-1}\check{\tau}(x,y)h,\label{transequiv}\\
\check{\tau}(x,y)\check{\tau}(y,z)&=&\check{\tau}(x,z),\\
\check{\tau}(y,x)&=&\check{\tau}(x,y)^{-1}.
\eea
Here $x,y,z\in X$, $g,h\in G$, and $e$ is the neutral element of~$G$. As will be
demonstrated later, the continuity of $\check{\tau}$ is guaranteed by the properness
of the action.
The noncommutative counterpart
of the translation map plays a fundamental role in the theory of Galois-type
extensions and coextensions.

There is some confusion in the literature concerning the definition of a principal action
(principal bundle). Some authors call principal an action that is free and induces
continuous translation map. While there are important results for which these two
properties suffice, they do not guarantee that the quotient space is Hausdorff
(see Example~\ref{rysiu}). This is why we adopt the more restrictive definition that
calls principal an action that is free and proper, and which implies that
the space of orbits is Hausdorff. As we shall argue later on, our definition is equivalent with
that of H.~Cartan, who requires not only the continuity of the translation map, but also
that $X\times_{X/G}X$ be closed in $X\times X$.

\subsection{Proper maps}

To begin with, let us recall the separation axioms of topology and other basic
definitions. Let $X$ be a topological space,  $O(X)$ its topology (the family
of open subsets) and $C(X,Y)$ the space of continuous functions from $X$ to
some topological space $Y$. We say that $X$ is
\begin{itemize}

\item[$T_{0\phantom{\frac{1}{2}}}$]
iff $\fa x,y\in X,\, x\not =y\;\exists\,\,U\in O(X):\;
(x\in U \mbox{ and } y\not\in U) \mbox{ or } (x\not\in U \mbox{ and } y\in U)$;

\item[$T_{1\phantom{\frac{1}{2}}}$]
 iff $\fa x,y\in X,\, x\not =y\;\exists\,\,U,V\in O(X):\;
(x\in U \mbox{ and } y\not\in U) \mbox{ and } (x\not\in V \mbox{ and } y\in V)$;

\item[$T_{2\phantom{\frac{1}{2}}}$]
  or Hausdorff iff $\fa x,y\in X,\, x\not =y\;\exists\,\,U,V\in O(X):\;
x\in U,\;  y\in V \mbox{ and } U\cap V=\emptyset$;

\item[$T_{3\phantom{\frac{1}{2}}}$]
  or regular iff it is $T_2$ and $\fa x\in X,\; \mbox{closed } C\inc X,\,
x\not\in C\;\exists\,\,U,V\in O(X):$\\
$x\in U,\;  C\inc V \mbox{ and } U\cap V=\emptyset$;

\item[$T_{3\frac{1}{2}}$]  or Tichonov or completely regular iff it is $T_2$ and\\
$\fa x\in X,\; \mbox{closed } C\inc X,\,
x\not\in C\;\exists\,\,f\in C(X,[0,1]):\; f(x)=0\mbox{ and } f(C)=\{1\}$;

\item[$T_{4\phantom{\frac{1}{2}}}$]
or normal iff it is $T_2$ and $\fa \mbox{closed } C,D\inc X,\,
C\cap D=\emptyset\;\exists\,\,U,V\in O(X):$\\
$C\inc U,\;  D\inc V \mbox{ and } U\cap V=\emptyset$.

\end{itemize}

A topological space is called {\em compact} if out of any of its
open covers one can choose a finite subcover, and {\em locally
compact} if any of its points is contained in an open set whose
closure is compact. (In some books, e.g. \cite{b-n71,e-r77},
(locally) compact means (locally) compact and Hausdorff.) Every
compact Hausdorff space is normal \cite[Theorem~3.1.9]{e-r77} and
every \lch\ space is completely regular
\cite[Theorem~3.3.1]{e-r77}, though not necessarily normal
\cite[Example~3.3.14]{e-r77}. Also, since every metrizable space is normal
\cite[Corollary~4.1.13]{e-r77}, the space of rational numbers is
normal although one can directly check that it is not locally
compact.
 Symbolically, we can write all this in the following
way:

\vspace*{0mm}\[
\xymatrix{
& \text{normal} \ar@{=>}[dl] \ar@{=}@<1ex>[d] & \\ 
\text{completely regular}  & \text{does not imply} \ar@{=>}@<1ex>[u] \ar@{=>}@<1ex>[d]
& \text{compact Hausdorff} \ar@{=>}[ul] \ar@{=>}[dl]\\
& \text{locally compact Hausdorff} \ar@{=>}[ul] \ar@{=}@<1ex>[u]&
}
\]\vspace*{0mm}

\note{
Since we want to describe
principal group actions in terms of \csa s, we focus our attention
to \lch\ spaces and continuous actions by \lc\ groups.
}
Having established the terminology, let us now review the relevant key facts.
\ble[Urysohn Lemma]\Label{ury}
Let $X$ be a normal topological space and $C,D$ disjoint closed subsets of $X$.
Then there exists a continuous function $f:X\ra [0,1]$ such that $f(C)\inc\{0\}$
and $f(D)\inc\{1\}$.
\ele
\begin{theorem}[Tietze-Urysohn Extension Theorem]\Label{tie}
A topological space $X$ is normal if and only if for any closed subset
$C\inc X$ any $\C$-valued bounded continuous function on $C$ can be extended
to a bounded continuous function on $X$ with the same sup norm.
\end{theorem}
\ble\Label{cel} (cf. \cite[Exercise~3.2.J]{e-r77})\phantom{.}
Let $Y$ be a \lch\ space and $L$ a compact subset of $Y$. Then any continuous
function $g:L\ra\C$ can be extended to a continuous function on $Y$ with compact
support and the same sup norm.
\ele
\bpf
\note{
Since $Y$ is \lch, it is a subset of its one-point compactification $Y^+=Y\cup\{\infty\}$.
(Recall that the open sets in $Y^+$ are the open sets in $Y$ and the sets of the form
$(X\setminus K)\cup\{\infty\}$, where $K$ is compact.)
 On the other hand, as $L$ is compact and $Y^+$ is Hausdorff,
it is closed in~$Y^+$. Consequently, since $Y^+$ is compact Hausdorff
and any compact Hausdorff space is normal,
by Theorem~\ref{tie}, any continuous bounded function on $L$ can be extended to
a continuous function on $Y^+$ with the same norm. Therefore, it can also be extended
to a continuous function on $Y$ with the same norm.
}
Since $L$ is compact and $Y$ is \lch , there exists an open set $U$ containing $L$
and whose closure $\overline{U}$ is compact \cite[Theorem~3.3.2]{e-r77}. Similarly,
there exists open $V$ such that $\overline{U}\inc V$ and $\ov{V}$ is compact.
Since $Y$ is Hausdorff, $\ov{V}$ and $L$ are closed. The set $L$ is also closed in
$\ov{V}$, as $\ov{V}\setminus L=\ov{V}\cap(Y\setminus L)$. The set $U$ is open
and $U=U\cap\ov{V}$, so that it is open in $\ov{V}$. Thus $\ov{V}\setminus U$
is closed in $\ov{V}$. Next, we have $L\cap(\ov{V}\setminus U)=\emptyset$, so that
we can define a continuous function $\widetilde{g}:L\cup(\ov{V}\setminus U)\ra\C$
by the formula
\[
\widetilde{g}(y)=\left\{\ba{cr}
g(y)&\mbox{for }y\in L,
\\&\\
0&\mbox{for }y\in \ov{V}\setminus U.
\\
\ea\right.
\]
(If the union of two closed subsets equals the space, then a function is
continuous if and only if its restrictions to the closed subets are continuous.)
Clearly, $\|\widetilde{g}\|=\|g\|$. Since $\ov{V}$ is a compact Hausdorff space,
it is normal. Thus, by Theorem~\ref{ury}, there exists a continuous extension
$\widetilde{\widetilde{g}}:\ov{V}\ra\C$  of $\widetilde{g}$ such that
$\|\widetilde{\widetilde{g}}\|=\|\widetilde{g}\|$. Now define
$g_Y:Y\ra\C$ by the formula
\[
g_Y(y)=\left\{\ba{cr}
\widetilde{\widetilde{g}}(y)&\mbox{for }y\in \ov{V},
\\&\\
0&\mbox{for }y\in Y\setminus U.
\\
\ea\right.
\]
This function is well defined because $(Y\setminus U)\cap\ov{V}=\ov{V}\setminus U$
and $\widetilde{\widetilde{g}}(Y\setminus U)\inc\{0\}$. It is continuous
because both $\ov{V}$ and  $Y\setminus U$ are closed and $\ov{V}\cup(Y\setminus U)
=Y$.  It extends $g$ and has compact support
by construction. Finally, its norm remains unchanged
($\|g_Y\|=\|g\|$).
\epf
A straightforward application of this lemma (or its simpler version), allows one
to prove a $C^*$-characterisation of properness. To this end, recall first that a function
$f:X\ra\C$ is said to be {\em vanishing at infinity} iff
\[
\fa\he>0\;\exists\,\text{compact}\, K\inc X\;\fa x\not\in K:\; |f(x)|<\he\,.
\]
The $C^*$-algebra of all continuous vanishing at infinity functions on a \lch\ space
is denoted by~$C_0(X)$. Now we can state:
\ble\Label{proper}
Let $X$ and $Y$ be \lch\ spaces. A continuous map $F:X\ra Y$ is {\em proper}
if and only if $F^*(C_0(Y))\inc C_0(X)$.
\ele
\bpf
Assume first that $F$ is proper. We want to show that $f\circ F$ vanishes at
infinity for any $f\in C_0(Y)$. Choose $\he>0$. There exists compact $K_\epsilon
\inc Y$ such that $y\not\in K_\epsilon\Rightarrow |f(y)|<\he$. Since $F$ is
proper, $F^{-1}(K_\epsilon)$ is compact. Moreover, we have
$x\not\in F^{-1}(K_\epsilon)\Leftrightarrow F(x)\not\in K_\epsilon
\Rightarrow |f(F(x))|<\he$. Hence $F^*(f)\in C_0(X)$, as needed.

Assume now that   $F^*(C_0(Y))\inc C_0(X)$.
 Since  $Y$ is Hausdorff,
any compact subset $K$ of $Y$ is closed. Therefore, so is $F^{-1}(K)$
 by the continuity of $F$. The case $K=\emptyset$ is trivial, so that,
without loss of generality, we may assume $K\not=\emptyset$. Then,
 by Lemma~\ref{cel},
there exists a compactly supported continuous function  $f$ on $Y$ such that
$f(K)=\{1\}$.
Due to our assumption $f\circ F$ vanishes at infinity. Hence, for $\he=\frac{1}{2}$,
there exists a compact set $L\inc X$ such that
$x\not\in L\Rightarrow |f(F(x))|<\frac{1}{2}$. Consequently, $F^{-1}(K)\inc L$,
and $F^{-1}(K)$ has to be compact as a closed subset of a compact set.
This proves that $F$ is proper.
\epf
For the sake of the general $C^*$-algebraic setting, it is useful to have the
following refinement of the above lemma:
\bco\Label{refproper}
Let $X$ and $Y$ be \lch\ spaces and $F:X\ra Y$ a continuous map.
Then $F$ is  proper
if and only if $F^*(A)\inc C_0(X)$ for some norm dense subset $A$ of $C_0(Y)$.
\eco
\bpf
Due to Lemma~\ref{proper}, it suffices to show that
$F^*(A)\inc C_0(X)$ for some norm dense subset $A$ of $C_0(Y)$ if and only if
$F^*(C_0(Y))\inc C_0(X)$.  To prove
the left-to-right implication, recall that for any set
$S$ and any continuous function
$f$, we have $f(\ov{S})\inc\ov{f(S)}$. (Indeed,
$S\inc f^{-1}(f(S))\inc f^{-1}(\ov{f(S)})$. As the last set is closed by the
continuity of $f$, we have $\ov{S}\inc f^{-1}(\ov{f(S)})$, whence
$f(\ov{S})\inc\ov{f(S)}$.)
Therefore, since $C_0(X)$ is norm closed and
$F^*$ is continuous in the norm topology, we have:
\[\Label{inc}
F^*(C_0(Y))=F^*(\ov{A})\inc \ov{F^*(A)}\inc\ov{C_0(X)}=C_0(X).
\]
The opposite implication is trivial.
\hfill\epf

\subsection{C*-characterisation of injectivity}

We can pass now to the most demanding part of our topological introduction,
 which a $C^*$-characterisation
of injectivity. Although such a result should have been proven decades ago, the
only reference we are aware of is \cite{e-da00}. The main difficulty is to
characterize injective but non-proper maps. (For instance, take
$
[0,1)\ni\theta\mapsto e^{2\pi i\theta}\in\C
$.)
To this end, we need the concept of strict topology on multiplier \csa s (see
\cite{lcqg}).  If
$A$ is a \csa\ and $M(A)$ is its multiplier $C^*$-algebra, then a subset $C\inc M(A)$
is {\em strictly dense} (dense in the strict topology) iff
$\fa x\in M(A)\;\exists\;\{c_\alpha\}_{\alpha\in\mbox{\scriptsize directed set}}
\inc C
\;\fa a\in A$:
\[\Label{strict}
\; \lim_\alpha\|(c_\alpha-x)a\|+\|a(c_\alpha-x)\|=0.
\]
In the commutative case, we have $M(C_0(X))=C_b(X)$ (the latter means the
algebra of continuous bounded functions on $X$), and our $C^*$ characterisation
of injectivity takes the following form:
\begin{theorem}[\cite{e-da00}]\Label{eda}
Let $X$ and $Y$ be \lch\ spaces. A continuous map $F:X\ra Y$ is {\em injective}
if and only if $F^*(C_0(Y))$ is strictly dense in $C_b(X)$.
\end{theorem}
\bpf
Suppose first that $F$ is not injective. Let $x\neq y$ be such that $F(x)=F(y)$.
Since $X$ is Hausdorff, there exists a closed set containing $y$ and not
containing $x$. Since $X$ is also locally compact, it is completely regular.
Therefore, there exists a  continuous function $f$ such that $f(x)\neq f(y)$.
On the other hand, as $\{x,y\}$ is compact, by Lemma~\ref{cel}, there exists
$g\in C_0(X)$ such that $g(x)=1=g(y)$. If $F^*(C_0(Y))$ is strictly dense in
$C_b(X)$, then $\lim_\alpha\|c_\alpha g-fg\|=0$ for some generalised sequence
$\{c_\alpha\}_{\alpha\in I}\inc F^*(C_0(Y))$. (Note that $c_\alpha(x)=c_\alpha(y)$.)
Thus, to derive the desired
contradiction, we can use the triangle inequality:
\begin{align}
\|c_\alpha g-fg\|&\;\geq\;\max\{|(c_\alpha g)(x)-(fg)(x)|,|(c_\alpha g)(y)-(fg)(y)|\}
\nonumber\\
&\;\geq\;\frac{1}{2}\left(|(c_\alpha g)(y)-(fg)(x)|+|(c_\alpha g)(y)-(fg)(y)|\right)
\nonumber\\
&\;\geq\;\frac{|(fg)(x)-(fg)(y)|}{2}>0.
\end{align}
(We simply exploited the fact that the evaluation map is continuous in the strict topology.)
Assume now that $F$ is injective. Let $f\in C_b(X)$. By Lemma~\ref{cel}, to each
compact subset $K$ of $X$ we can assign $c_K\in C_c(X)$ (compactly supported
continuous function) such that $c_K(K)\inc\{1\}$. Define $f_K=fc_K$ and assume
$K\neq\emptyset$.
Since $K$ is compact, $F$ is injective and $Y$ is Hausdorff, $F|_K:K\ra F(K)$ is
a homeomorphism (e.g., see \cite[p.87]{b-n71}). Thus $f_K\circ(F|_K)^{-1}$ is
continuous on the compact subset $F(K)$ of a \lch\ space $Y$. Hence,
by Lemma~\ref{cel}, it can be extended to $g_K\in C_c(Y)$ such that
$\|g_K\|=\|f_K\circ(F|_K)^{-1}\|$. We want to show that
$\lim_{\emptyset\neq K\subseteq X}F^*(g_K)=f$ (in the strict topology).
This means that
\[
\lfa a\in C_0(X): \lim_{\emptyset\neq K\subseteq X}\|(F^*(g_K)-f)a\|=0.
\]
If $X$ is compact, we are done. Otherwise, choose $a\in C_0(X)$ and $\he>0$.
For any non-empty compact subset $K$ of $X$, we have
\bea
\|F^*(g_K)-f\| &\leq& \|g_K\|+\|f\|
\nonumber\\ &\leq&
\|f_K\|+\|f\|
\nonumber\\ &\leq&
\|f\|~\|c_K\|+\|f\|
\nonumber\\ &=&
2\|f\|.
\eea
On the other hand, as $\llp(F^*(g_K)-f)a\lrp(K)=\{0\}$, we have
\bea
\|(F^*(g_K)-f)a\| &=& \sup_{x\not\in K}|F^*(g_K)-f|~|a|
\nonumber\\ &\leq&
\sup_{x\not\in K}|F^*(g_K)-f|\sup_{x\not\in K}|a|
\nonumber\\ &\leq&
\|F^*(g_K)-f\|\sup_{x\not\in K}|a|
\nonumber\\ &\leq&
2\|f\|\sup_{x\not\in K}|a|.
\eea
Next, since $a$ is a function vanishing at infinity, there exists a non-empty
compact set $K$ such that $2\|f\|\sup_{x\not\in K}|a|<\he$. Therefore,
\[
\lfa a\in C_0(X),\, f\in C_b(X):\;
\lim_{\emptyset\neq K\subseteq X}\|(F^*(g_K)-f)a\|=0,
\]
i.e., $F^*(C_0(Y))$ is strictly dense in $C_b(X)$.
\epf
Again, we can refine the assertion of this theorem to suit our later purpose.
To this end, let us first prove the following general lemma:
\ble\Label{gl}
Let $\phi:S\ra T$ be continuous and let $A\inc B\inc\ov{A}\inc S$. Assume also
that $T$ is equipped with an additional coarser topology
(i.e., the closure in the coarser topology is bigger: 
$\ov{C}\inc\ov{C}^c\inc T$). Then $\ov{\phi(A)}^c=\ov{\phi(B)}^c$.
\ele
\bpf
Since $\phi$ is continuous, we have $\phi(\ov{A})\inc\ov{\phi(A)}$. Hence
\[
\phi(A)\inc\phi(B)\inc\phi(\ov{A})\inc\ov{\phi(A)}\inc\ov{\phi(A)}^c.
\]
Consequently, $\ov{\phi(A)}^c\inc\ov{\phi(B)}^c\inc\ov{\phi(A)}^c$.
\epf
The strict topology is coarser than the norm topology. Therefore,
taking $B=\ov{A}=C_0(Y)$
 in the preceding lemma and combining it with Theorem~\ref{eda},
yields:
\bco
Let $X$ and $Y$ be \lch\ spaces and $F:X\ra Y$ a continuous map. If $F$
is  injective,
then $F^*(A)$ is strictly dense in $C_b(X)$ for any norm dense subset
$A$ of $C_0(Y)$. If $F^*(A)$ is strictly dense in $C_b(X)$ for some norm dense subset
$A$ of $C_0(Y)$, then $F$
is  injective.
\eco

\subsection{Free and proper actions}

We have just derived the $C^*$-characterisations of injectivity and propernesss
as independent conditions. Applying it to the map  (\ref{pmap}), we obtain
a $C^*$-definition of free actions and proper actions. Since we are interested
here in actions that are both free and proper (principal), we would like to know what is
a good $C^*$-characterisation of injectivity and properness treated simultaneously,
as one condition. In other words, the question is what  strict density of
$\im F^*$ in $C_b(X)$ (injectivity of $F:X\ra Y$) and the inclusion
$\im F^*\inc C_0(X)$ (properness of $F$) add up to. The answer is: the
 surjectivity of $F^*$. To prove this assertion by direct topological arguments,
let us first recall two more lemmas of general topology.
\ble \Label{closed}
Let $X$ and $Y$ be Hausdorff spaces, and let $Y$ be locally compact.
A continuous injection $F:X\ra Y$ is  proper
if and only it is closed.
\ele
\bpf
Assume first that $F$ is closed. Then $F(X)$ is closed in $Y$.
 Since $Y$ is Hausdorff,  any compact $K\inc Y$ is closed. Consequently,
$F(X)\cap K$ is closed in $Y$ and in $K$. Moreover, it is compact as a closed
subset of a compact set. On the other hand, since $F$ is injective and closed,
the map $F^{-1}:F(X)\ra X$ is continuous.
Hence the set $F^{-1}(F(X)\cap K)=F^{-1}(K)$ is also compact. This shows that
 $F$ is proper.

Suppose now that $F$ is proper but not closed. Then there exists a closed
$C\inc X$ such that $F(C)$ is not closed. This means that there exists
$y\in \ov{F(C)}\setminus F(C)$. Since $Y$ is locally compact, there exists
a neighbourhood $U$ of $y$ such that $\ov{U}$ is compact. The set $F^{-1}(\ov{U})$
is also compact due to the properness of $F$. As $X$ is Hausdorff,
$C\cap F^{-1}(\ov{U})$ is again compact by the same argument as in the first part of
the proof. Therefore, $F(C\cap F^{-1}(\ov{U}))$ is compact, whence closed by the
Hausdorffness of $Y$. Clearly, $F(C\cap F^{-1}(\ov{U}))\inc F(C)$, so that
\[\Label{contr}
y\not\in F(C\cap F^{-1}(\ov{U})).
\]
 We want to show that
$y\in\ov{F(C\cap F^{-1}(\ov{U}))}=F(C\cap F^{-1}(\ov{U}))$ and thus derive the
desired contradiction. To this end, recall that
\[
x\in\ov{A}\;\Leftrightarrow\;
\fa \mbox{open }U\ni x:\, A\cap U\neq\emptyset.
\]
 (Indeed,
$x\not\in\ov{A}\Leftrightarrow\,\,\exists\,\,\mbox{closed }C\not\ni x,\, A\inc C
\Leftrightarrow\,\,\exists\,\,\mbox{closed }C,\, x\in X\setminus C,
\,(X\setminus C)\cap A=\emptyset
\Leftrightarrow\,\,\exists\,\,\mbox{open } U\ni x,\, U\cap A=\emptyset$.)
Now, take any neighbourhood $V$ of $y$. Since $U\cap V$ is also a neighbourhood
of $y\in\ov{F(C)}$, we have  $U\cap V\cap F(C)\neq\emptyset$. Hence there exists
$c\in C$ such that $F(c)\in U\cap V\inc\ov{U}$. Thus $c\in C\cap F^{-1}(\ov{U})$,
whence $F(c)\in F(C\cap F^{-1}(\ov{U}))$. On the other hand, as $F(c)\in V$, we have
$F(C\cap F^{-1}(\ov{U}))\cap V\neq\emptyset$. However, this means that
$y\in\ov{F(C\cap F^{-1}(\ov{U}))}=F(C\cap F^{-1}(\ov{U}))$, which contradicts
(\ref{contr}).
\epf
\bre
In the second part of the proof we have not used the assumption that $F$ is
injective. Therefore, we can claim that any proper map from a Hausdorff space
to a \lch\ space is closed.
\ere
\ble \Label{ext}
Let $Y$ be a \lch\ space and $X$ a closed subset of $Y$. Then any $f\in C_0(X)$
can be extended to an element of $C_0(Y)$.
\ele
\bpf
Note first that any closed subset of a \lch\ space is itself a \lch\ space.
Indeed, let $x\in X$ and $U$ be a neighbourhood of $x$ in $Y$ such that $\ov{U}$
is compact. Then $X\cap U$ is a neighbourhood of $x$ in $X$ and
\[\Label{cap}
\ov{X\cap U}=\bigcap_{C=\ov{C},\,X\cap U\subseteq C\cap X}\!\!\!\!\!\!C\cap X\;\;\inc
\bigcap_{C=\ov{C},\,U\subseteq C}\!\!\!C\cap X=
\left(\bigcap_{C=\ov{C},\,U\subseteq C}\!\!\!C\right)
 \cap X=\ov{U}\cap X.
\]
Since $X$ is closed, $\ov{U}\cap X$ is compact. Hence it follows from (\ref{cap})
that $\ov{X\cap U}$ is compact. Consequently, $X$ is a \lch\ space. Let $X^+=X\cup
\{\infty\}$
and $Y^+=Y\cup
\{\infty\}$ denote the one-point compactifications of $X$ and $Y$, respectively.
(Recall that the open sets in $X^+$ are the open sets in $X$ and the sets of the form
$(X\setminus K)\cup\{\infty\}$, where $K$ is compact.)
Let us now show that any $f\in C_0(X)$ extends to $f^+\in C(X^+)$ by the formula
$f^+(\infty)=0$. If $D\inc\C$ is an open disc not containing 0, then $\infty\not\in
(f^+)^{-1}(D)$, and we have $(f^+)^{-1}(D)=f^{-1}(D)$. By the continuity of $f$,
the set $f^{-1}(D)$ is open in $X$, and consequently in $X^+$. If $D_\epsilon\inc\C$
is an open disc of radius $\epsilon$ centered at 0, then, as $f$ is a function
vanishing at infinity, there exists a compact subset $K_\epsilon$ of $X$ such that
$x\not\in K_\epsilon\Rightarrow f(x)\in D_\epsilon$. Hence the  set
$(f^+)^{-1}(D_\epsilon)=
\{\infty\}\cup(X\setminus K_\epsilon)\cup(f|_{K_\epsilon})^{-1}(D_\epsilon)$
is open in $X^+$ due to the continuity of $f$. Since any open set in $\C$ can
be decomposed into the union of open discs belonging to the above discussed family,
we can conclude that $f^+$ is continuous on $X^+$, as claimed. Next, since
$Y^+\setminus X^+=Y\setminus X$ is open in $Y$, it is open in $Y^+$. Hence $X^+$
is a closed subset of the compact Hausdorff space $Y^+$. As the latter is normal,
we can apply Theorem~\ref{tie} to conclude that $f^+$ can be extended to $f^+_{Y^+}
\in C(Y^+)$. We need to prove now that
the restriction of  $f^+_{Y^+}$ to $Y$ is an element of $C_0(Y)$.
First, note that $f^+_{Y^+}(\infty)=f^+(\infty)=0$. Moreover, as $f^+_{Y^+}$ is
continuous at $\infty$, we have
\[
\fa\he>0\,\,\exists\,\,\mbox{compact } L_\epsilon:\,
y\in (Y\setminus L_\epsilon)\cup\{\infty\}\Rightarrow |f^+_{Y^+}(y)|<\he.
\]
Since $L_\epsilon$ is compact by construction, we have that $f^+_{Y^+}|_Y$
vanishes at infinity, as needed. Finally, $f^+_{Y^+}|_Y$ is an extension of $f$
because $f^+_{Y^+}|_Y(x)=f^+_{Y^+}(x)=f^+(x)=f(x)$.
\epf
Combining Lemma~\ref{closed} and Lemma~\ref{ext} we obtain:
\begin{theorem}\Label{}
Let $X$ and $Y$ be \lch\ spaces. A continuous map $F:X\ra Y$ is {\em injective
and proper}
if and only if $F^*(C_0(Y))=C_0(X)$.
\end{theorem}
Reasoning along the lines of (\ref{inc}) proves the following refinement of
this theorem:
\bco\Label{injclo}
Let $X$ and $Y$ be \lch\ spaces and $F:X\ra Y$ a continuous map.
Then $F$ is injective and proper
if and only if $\ov{F^*(A)}= C_0(X)$ for some norm dense subset $A$ of $C_0(Y)$,
and if and only if $\ov{F^*(A)}= C_0(X)$ for any norm dense subset $A$ of $C_0(Y)$.
\eco

Thus we are able to define a principal action of $G$ on $X$ entirely in terms of
their $C^*$-algebras --- it suffices to take $F=F^G$. However, as was already
mentioned at the beginning of this section, for non-compact groups it is not
sufficient to consider $F^G_X$. Indeed, with the help of Lemma~\ref{closed}, one can check that
$F^G$ is injective and proper if and only if $F^G_X$ is a homeomorphism and its image
is closed.  On the other hand, $F^G_X$ is a homeomorphism if and only if  the translation
map (\ref{classtrans}) is continuous. Therefore, the action is principal if and only if
the translation map is continuous and $X\times_{X/G}X$ is closed in $X\times X$
\cite[Proposition 6, p.255]{b-n71}. This way we have arrived at the definition of
principal action provided by H.~Cartan \cite[condition (FP), p.6-05]{c-h49/50}.

Recall that a continuous injection
need not map  homeomorphically its domain to its image, e.g., take
$
[0,1)\ni\theta\mapsto e^{2\pi i\theta}\in\C
$.
The ergodic action of $\R$ on $T^2$ gives an example of     $F^G_X$ which is a continuous
bijection but not a homeomorphism (discontinuous translation map).
Similarly, a continuous
map that maps  homeomorphically its domain to its image need not be proper,
e.g., take
$
\R\ni x\mapsto \arctan(x)\in\R
$.
It is not proper because its image is not closed.
This is precisely why it is not sufficient to assume that $F^G_X$ is a homeomorphism.
Let us consider the following example in which $F^G_X$ is a homeomorphism
and the action is not principal. (This example
is essentially the same as in   \cite[p.298]{p-rs61}.)

\bex\Label{rysiu}
Take $X=\R^2$ and $G=\R$. The idea is to define the action of $\R$ on $\R^2$
by the flow of the (unital) smooth vector field on $\R^2$ given by the formula
$v(x,y)=(\cos x, \sin x)$. This vector field is invariant with respect to the
vertical translations. Its integral curves are vertical lines (special orbits)
for
$x=\frac{\pi}{2}+\mu \pi$, $\mu\in\Z$, and otherwise satisfy the equation
$\gamma'(x)=\tan(x)$ (regular orbits). Hence $\gamma(x)=\log|\!\cos x|^{-1} +a$,
$a\in\R$, or
$x=\frac{\pi}{2}+\mu \pi$. Roughly speaking, what happens is that regular orbits
converge to both of their neighbouring special orbits. Therefore, we can choose
a sequence of pairs of points from the same regular orbit that converges to a pair of
points belonging to two different orbits. This shows that the image of $F^G$ is not
closed. On the other hand, contrary to the Kronecker foliation case, if such a
sequence of pairs  converges to a pair of points on the same orbit, then the induced
sequence of group elements determined by each pair also converges to the group
element corresponding to the limit pair. This means that $F^G_X$ is a homeomorphism
despite $F^G$ not being proper. Let us explicitly prove all of this.

\begin{figure}[h]
\[
\includegraphics[width=70mm]{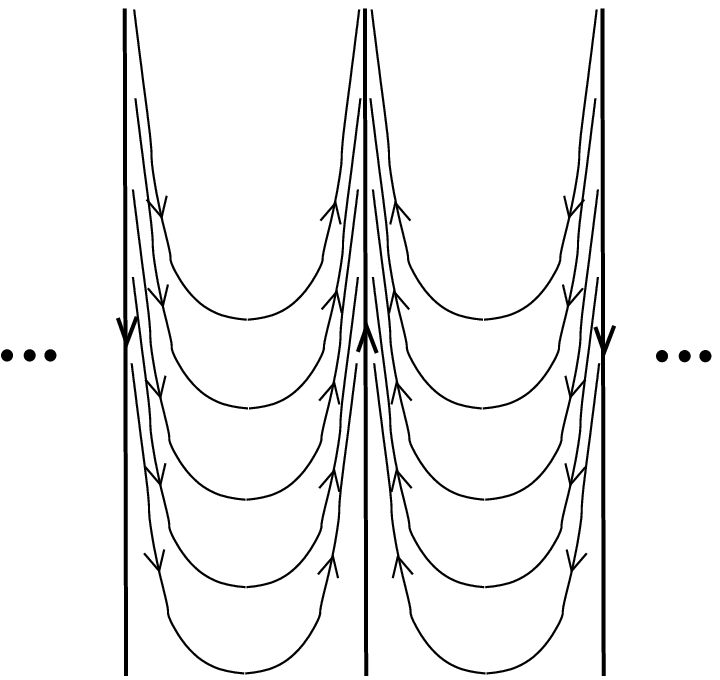}
\nonumber
\]
\end{figure}

First, we want to derive an explicit formula describing the action.
On the special orbits, it is clear
that the action is simply the vertical translation in one or the other direction.
 For a regular orbit, we have to determine the position
of the moved point as a function of the initial point and the length of the
curve joining the two points. (The latter is, by construction, the element of
the group $\R$ by which   the initial point has been moved.) Since all regular
orbits can be obtained by a vertical shift of a chosen orbit, the first coordinate
of the moved point depends only on the first coordinate $x$ of the initial point
and the group element $t$. Let us denote the value of the first coordinate
of the moved point by $g(x,t)$. It is determined by the following
 equation
\[
\int_x^{g(x,t)}\!\!\!\sqrt{1+(\gamma'(s))^2}\;\d s=t.
\]
Solving this equation yields
\bea
g(x,t)&=&\left\{\ba{cr}
\arcsin g_0(x,t)+2\mu\pi&
\mbox{for }x\in (-\frac{\pi}{2}+2\mu\pi,\frac{\pi}{2}+2\mu\pi),\;\mu\in\Z,
\\&\\
\pi-\arcsin g_0(x,t)+2\mu\pi&
\mbox{for }x\in (\frac{\pi}{2}+2\mu\pi,\frac{3\pi}{2}+2\mu\pi),\;\mu\in\Z,
\\
\ea\right.\nonumber\\ &&\\
g_0(x,t)&:=&\frac{e^t(1+\sin x)-e^{-t}(1-\sin x)}{e^t(1+\sin x)+e^{-t}(1-\sin x)}.
\nonumber
\eea
Note that $g$ remains well defined also at the special values of  $x$, i.e.,
for $x=\frac{\pi}{2}+\mu\pi,\;\mu\in\Z$, it is continuous at all these points
and independent of $t$. The latter agrees with the formula for the action on
the special orbits. To find out the behaviour of the second
coordinate under the group action, observe that the second
coordinate of the moved point can be written as $y+\widetilde{g}(x,t)$, where
$(x,y)$ stands for the initial point,  $t$ is the element of
the group by which   the initial point has been moved, and
\[
\widetilde{g}(x,t):=\gamma(g(x,t))-\gamma(x)=
\log\frac{e^t(1+\sin x)+e^{-t}(1-\sin x)}{2}.
\]
Again, the vertical invariance makes it irrelevant which regular orbit $\gamma$
we take. Thus we  have arrived  at the explicit formula for the
continuous action:
\[
\R^2\times\R\ni ((x,y)\;,\;t)\longmapsto (g(x,t)\;,\;y+\widetilde{g}(x,t))\in\R^2.
\]

It is entertaining to verify directly that this formula indeed gives an action,
i.e., that $g(g(x,t),t')=g(x,t+t')$ and
$\widetilde{g}(x,t)+\widetilde{g}(g(x,t),t')=\widetilde{g}(x,t+t')$. Similarly,
one can check that, for any $x\in\R$, if
$
g(x,t)=x
$
and $\widetilde{g}(x,t)=0$, then  $t=0$. Hence the action is free, and
$F^G$ is a continuous injection. On the other hand, the sequence of
pairs of points $((x_n,y_n)\,,\,(x'_n,y'_n))$ given by the formulas
\[
x_n=\frac{1}{n}-\frac{\pi}{2},\;\;\;y_n=0,
\;\;\;
x'_n=-\frac{1}{n}+\frac{\pi}{2},\;\;\;y'_n=0,
\]
is contained in $\R^2\times_{\R^2/\R}\R^2$, but it converges to
$((-\frac{\pi}{2},0)\,,\,(\frac{\pi}{2},0))\not\in\R^2\times_{\R^2/\R}\R^2$.
Consequently, the image of $F^G$ is not closed, so that, by Lemma~\ref{closed},
 $F^G$ is not a proper
map.

Finally,  to show that $F^G_X$ is a homeomorphism, it suffices
to prove that the translation map $\check{\tau}$  is continuous.
 Assume that $(x,y)$ and $(x',y')$ are two points
on the same regular orbit. Solving the equation $g(x,t)=x'$ for $t$ gives
\[\Label{tfor}
t=\log\sqrt{\frac{1-\sin x}{1+\sin x}\frac{1+\sin x'}{1-\sin x'}}
=\log\left(\left|\frac{\cos x}{\cos x'}\right|\frac{1+\sin x'}{1+\sin x}\right)
=\log\left(\left|\frac{\cos x'}{\cos x}\right|\frac{1-\sin x}{1-\sin x'}\right).
\]
Thus $\check{\tau}$  is evidently continuous on the (dense) subset of
$\R^2\times_{\R^2/\R}\R^2$ consisting of all pairs from regular orbits. Remembering
the relationship $y'-y=\log\left|\frac{\cos x}{\cos x'}\right|$, we can re-write
(\ref{tfor}) in the form
\[\Label{tfor2}
t=y'-y+\log\frac{1+\sin x'}{1+\sin x}=y-y'+\log\frac{1-\sin x}{1-\sin x'}.
\]
Now, let
$(-\frac{\pi}{2}+2\mu\pi,y_0)$, $(-\frac{\pi}{2}+2\mu\pi,y'_0)$, $\mu\in\Z$,
be a limit pair of some sequence
\[
\{((x_n,y_n)\,,\,(x'_n,y'_n))\}_{n\in\N}\inc
\R^2\underset{\R^2/\R}{\times}\R^2.
\]
 Then, for any
sufficiently big $n$, both $\sin x_n$ and $\sin x_n'$ are different from
1. Therefore, due to (\ref{tfor2}), for any $x_n$ and $x'_n$
in a small neighbourhood of
$-\frac{\pi}{2}+2\mu\pi$, we have:
\bea
&&
\lim_{n\ra\infty}\check{\tau}((x_n,y_n)\,,\,(x'_n,y'_n))
\nonumber\\ &&=
\lim_{n\ra\infty}\left(y_n-y'_n+\log\frac{1-\sin x_n}{1-\sin x'_n}\right)
\nonumber\\ &&=
y_0-y'_0
\nonumber\\ &&=
\check{\tau}
\left((-\frac{\pi}{2}+2\mu\pi\,,\,y_0)\,,\,(-\frac{\pi}{2}+2\mu\pi\,,\,y'_0)\right).
\eea
Much as above, if $(\frac{\pi}{2}+2\mu\pi,y_0)$, $(\frac{\pi}{2}+2\mu\pi,y'_0)$,
$\mu\in\Z$,
is a limit pair of some sequence  from
$\R^2\times_{\R^2/\R}\R^2$, then
\bea
&&
\lim_{n\ra\infty}\check{\tau}((x_n,y_n)\,,\,(x'_n,y'_n))
\nonumber\\ &&=
\lim_{n\ra\infty}\left(y'_n-y_n+\log\frac{1+\sin x'_n}{1+\sin x_n}\right)
\nonumber\\ &&=
y'_0-y_0
\nonumber\\ &&=
\check{\tau}
\left((\frac{\pi}{2}+2\mu\pi\,,\,y_0)\,,\,(\frac{\pi}{2}+2\mu\pi\,,\,y'_0)\right).
\eea

Summarising, we have shown that $F^G_X$ is a homeomorphism, whereas $F^G$
is a continuous but not proper injection. The lack of properness manifests
itself in the non-Hausdorffness of the quotient space $\R^2/\R$: the neigbouring
special orbits cannot be separated by any open sets.

\vspace*{2mm}\begin{figure}[h]
\[
\includegraphics[width=90mm]{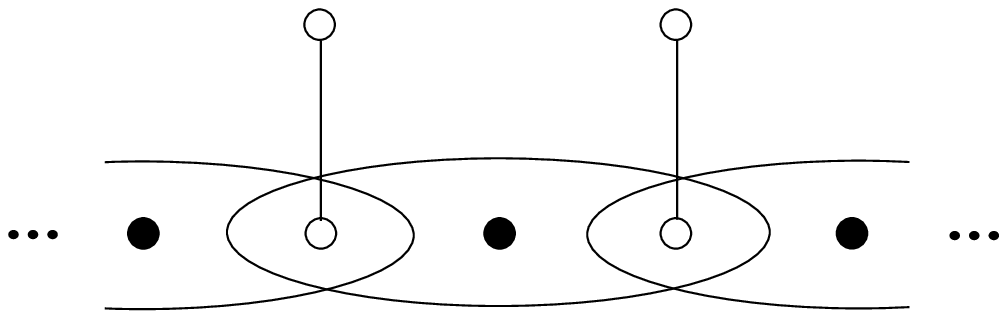}
\nonumber
\]
\end{figure}
\eex

\subsection{Associated bundles}

Let us first shortly recall some basic
 terminology related to topological  bundles.
In great generality, a {\em bundle} is a triple $(E,\pi,M)$, where $E$ and $M$ are topological spaces
and $\pi:E\ra M$ is a contiuous surjective map. Here $M$ is called the base space,
$E$  the total space, and $\pi$  the projection of the bundle.
For $p\in M$, the
fibre  over~$p$ is the topological
space $\pi^{-1}(p)$. A {\em local section} of a bundle is a
continuous map $s:U\ra E$ with $\pi\ci s=\id$, where $U$ is an open
subset of $M$. If each fibre of a bundle is endowed with a 
 vector space structure
such that the addition and scalar multiplication are continuous,
we call it a  {\em bundle of vector spaces}. (If in addition each
fibre is finite dimensional, this coincides with the notion of a
family of vector spaces \cite[p.1]{a-mf89}.)
A bundle $(E,\pi,M)$ all of whose fibres are homeomorphic to a space $F$
  is called
a {\em fibre bundle} with typical fibre $F$, and is denoted by $(E,\pi,M,F)$.

Let a topological group $G$ act from the right on a topological 
space $X$.
Then the triple $(X,\pi,X/G)$, where $X/G$ is the orbit space (with the standard quotient topology)
and $\pi$ is the
natural projection, is a bundle in the above sense. More generally,
we call a bundle a {\em $G$-bundle} iff all the fibres are orbits
of the $G$-action and its base space is homeomorphic to the orbit space.
The projection of a $G$-bundle is necessarily an open map. Indeed, for any open
 $U\inc X$, we have $\pi^{-1}(\pi(U))=\bigcup_{g\in G}Ug$. Therefore, as $G$
 acts by homeomorphisms, $\pi^{-1}(\pi(U))$ is a union of open sets, whence itself
 open. Thus, by the definition of quotient topology, $\pi(U)$ is open, as claimed. 
Finally, if the action is principal, we arrive
at the following fundamental definition:
\bde\label{prinb}
A {\em principal bundle} is a quadruple $(X,\pi,M,G)$ such that
\vspace{-.1cm}
\begin{itemize}
\item[(i)]
$(X,\pi,M)$ is a bundle and $G$ is 
a topological group acting continuously on $X$ from the right,
\item[(ii)]
the action of $G$ on $X$ is principal (i.e., free and proper),
\item[(iii)]
$\pi(x)=\pi(y)$ if and only if $\exists\, g\in G:y=xg$
(the fibres are the orbits of $G$),
\item[(iv)]
the induced map $X/G\ra M$ is a homeomorphism.
\end{itemize}
\ede

Observe 
 that the induced map $X/G\ra M$, existing by (iii), is always a continuous
bijection but need not be a homeomorphism (e.g., if $M$ has the indiscrete 
topology). However, it is an open map (and hence a homeomorphism) if and only if
the bundle projection $\pi$ is an open map.
Next, note that
if $G$ is a topological group acting properly
on a locally compact Hausdorff space $X$, then $X/G$ is also locally compact and
Hausdorff \cite[Proposition~3, p.253; Proposition~9, p.257]{b-n71}.
We say that a principal bundle is (locally) compact iff the involved
spaces are (locally) compact Hausdorff. 

Clearly, a principal action of $G$ on $X$ automatically makes
the bundle $(X,\pi,X/G)$ a principal bundle. However, not every
principal bundle has to be of this form. If we replace $X/G$
by a homeomorphic space, not only we formally define a different
 bundle, but also it might happen that such a new bundle is not
equivalent to $(X,\pi,X/G)$ \cite[p.157]{f-t00}. Thus equivalent
(equivariantly homeomorphic)
$G$-spaces  might lead to inequivalent principal bundles.

Since sometimes one can relax the proper-action assumption to
the continuity-of-translation-map condition, we define a 
{\em quasi principal bundle} exactly  as above except for 
point (ii) that is replaced by:

(ii')  the action of $G$ on $X$ is free and the translation map
is continuous.

\noindent
As argued below Corollary~\ref{injclo}, a quasi principal bundle
is a principal bundle if and only if the fibre product $X\times_MX$
is closed in $X\times X$. Note also that,
due to the continuity of the translation map, each fibre of 
a quasi principal bundle is homeomorphic
to~$G$, so that it is always a fibre bundle. (See (\ref{locsec}) and replace $U$ by a one-point set.)

Let $(X,\pi,M)$ be a $G$-bundle, and let $G$ act on the left on
another topological space~$F$.
Then
$(x,v,g)\mapsto (x g,g^{-1} v)$ defines a right action
of $G$ on $X\times F$. The map $\pi_F: (X\times F)/G\ra M$,
$[(x,v)]\mapsto \pi(x)$, is well defined and continuous.
Much as for a $G$-bundle, this projection
is an open map. To show this, consider the commutative diagram of 
continuous
surjections:
\[
\xymatrix@+1pc{
X\times F \ar[r]^-{\pi_E} \ar[d]_{\text{pr}_1}& (X\times F)/G \ar[d]^{\pi_F}\\
X \ar[r]^{\pi} & M
}
\]
(Here the left vertical and upper horizontal arrows are the obvious surjections.) We already
know that the projections of $G$-bundles are open, and the canonical surjection
$\text{pr}_1$ is open. Taking advantage of this and the surjectivity and
continuity of $\pi_E$, we obtain that $\pi_F(U)$ is open for any open
$U\inc (X\times F)/G$:
\[
\pi_F(U)=\pi_F\llp\pi_E(\pi_E^{-1}(U))\lrp=(\pi\ci\text{pr}_1)(\pi_E^{-1}(U)).
\]

For a quasi principal bundle one can check that
the assignment
$v\mapsto [(x,v)]$ defines a homeomorphism from $F$ to the fibre over~$\pi(x)$.
(Replace $U$ by $\{\pi(x)\}$ in (\ref{ltafb}).)
Thus, $((X\times F)/G, \pi_F,M,F)$ 
is a fibre bundle with typical fibre $F$.
It is called an {\em associated fibre bundle}.

Global continuous sections of an associated fibre bundle can be identified with
continuous equivariant maps defined on $X$ with values in $F$.
Let
\[
\hom_G(X,F):=\{\,f\in C(X,F)\;|\;f(x g)=g^{-1} f(x)\,\},
\]
 and
let $\Gamma(E)$ denote the space of continuous global sections of a bundle $(E,\pi,M)$.
We have:
\ble\label{homma}
Assume that $(X,\pi,M,G)$ is
a quasi principal bundle ($G$ acts freely on $X$ and  the translation map is continuous). Let
$E=(X\times F)/G$ be a fibre bundle associated to  $(X,\pi,M,G)$. Then the formulas
$$
s_f:M\ni \pi(x)\longmapsto[(x,f(x))]\in E,\quad
f_s:X\ni x\longmapsto v \in F,\; s(\pi(x))=[(x,v)],
$$
define mutually inverse  bijections between
$\hom_G(X,F)$ and $\Gamma(E)$.
\ele
\bpf
One easily verifies that  $f_{s}$ and $s_{f}$ are well defined. It is also clear that  $f_{s}$
is equivariant, $\pi\ci s_{f}=\id$, and $f_{s_f}=f$, $s_{f_s}=s$. In order to see that the continuity of
$f$ entails the continuity of $s_f$,  note that $\pi_E\ci(\id,f)\ci\text{diagonal map}=s_f\ci\pi$, where
 $\pi_E$ is the canonical quotient map $X\times F\ra E$. The left hand side is evidently
 continuous, so that the continuity of $s_f$ follows from the fact that $\pi$ is an open surjection (see the second paragraph).

Indeed, for any commutative diagram
\[\label{toptrio}
\xymatrix@+1pc{
A \ar[r]^{f_{AB}} \ar[rd]_{f_{AC}}& B \ar[d]^{f_{BC}}\\
& C
}
\]
where $A,B,C$ are topological spaces, $f_{AC}$ is continuous and $f_{AB}$
is a surjection that maps open sets to open sets (not necessarily continuous),
we can conclude that $f_{BC}$ is continuous. To verify this, we employ the
surjectivity of $f_{AB}$ and compute:
\[
f_{BC}^{-1}(U)=f_{AB}\llp f_{AB}^{-1}(f_{BC}^{-1}(U))\lrp=f_{AB}(f_{AC}^{-1}(U)).
\]
Hence, for any open $U\inc C$, the continuity of $f_{AC}$ and the assumption that
$f_{AB}$ maps open sets to open sets entail that $f_{BC}^{-1}(U)$ is open. 
Thus $f_{BC}$ is continuous, and substituting $s_f$ for $f_{BC}$ 
and $\pi$
for $f_{AB}$ yields the desired conclusion.

Conversely, assume that $s$ is continuous.
 Note first that, since the map
\[
\id\times\pi_E:X\times X\times F\lra X\times E 
\]
is open (see above
Definition~\ref{prinb}) and
$
(\id\times\pi_E)^{-1}(X\times_{M} E)=X\times_{M} X\times F,
$
by \cite[Proposition 2a), p.51]{b-n71}, the restriction of $\id\times\pi_E$
to $X\times_{M} X\times F$ is an open map. Thus we obtain an open surjection
$\overline{\pi_E}:X\times_{M} X\times F\ra X\times_{M} E$. On the other hand,
 the continuity of the translation map implies that
\[
T :\: X\!\underset{M}{\times}\! X\times F\ni (x,y,v)\longmapsto \check{\tau}(x,y)v\in F
\]
is continuous. 
Furthermore, the property (\ref{transequiv}) of $\check{\tau}$ 
implies that $T$ is well defined on the quotient
$X\times_{M} E$. This yields the commutative diagram
\[\label{}
\xymatrix@+1pc{
X\!\underset{M}{\times}\! X\times F \ar[r]^-{\overline{\pi_E}} \ar[dr]_T 
& X\!\underset{M}{\times}\! E \ar[d]^{\overline{T}}\\
& F
}
\]
Arguing as with (\ref{toptrio}), we infer that $\overline{T}$ is
continuous.
 Finally, since $s$ is continuous by assumption and
$f_s(x)=\overline{T}(x,s(\pi(x)))$, we can conclude that $f_s$ is continuous,
as desired.
\epf

If in addition $F$ is equipped with  a vector space structure compatible
with its topology, and the  action of $G$ is given by
a linear representation, the fibres of the associated fibre bundle
\[
((X\times F)/G\,,\, \pi_F\,,\,M\,,\,F)
\]
 carry a natural vector space structure
such that the homeomorphisms $v\to [(x,v)]$ between the typical fibre  and the fibres are linear.
This way we obtain a  bundle of vector spaces as an associated fibre bundle.

\subsection{Local triviality and piecewise triviality}

Local triviality is a notion  commonly
assumed in the classical setting of  bundles --- see standard definitions in
\cite[16.12.1]{d-j74}, \cite[p.40]{h-f78}, \cite[p.13]{ms74}, \cite[p.50]{kn63}.
\bde\label{ltr}
A  bundle $(E,\pi,M)$ is said to be {\em locally trivial},
 if for every $p\in M$ there exists a neighbourhood
$U$, a topological space $F$ and a homeomorphism $\phi:U\times F\ra \pi^{-1}(U)$ that is
fibre preserving, i.e. $\pi\ci\phi:U\times F\ra U$ is the canonical projection.
\ede
 \noindent
A very simple example of a bundle which is not locally trivial is exhibited by
the following picture:
\[
\includegraphics[width=40mm]{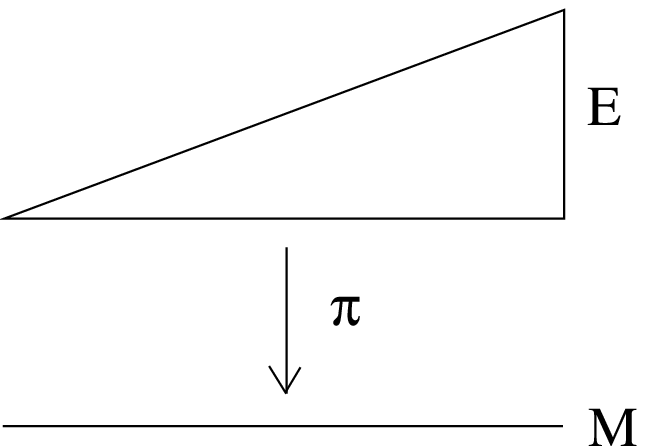}
\nonumber
\]

Among locally trivial bundles, particularly important are
vector bundles --- they are a starting point
for $K$-theory.
\bde
A bundle of vector spaces is called a {\em vector
bundle} iff it is locally trivial, 
the trivialising maps are compatible with the linear structure,
and each of the fibres is finite dimensional.
\ede
\noindent
Note that a vector bundle need not be a fibre bundle. For
instance, one might have 0-dimensional fibres over one connected
 component and 1-dimensional fibres over another connected
 component. However, a vector bundle is always a fibre bundle if
its base space is connected \cite[p.3]{a-mf89}.

Complex vector bundles (i.e.\ 
 vector bundles
  whose  fibres are  complex vector spaces)
are algebraically characterized as finitely generated projective modules over 
commutative $C^*$-algebras.
This is one of the fundamentals of noncommutative geometry 
 \cite[pp.34--36]{r-j94}.
\bth
[Serre \cite{s-jp57}, Swan \cite{s-rg62}]
Let $(E,\pi,M)$ be a  complex vector bundle
and
$M$ a compact Hausdorff space. Then the set $\Gamma(E)$ of all
continuous global sections of $E$
is a finitely generated projective
module over the $C^*$-algebra $C(M)$ of continuous functions on $M$
(with the natural module structure given by pointwise addition and
multiplication).
Conversely, any finitely generated projective module over
$C(M)$ is isomorphic to $\Gamma(E)$ for some 
complex vector
bundle $(E,\pi,M)$.
\ethe

There are generalisations of this theorem for
 general topological  spaces
using the notion of finite-type  bundles of vector spaces. These are finite-dimensional
bundles that are assumed to be locally trivial with respect to a finite covering  
of the base space admitting a partition of unity. 
(The partition-of-unity assumption
allows one to go beyond paracompact spaces.) 
By \cite[Theorem 1]{v-ln86}, the category of 
finite-type  bundles of vector spaces over
$M$ is equivalent to the category of finitely generated projective modules over
the algebra of
all continuous functions on $M$. By \cite[Theorem 2(3)]{v-ln86} and the 
subsequent remark (see also
\cite{s-rg77}),
the same is true for the algebra of bounded functions and the algebra of 
smooth functions, if $M$ is a
smooth manifold.

If $E$ is a vector
bundle, the trivialising map has to be compatible with its linear structure.
For a $G$-bundle $(X,\pi,M)$, the trivialising map has to be $G$-equivariant. A traditional definition of a principal bundle
assumes that it is locally trivial in this sense but does not require
that the group action  is proper
\cite[p.156]{f-t00}. For instance, the bundle studied in 
Example~\ref{rysiu} is locally trival (by Theorem~\ref{lotri} or
direct inspection), so that it fits the 
aforementioned definition,
 but the group action is not proper.
On the other hand, under the assumption of local triviality,
one can remove from Definition~\ref{prinb}
the condition  that   the quotient space
is homeomorphic with the base space by the induced map.  

Indeed, we already know that the induced map $f:X/G\ra M$ is a continuous
bijection. Thus it remains to see that it is an open map. 
Let $\{U_i\}_i$ be an open cover of $M$ locally trivialising $X$, and let
$V$ be an open subset of $X/G$. Then each
$V_i:=V\cap f^{-1}(U_i)$ is an open set, and 
$V=\bigcup_i V_i$. Furthermore, $\pi^{-1}(f(V_i))$ is open
because it is the 
preimage of $V_i$ under the canonical quotient map.
On the other hand, since $f(V_i)\inc U_i$, the preimage
 $\pi^{-1}(f(V_i))$ is mapped by a  homeomorphism onto
$f(V_i)\times G$. Therefore,
$f(V_i)\times G$ is an open subset of $U_i\times G$. 
Consequently,  $f(V_i)$ is an open subset of $U_i$ (and hence $M$),
 because
the projections from the Cartesian product of two sets onto 
its components
are always open maps.
 Finally, it follows from 
the general fact $f(\bigcup_i A_i)=\bigcup_i f(A_i)$ (for any map and 
any family of sets) that $f(V)=\bigcup_i f(V_i)$. This is an open set, 
as desired.

\bpr\label{locsecpr}
If $(X,\pi,M,G)$ is a quasi principal bundle ($G$ acts freely and the translation map is continuous), then local triviality
is equivalent to the existence of local sections at each point of~$X$. 
\epr
\bpf
If there exists a local section $\hs:U\ra X$, then we can construct a continuous map
\[\label{locsec}
\phi_\sigma: U\times G\ni (u,g)\longmapsto \hs(u)g\in\pi^{-1}(U).
\]
Its inverse is given by $x\to(\pi(x),\check{\tau}(\hs(\pi(x)),x))$, so that $\phi_\sigma$
is evidently a trivialising map if $\check{\tau}$ is continuous. Conversely, if we have a
trivialising map $\phi:U\times G\ra \pi^{-1}(U)$, then $u\to\phi(u,e)$ is clearly a local
section.
\epf

The following theorem guarantees that, if $(X,\pi,M,)$ is a principal
bundle and $G$ is a Lie group, then it is a locally trivial bundle.
The compact Lie group version of this result is attributed to A.M.~Gleason \cite{g-am50}, the
general case is announced
 by J.-P.~Serre \cite[Th\'eor\`eme~1]{s-jp50}, and its proof can be
found in \cite[p.315]{p-rs61}.
\bth\label{lotri}
Let $G$ be a Lie group acting freely on a completely regular space~$X$. Then
there exists a local section through each point of $X$ if and only if the translation
map is continuous.
\ethe

An immediate corollary of this theorem is that the quotient of a locally compact group
by a closed Lie subgroup yields a locally trivial bundle. This is because in this situation
the translation map is automatically continuous: $\check{\tau}(g,g')=g^{-1}g'$. Also,
note that, since the translation map in Example~\ref{rysiu} is continuous, it is an example
of a locally trivial bundle over a non-Hausdorff space. On the other hand, the translation map
for the Kronecker foliation (ergodic action of $\R$ on $T^2$) is not continuous. Thus this is
an example of a free action of a Lie group that does not allow local triviality.

\bex\noindent\label{nlotri}
In the light of Theorem~\ref{lotri}, to construct an example of a principal
bundle which is not locally trivial,  we
have to stay clear of Lie groups. What follows is a simple example of a compact
principal bundle where we avoid the Lie group structure by constructing a group with the
topology of the Cantor set.
Let $X$ be the infinite product of unitary groups $\prod_{\mathbb{N}} U(1)$,
and $G$ its subgroup that is the infinite product $\prod_{\mathbb{N}}(\mathbb{Z}/2\mathbb{Z})$ 
of the 2-element group $\mathbb{Z}/2\mathbb{Z}$ viewed as the subset $\{-1,1\}$ of $U(1)$. 
Recall that the standard product
topology is generated by the subsets 
$\prod_{i\in I} U_i\inc\prod_{i\in I} X_i $, where there exists a finite subset $J\inc I$
such that $U_i$ is open
in $X_i$ for all $i \in J$  and $U_i=X_i$  for all $i \not\in J$. With this topology, the
group $X$
is a compact Hausdorff space by the Tichonov Theorem, it is indeed a compact group
(see \cite[Problem 8.5.4.]{e-r77}) and $G$ is a closed subgroup. 
Since the natural action 
of any subgroup on a group is free and any action of a compact group is proper,
$(X,\pi,X/G,G)$ is automatically a principal $G$-bundle.

Suppose now that this bundle is locally trivial. Then there exists an open subset $U$ of
$X/G\cong\prod_{\mathbb{N}}S^1$ over which we have a local section. On the other
hand, $U$  is a union of generating 
subsets defined above. Hence it must contain an open
subset  of the form
$\prod_{i=0}^nU_i \times \prod_{ i=n+1}^\infty S^1$, where each $U_i\inc S^1$ is open.
However, this implies that the principal $\mathbb{Z}/2\mathbb{Z}$-bundle 
$U(1)\ra S^1$ has a global section. Thus we arrive at the desired contradiction proving
that the principal bundle $(X,\pi,X/G,G)$ is not locally trivial. 

Let us now prove that 
$\prod_{\mathbb{N}}(\mathbb{Z}/2\mathbb{Z})$
 with the  (product) topology inherited from
 $\prod_{\mathbb{N}}U(1)$
 is homeomorphic to the Cantor set with the topology it receives as a subspace
of $[0, 1]$.
The Cantor set $\mathrm{C}$ is, by definition,
obtained by removing all open ``middle-third'' intervals from $[0, 1]$. 
Equivalently, if we write an arbitrary $t\in [0, 1]$ as 
$t = \sum_{j=1}^{\infty}{a_j}3^{-j}$,
where $a_j\in\{0,1,2\}$, then $\mathrm{C}$ consists of all those $t\in[0, 1]$ such that
all the $a_j$ are either 0 or 2, i.e., 1 does not occur as an $a_j$.
Since the coefficients $a_j$ of $t\in\mathrm{C}$ are uniquely determined,
any such $t$ can be viewed as a sequence of 0's and
2's. This establishes an evident bijection from the $\mathrm{C}$ to
$\prod_{\mathbb{N}}(\mathbb{Z}/2\mathbb{Z})$  
mapping a sequence of 0's and 2's to the sequence of $-1$'s and $1$'s
obtained by replacing each 0 by 1 and each 2 by $-1$.

It now remains to prove that this bijection is a homeomorphism. Note that, since both $\mathrm{C}$
and $\prod_{\mathbb{N}}(\mathbb{Z}/2\mathbb{Z})$ are compact Hausdorff
spaces, we need only prove continuity in one direction. Also, it is clear that it
suffices to check that the preimages of the generating open sets are open.
If $U$ is a generating open set in $\prod_{\mathbb{N}}(\mathbb{Z}/2\mathbb{Z})$,
then its preimage  consists of all those $t\in\mathrm{C}\subset [0, 1]$
satisfying conditions on finitely many $a_j$'s. Thus, the preimage of $U$ is the
union of subsets of $\mathrm{C}$ defined by fixing the first $n$ $a_j$'s and
taking all possible $a_j$'s for $j>n$.
Since such sets are of the form $\mbox{(open interval)}\cap\mathrm{C}$,  they are open.
Hence the preimage of $U$ is open, as needed.
\eex

In the setting of $C^*$-algebras, the concept of fibre products
\cite{p-gk99} replaces the notion of gluing of topological spaces. 
 In general, a gluing of Hausdorff spaces need not be Hausdorff.
For instance, the gluing of two closed unit intervals along the open
 unit interval is not a Hausdorff space.
On the other hand, 
a gluing of  Hausdorff spaces along closed subsets is always a 
Hausdorff space (\cite[p.135, Exercise 8]{b-n71}).
Since the disjoint union of finitely many 
compact spaces is compact and the quotient map is always continuous, the gluing of compact spaces is always compact.
Thus 
we can conclude that the
gluing of compact Hausdorff spaces over a closed subset is always
compact Hausdorff.

Spaces obtained by gluing compact Hausdorff spaces along closed subsets can be 
viewed as covered by 
the spaces from which they were glued. 
Such coverings, however, are coverings by finitely many closed subsets. They are different from the usual
open covers.
We now discuss the relationship of these two concepts
in the context of local triviality of principal bundles. 
\bde
A principal bundle $(X,\pi,M,G)$ is said to be {\em piecewise trivial},
 if there exists 
a covering of $M$ by finitely many closed sets $W_1,\ldots,W_n$,  
and fibre-preserving $G$-equivariant homeomorphisms 
$\phi_i:\pi^{-1}(W_i)\ra W_i\times G$, 
$\forall i\in\{1,\ldots,n\}$.
\ede 

The assumption of piecewise triviality
implies that the induced map $f:X/G\ra M$ is a homeomorphism. Since
it is a continuous bijection, it suffices to show that it is a 
closed map. 
This can be shown reasoning much as under the assumption of local 
triviality. Instead of the openess of the projections onto components,
we use the fact that a non-empty
subset of a Cartesian product is closed if
and only if each of the factors is closed \cite[p.48]{b-n71}. Then,
in the final argument, we take  advantage of the fact that
a finite union of closed sets
is closed.
Also, the analogue of Proposition \ref{locsecpr} is true in this setting.

Let $M$ be a paracompact space. Then the existence of a partition of
 unity implies that
every open covering has a subordinate covering by closed sets.
As a consequence, every locally trivial bundle over a compact Hausdorff space
is piecewise trivial with respect to a finite closed cover. However,
the converse is not true. 

\bex\label{clnop}
Let $\mathrm{C} \subset [0,1]$ be the Cantor set. As in 
Example~\ref{nlotri}, 
it can be given the structure of a compact abelian group.
Define a base space $M$ to be 
the gluing along $\mathrm{C}$ of two copies $I_1$, $I_2$, of $[0,1]$. 
Denote by $ x \mapsto [x]$ be the quotient map from the disjoint
union of $I_1$ and $I_2$ to $M$.
Then $M$ is a compact Hausdorff space with a closed covering 
given by the two sets $[I_1]$, $[I_2]$. 
Define now a total space $X$ as the gluing of $I_1\times \mathrm{C}$ and 
$I_2\times \mathrm{C}$
by the identifying homeomorphism  
\[
I_1\times \mathrm{C} \supseteq\mathrm{C}\times\mathrm{C}\ni
(x,y)\longmapsto(x,xy)
\in\mathrm{C}\times\mathrm{C}\subseteq I_2\times \mathrm{C}.
\]
 Since $\mathrm{C}\times\mathrm{C}$ is a closed subset of 
$[0,1]\times C$, the total space $X$ is also compact Hausdorff.
 It is straightforward to verify that
the obvious right action of $\mathrm{C}$ on the disjoint union of
$I_1\times \mathrm{C}$ and $I_2\times \mathrm{C}$ descends to a 
continuous free
action on $X$. The action is obviously proper because C is a compact 
group. The obvious homeomorphism from the orbit space $X/\mathrm{C}$
to $M$ defines a continuous surjection $\pi:X\ra M$.
Thus $(X,\pi,M,\mathrm{C})$ is a piecewise trivial compact principal
bundle.

Suppose now that $X$ admits a  a global
section. Then there exist continuous functions $f_1:I_1\ra \mathrm{C}$
 and 
$f_2:I_2\ra \mathrm{C}$  satisfying the compatibility condition 
$f_2(x)=xf_1(x)$
for all $x\in \mathrm{C}$. Since $I_1$ and $I_2$ are connected and 
the Cantor set
is totally disconnected, both $f_1$ and $f_2$ must be constant maps.
Therefore, the compatibility condition evidently cannot be satisfied,
so that $X$ is a non-trivial principal bundle. To prove that it is also
not locally trivial, suppose that it is trivial over
 an open subset $U$ of $M$
containing an element of $\mathrm{C}$.
Then there exists a local section over $U$, and the same argument as before  applies to prove that this is impossible. (Any open set
containing an element of $\mathrm{C}$ contains a copy of $M$.)
Summarising, $(X,\pi,M,\mathrm{C})$ is a piecewise trivial compact principal that is not locally trivial.

To end with, let us show that
the space $M$ is homeomorphic to a subset of
the plane that can be visualised as follows (the bubble space):

\vspace*{-0mm}\begin{figure}[h]\label{bub}
\[
\includegraphics[width=90mm]{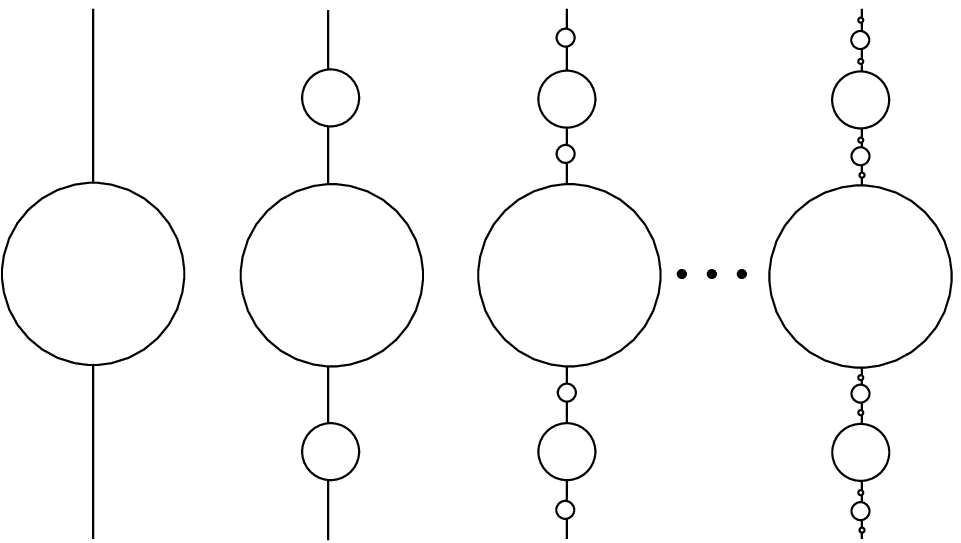}
\nonumber
\]
\end{figure}\vspace*{-0mm}
\noindent
It is created the same way one creates the Cantor set, only now instead
of removing middle third intervals we replace them by circles whose
diameters are equal to the lengths of replaced intervals. By cutting
the bubble space vertically along the middle axis we obtain two 
homeomorphic subsets of a plane. They can be identified with the
graph of a continuous function on $[0,1]$. Indeed, the iterating
procedure defining the bubble space can be easily translated into
a suquence of uniformly convergent functions whose limit is the
aforementioned function. 
\note{
In order to prove that $M$ is indeed homeomorphic to the above subset of the plane,
let us first see that the unit interval $[0,1]$ is homeomorphic to the set $L$ defined as
the bubble space minus the union of all the open semicircles on the right. 
To this end, identify the plane with $\C$
and define a sequence of functions
$f_n:[0,1]\ra\C$ as follows: 
\beq
f_0(t)=(0,it), ~~f_1(t)=\left\{\ba{ll}f_0(t)& t\in[0,\frac{1}{3}]\cup[\frac{2}{3},1]\\
(-\sqrt{\frac{1}{36}-(t-\frac{1}{2})^2},it)& t\in[\frac{1}{3},\frac{2}{3}]\ea\right.,
\eeq
\beq
f_2(t)=\left\{\ba{ll}f_1(t)& t\in [0,\frac{1}{9}]\cup[\frac{2}{9},\frac{7}{9}]\cup[\frac{8}{9},1]\\
(-\sqrt{(\frac{1}{18})^2-(t-\frac{3}{18})^2},it)& t\in[\frac{1}{9},\frac{2}{9}]\\
(-\sqrt{(\frac{1}{18})^2-(t-\frac{15}{18})^2},it)& t\in[\frac{7}{9},\frac{8}{9}]
\ea\right..
\eeq 
For defining $f_n$,
let $R_n$ denote the part of $[0,1]$ that is removed in the $n$-th step of
the construction of the Cantor set. $R_n$ consists of subintervals of length $(1/3)^n$.
Let $S_{n,k}$, $k=1,\ldots,2^{n-1}$, denote the $k$-th of those subintervals, and $t_{n,k}$ its
midpoint. Then
\beq
f_n(t)=\left\{\ba{ll}f_{n-1}(t)& t\in [0,1]\setminus R_n\\
(-\sqrt{(\frac{1}{2\cdot 3^n})^2-(t-t_{n,k})^2},it)& t\in S_{n,k},~k=1,\ldots,2^{n-1}
\ea\right..
\eeq
Thus, the function $f_n$ sends elements of $R_n$
  to points on little semicircles of radius $1/2(1/3)^n$ based on corresponding
subintervals of the part $\{0\}\times R_n$ of the imaginary axis. On $[0,1]\setminus R_n$ it coincides
with $f_{n-1}$. All $f_n:[0,1]\ra\C$ are continuous and form a Cauchy sequence with respect
to the sup-norm. Indeed, 
$\sup_{t\in[0,1]}|f_n(t)-f_k(t)|\leq 1/2(1/3)^n
$ $\forall k\geq n$.
Thus there exists the uniform limit $f=\lim_{n\to\infty}f_n$ which is a 
continuous function
whose image is exactly $L$. 
}
Next, let $F_1:I_1\ra\C$ and $F_2:I_2\ra\C$ denote two copies of this 
function
mapping onto the left and right ``halves" of the bubble space,
 respectively. Since $F_1$ and $F_2$ agree only on the Cantor subsets
of $I_1$ and $I_2$, they define a continuous bijection from $M$ onto
the bubble space. This continuous bijection is a homeomorphism because
$M$ is compact and the bubble space is Hausdorff.  
\eex

To complete the picture, we show that 
the principal bundle of Example~\ref{nlotri} not only is not locally
trivial, but also it is not piecewise trivial. Note first that for
any finite  cover $\{W_i\}_i$ of $M$ there  exists a set $W_j$ such that
$\bigcup_{i\neq j}W_i\neq M$. If the cover is also closed, then
$M\setminus\bigcup_{i\neq j}W_i$ is a non-empty open set contained in
$W_j$. Thus every finite closed cover has an element 
with a non-empty interior.\footnote{We owe this simple argument to 
S.~L.~Woronowicz. One can also argue using 
\cite[Proposition 5, p.24]{b-n71}.
}
In the case of Example~\ref{nlotri}, local triviality over a closed set 
with non-empty
interior would imply the existence of a local section over 
a non-empty open set.
This, however, is explicitly proven in Example~\ref{nlotri} to be
 impossible.

Finally, let us  comment on the hierarchy of examples we have considered.
Our first example was the action of $\R$ on the 2-torus by translation
 under an irrational angle
(leading to the Kronecker foliation).
This is a free action which is not proper and has no continuous
 translation map. 
In Example~\ref{rysiu}, we have a free action, which is not proper, 
but has a continuous translation map. This yields a 
quasi principal bundle. (Notice that by removing 
the vertical lines where the orbits accumulate we get a  
trivial principal bundle, and that the whole bundle is locally trivial.) 
Example~\ref{nlotri} provides an  action that is principal, but
 not piecewise trivial or locally trivial. (The base space in this
example is compact Hausdorff, so that the fact that it is not
piecewise trivial implies that it is also not locally trivial.)
Finally,
Example~\ref{clnop} shows  a principal bundle 
that is piecewise trivial but not locally trivial.

If $(E,\pi,M,F)$ is a locally trivial fibre bundle
and  $U_1$ and $U_2$ are trivialising neighbourhoods with homeomorphisms
$\phi_1:U_1\times F\ra \pi^{-1}(U_1)$ and  $\phi_2:U_2\times F\ra \pi^{-1}(U_2)$,
then $\phi_{21}:=\phi_2^{-1}\ci\phi_1$ is defined as a map
$(U_1\cap U_2)\times F\ra (U_1\cap U_2)\times F$. This map is a homeomorphism
of the form
$
\phi_{21}(p,v)=(p,\vartheta_{21}(p,v)),
$
with continuous $\vartheta_{21}:(U_1\cap U_2)\times F\ra F$.
For three trivialising neighbourhoods, the restrictions to $(U_1\cap U_2\cap U_3)\times F$
satisfy
$
\phi_{31}=\phi_{32}\ci\phi_{21}.
$
In particular, if $\{U_i\}_{i\in I}$ is an open covering of $M$, where $U_i$
are trivialising neighbourhoods with homeomorphisms $\phi_i:U_i\times F\ra
\pi^{-1}(U_i)$, then there is an analogous equation for every triple of
indices, i.e.
\beq\label{cocy}
\fa i,j,k\in I:\;\phi_{ki}=\phi_{kj}\ci\phi_{ji}\,.
\eeq
(This includes $\phi_{ii}=\id$
and $\phi_{ij}=\phi_{ji}^{-1}$.) Let us call the maps $\phi_{ij}$ the gluing functions of
the bundle.

On the other hand, if one is given an open covering $\{U_i\}_{i\in I}$ of $M$ and
homeomorphisms $\phi_{ij}:(U_i\cap U_j)\times F\ra (U_i\cap U_j)\times F$ fullfilling
the ``cocycle condition'' (\ref{cocy}), then one can construct, by the standard gluing
procedure,
a locally trivial fibre bundle over $M$ with
typical fibre $F$. The bundle comes with canonical local trivialisations
whose gluing functions coincide with the initial maps $\phi_{ij}$ (cf.
\cite[16.13]{d-j74}). 

Consider now  a locally trivial  $G$-bundle given by a free action inducing
a continuous translation map.
The locally trivialising maps $\phi_i$ are $G$-equivariant and
thus of the form
$
\phi_i(p,g)=\sigma_i(p) g,
$
 where $\sigma_i:
U_i\ra \pi^{-1}(U_i)$ is a local section (see (\ref{locsec})).
 The corresponding gluing functions have the
form 
\[
\phi_{ij}(p,g)=\phi_i^{-1}(\phi_j(p,g))=\phi_i^{-1}(\hs_j(p)g)=
\llp p,\check{\tau}\llp\hs_i(p),\hs_j(p)\lrp g\lrp
=:(p,\theta_{ij}(p)g).
\]
The continuous maps
$
\theta_{ij}:=\check{\tau}\ci(\hs_i\times\hs_j)\ci\text{diagonal map}
:U_i\cap U_j\ra G
$
 are called the {\em transition functions} of the
bundle. For any triple of indices $i,j,k$, they satisfy  the cocycle condition
(cf. \cite[Exercice 1, 16.14]{d-j74})
\beq\label{cocp}
\theta_{ij}(p)=\theta_{ik}(p)\theta_{kj}(p),\quad \fa p\in U_i\cap U_j\cap U_k\,.
\eeq
Conversely, one can reconstruct the bundle from a given
covering and transition functions $\theta_{ij}$ satisfying the
cocycle condition (\ref{cocp}) (again cf. \cite[Exercice 1, 16.14]{d-j74}).

For a $G$-bundle $(X,\pi,M)$ as above, the associated 
bundle $(X\times_GF,\pi_F,M,F)$  is  also locally
trivial. It has natural locally trivialising maps
\[\label{ltafb}
\phi_F:U\times F\lra \pi_F^{-1}(U),\quad
\phi_F(p,v)=[(\sigma(p),v)],
\]
 where $\sigma:U\ra X$
is a local section.
Indeed, $\phi_F$ is a continuous fibre-wise map and its inverse is given by
\[
\phi_F^{-1}([(x,v)])=(\pi(x),\check{\tau}(\hs(\pi(x)),x)v).
\]
Hence, since $\check{\tau}$ is continuous, $\phi_F$ is a homeomorhism,
and thus a locally trivialising map.
The transition functions of the associated
bundle have the form
$(p,v)\mapsto (p,\theta_{ij}(p)v)$, where $\theta_{ij}$ are the
transition functions of $(X,\pi,M,G)$.

\setcounter{equation}{0}
\section{Elements of general algebra}

Introducing geometric structures on topological spaces needs some algebraic concepts. Differential objects such as differential forms,
connections and  sections of vector bundles have their 
algebraic backbones in the universal differential algebra, splittings
of surjections on projective modules, and the 
cotensor product, respectively.
In this section, we recall the aforementioned elements of algebra 
and explain in detail certain basic mechanisms. Throughout this part,
we work with arbitrary associative unital algebras over a ground
field $k$. We assume that all our vector spaces and unadorned tensor 
products are over~$k$.

\subsection{Modules and comodules}

 Herein, we gather a number of basic facts about actions and coactions.
 For more details, the interested reader is referred to
  textbooks, e.g., for modules and comodules
  in the context of Hopf algebras  to \cite{s-me69,a-e04},
   and  in the context of corings to \cite{bw03}.

\note{ Recall that an {\em algebra} (properly, a $k$-algebra) is a
$k$-vector space $A$ with an associative product $m: A\otimes A\to
A$. A typical example is the algebra of functions (polynomial,
continuous, smooth etc.) on a space (topological space, manifold
etc.) with a product given by the point-wise multiplication.
Unless stated otherwise, we assume that algebra has a unit element
(e.g.\ the algebra of continuous functions on a compact space).
The product of elements of an algebra $A$ is denoted by
juxtaposition, i.e., $aa' := m(a\otimes a')$. A vector space $M$
with an associative (and unital) left (resp.\ right) action of an
algebra $A$ is called a {\em left} (resp.\ {\em right}) {\em
$A$-module}. Frequently, an $A$-module is called a {\em
representation} of $A$.  On elements, the action is denoted by a
dot, i.e.\ $a\cdot x$ denotes the left action of $a\in A$ on $x\in
M$, and $x\cdot a$ denotes a right action of $a\in A$ on $x\in M$.
}

A vector space $M$ that, at the same time, is a left module of an
algebra $A$ and a right module of an algebra $B$ with mutually
commuting actions is called an $(A,B)$-bimodule.
 Recall that, given a right $A$-module $M$ and a left $A$-module $N$,
 the {\em tensor product of $A$-modules} $M$ and $N$  is a vector space
 $M\otimes_AN$ defined by the following exact sequence
\[
M\otimes A\otimes N \st{\omega_A}{\lra} M\otimes N \lra M\ten{A}
N\lra 0,
\]
where $\omega_A$ is given by $\omega_A: m\otimes a\otimes n\mapsto
ma\otimes n - m\otimes an$. This means that, for all $a\in A$,
$m\in M$ and $n\in N$, $ma\otimes_A n= m\otimes_A an$ (cf.
(\ref{fibred}). The tensor product is a bifunctor from the
categories of modules to the category of vector spaces. In
particular, this means that, for any left $A$-module $M$ and any
right $A$-linear map $f: N\ra\tilde{N}$, the map
\[
f\ten{A} \id : N\ten{A} M \lra \tilde{N}\ten{A} M, \qquad n\ten{A}
m\longmapsto f(n)\ten{A} m
\]
is a homomorphism of vector spaces. A left $A$-module $M$ is said
to be {\em flat}, provided for any monomorphism (injection) $f$ of
right $A$-modules, the map $f\otimes_A \id_M$ is also injective.
This is equivalent to the statement that $M$ is
 flat if and only if any (short) exact sequence of right $A$-linear maps remains
 exact after tensoring with $M$.

Recall that the dual concept to that of an algebra is the notion
of a {\em coalgebra}. Thus, a coalgebra is a vector space $C$ with
a coassociative {\em coproduct} $\Delta: C\ra C\otimes C$ and a
{\em counit} $\eps: C\ra k$. Henceforth, to denote the action of
$\Delta$ on elements we use the {\em Sweedler notation}
\[
\Delta(c) = c\sw 1\otimes c\sw 2, \quad
((\Delta\otimes\id)\circ\Delta)(c) = ((\id\otimes
\Delta)\circ\Delta)(c) = c\sw 1\otimes c\sw 2\otimes c\sw 3,
\]
(summation understood) etc. Given coalgebras $C$, $D$, a {\em
coalgebra map}\index{coalgebra map} is a $k$-linear map $f: C\ra
D$ such that $\Delta\circ f = (f\otimes f)\circ \Delta$ and
$\eps\circ f = \eps$.

A {\em left} (resp.\ {\em right}) {\em $C$-comodule} or {\em
corepresentation} is a vector space $M$ with a counital and
coassociative coaction ${}_M\Delta: M\ra C\otimes M$ (resp.\
$\Delta_M: M\ra M\otimes C$). On elements the left $C$-coaction is
denoted by ${}_M\Delta(x) = x\sw{-1}\otimes x\sw{0}$, and the
right $C$-coaction is denoted by $\Delta_M(x) = x\sw 0\otimes x\sw
1$ (summation understood). A vector space $M$ that is at the same
time a left $C$-comodule and a right $D$-comodule with mutually
commuting coactions is called a {\em $(C,D)$-bicomodule}.

For any pair of right (resp.\ left) $C$-comodules $M$, $N$, a
linear map $f: M\ra N$ is said to be {\em colinear}\index{colinear
map} if it commutes with the coactions, i.e.\ $\Delta_N\circ f =
(f\otimes\id)\circ \Delta_M$ (resp.\ ${}_N\Delta\circ f =
(\id\otimes f)\circ {}_M\Delta$). The $k$-vector space of all
right (resp.\ left) $C$-colinear maps $M\ra N$ is denoted by
$\hom^C(M,N)$ (resp.\ ${}^C\!\hom (M,N)$).

Next, let us recall the notion of a cotensor product. Given a
coalgebra $C$,  let $M$ be a right $C$-comodule and $N$ be a left
$C$-comodule. The {\em cotensor product} $M\Box_{C}N$ is defined
by the exact sequence
\[
0\lra M\coten{C}N\lra M\ot N\st{\omega^C} {\lra}M\ot C\ot N.
\]
Here $\omega^C$ is the coaction equalising map $\Delta_M\ot
\id_N-\id_M\ot\, _N\Delta$ (cf. (\ref{fibred})).

Since there are no unit elements  in coalgebras,  it is not clear
how to define elements that are invariant under coalgebra
coactions (cf. (\ref{hopfinv})). However, if a coalgebra $C$
coacts (on the right) on an algebra $P$, then, following
M.~Takeuchi, we can define the {\em subalgebra of invariants}
\index{subalgebra of invariants} $P^{coC}\subseteq P$ by
\[\label{invsub}
P^{coC} := \{b\in P\;|\; \forall p\in P: \; \Delta_P(bp) =
b\Delta_P(p)\}.
\]
The subalgebra of invariants of an algebra and a left $C$-comodule is
defined in an analogous way.
\note{
 subalgebra of $P$ follows by the following simple calculation,
for all $b, b^{\prime} \in P^{coC}$ and $p \in P$, \beq \Delta_P(b
b^{\prime} p) =b \Delta_P(b^{\prime} p)=b b^{\prime} \Delta_P(p).
\eeq Thus $bb'\in P^{coC}$, and since, clearly, $1\in P^{coC}$, we
conclude that $P^{coC}$ is a subalgebra of~$P$. }

This definition of invariants immediately implies that, for any
left $C$-comodule $V$, the cotensor product $P\Box_CV$ is a left
$P^{coC}$-module with the action given by
$ b \sum_i p_i\otimes v_i  =  \sum_i bp_i\otimes v_i. $ Similarly,
for any right $C$-comodule $W$, the space $\hom^C(W,P)$ is a left
$P^{coC}$-module with the action $(bf)(w) = bf(w)$. The
relationship between these two modules is described in the
following: 
\ble\label{cote} 
Let $C$ be a coalgebra and $V$ a
finite-dimensional left $C$-comodule.  Let $v_1,\ldots,v_n$ be a
basis of $V$ and $v^1,\ldots,v^n$ the corresponding dual basis of
$V^* := \hom(V,k)$. Then  $V^*$  is a right $C$-comodule by the
formula
$$
\Delta_{V^*}(f) = \sum_i v^i\otimes {v_i}\sw{-1}f({v_i}\sw{0}).
$$
Moreover, for any algebra and right $C$-comodule $P$ and for any
finite-dimensional left $C$-comodule $V$, we have
 $P\Box_C V\cong \hom^C(V^*,P)$
as left $P^{coC}$-modules. 
\ele 
\bpf
 The coaction formula is a
standard result in comodule theory (cf.\ \cite[3.11]{bw03}). Next,
one easily checks that the isomorphism of vector spaces
\begin{equation}\Label{hom-ten}
P\otimes V\lra \hom(V^*,P), \qquad p\otimes v\longmapsto [f\mapsto
pf(v)],
\end{equation}
with the inverse $\phi\mapsto \sum_i\phi(v^i)\otimes v_i$,
restricts to the isomorphism
 $P\Box_{C} V\cong \hom^C(V^*,P)$ (cf.\ \cite[10.11]{bw03}).
 The form of the map  \eqref{hom-ten} immediately implies that this is an isomorphism
 of left $P^{coC}$-modules.
\epf

Having discussed basic facts concerning actions of algebras and
coactions of coalgebras, let us pass to 
coactions  of bialgebras.  
Recall first that an algebra $H$ with a coalgebra structure such
that the coproduct and counit are algebra maps is called a {\em
bialgebra}.\index{bialgebra} Given a bialgebra $H$, a {\em right
$H$-comodule algebra}\index{comodule algebra} is an algebra and a
right $H$-comodule $P$ such that the coaction is an algebra map.
In this case the subalgebra of invariants of $P$ can be defined
as
\[\label{hopfinv}
P^{co H} = \{b\in P\; |\; \Delta_P(b) = b\otimes 1_H\}.
\]
One easily checks that this definition coincides with that of the 
subalgebra of invariants of $P$ with respect to the coaction of {\em
coalgebra} $H$ in equation~\eqref{invsub}. A left $H$-comodule
algebra and its subalgebra of invariants are defined in an analogous
way.

\subsection{The universal differential calculus and algebra}
\label{uda}

The {\em universal differential calculus $\hO^1A$}\index{universal
calculus} is a bimodule defined by the exact sequence \beq\Label{bes}
0\lra\hO^1A\lra A\ot A\lra A\lra 0\, , \eeq i.e.\ as the kernel of
the multiplication map. The differential is given by $\d a:=1\ot
a-a\ot 1$. We can identify $\hO^1A$ with $A\ot A/k$ as left
$A$-modules via the maps \beq\Label{iden} \hO^1A\ni\sum_i a_i\ot
a'_i\to\sum_i a_i\ot \pi_A(a'_i)\in A\ot A/k \ni x\ot \pi_A(y)\to
x\d y\in\hO^1A\, , \eeq where $\pi_A:A\ra A/k$ is the canonical
surjection. Similarly, one can identify $\hO^1A$ with $A/k\ot A$
as right $A$-modules ($\sum_i a_i\ot a'_i\to\sum_i \pi_A(a_i)\ot
a'_i$).

Note that the sequence \eqref{bes} is split-exact as a sequence of
right (or left) $A$-modules. The multiplication map has a right
$A$-linear splitting $a\mapsto 1\ot a$ and a left $A$-linear
splitting $a\mapsto a\ot 1$. This implies that, for any left
$A$-module $N$ (no need for $N$  to be flat), also the sequence
\beq\label{nseq} 0\lra\hO^1A\ten{A}N\lra A\ot N\lra N\lra 0\, \eeq
is exact. (The other-sided version is analogous.)

Note next that if $B$ is a subalgebra of $P$, then we can write
$(\hO^1B)P$ for the kernel of the multiplication map $B\ot P\ra
P$. Indeed, $m((\hO^1B)P)=0$, and if $\sum_i b_i\ot
p_i\in\ker(B\ot P\st{m}{\ra} P)$, then 
\beq \sum_i b_i\ot p_i
=\sum_i (b_i\ot p_i-1\ot b_i p_i) =-\sum_i (\d
b_i)p_i\in(\hO^1B)P\, . 
\eeq 
Furthermore, it follows from the
exactness of (\ref{nseq}) that 
\beq\Label{id2} \hO^1B\ten{B}
P\cong(\hO^1B)P\quad \text{(cf.~\cite[p.251]{hm99})}. 
\eeq

In order to define curvature and also for some problems in cyclic cohomology,
it is helpful to introduce the concept of universal differential algebra. First,
let us recall that the $n$-th order universal differential calculus is
\[
\hO^n\!A:=
\underset{\mbox{$n$-times}}
{\underbrace{\hO^1\!A\underset{A}{\ot}\cdots\underset{A}{\ot}\hO^1\!A}}.
\]
We assume that $\hO^0\!A=A$ and define the {\em universal differential algebra}
of $A$ \cite{k-m82,k-m83} as
\[
\hO A:=\bigoplus_{n\in\N}\hO^n\!A.
\]
It is understood as the tensor algebra of the $A$-bimodule $\hO^1\!A$, so that its product
structure is given by $\ha\hb=\ha\ot_A\hb$.

Now, to turn it into a differential algebra, we need to extend our differential
$\d:A\ra\hO^1\!A$ to the higher order calculi. To this end, we first generalize
the identification (\ref{iden}) to an arbitrary degree. One can easily check that
the composite maps
\[
\xymatrix@+1pc{
\hO^n\!A \ar[r]^-{\subseteq\quad} 
& *+<1ex, 0ex>!<0ex, -1.3ex>{\overset{\mbox{$n$-times}}
{\overbrace{A^{\otimes 2}\underset{A}{\ot}\cdots\underset{A}{\ot}A^{\otimes 2}}}}
\ar[r]^-{\cong} & A^{\otimes n+1} \ar[r] &  A\ot (A/k)^{\otimes n} & *!<5ex,0ex>{\text{and}} 
}
\]
\[
\xymatrix@+1pc{
{A\ot (A/k)^{\otimes n}} \ar[r]^-{\id\otimes\bar\d^{\otimes n}}&
A\ot(\hO^1\!A)^{\otimes n} \ar[r]^-{\quad\phantom{\subset}\quad}
& *+<1ex, 0ex>!<0ex, -1.3ex>{\overset{\mbox{$n$-times}}
{\overbrace{\hO^1\!A\underset{A}{\ot}\cdots\underset{A}{\ot}\hO^1\!A}}}
\ar[r]^-{\quad\;\cong\quad } & \hO^n\!A 
}
\]
are mutually inverse isomorphisms of left $A$-modules. Here 
$\bar\d:A/k\ni[a]\to1\ot a-a\ot1\in\hO^1\!A$, and other maps are constructed
from the canonical surjections and the multiplication isomorphism 
$A\ot_A M\cong M$, for any left $A$-module $M$.

As the next step, we define a map $A\ot(A/k)^{\otimes n}\ra A\ot(A/k)^{\otimes n+1}$,
\[
a_0\ot [a_1]\ot\cdots\ot [a_n]\longmapsto
1\ot [a_0]\ot [a_1]\ot\cdots\ot [a_n],
\]
and transform it via the just described identification $\hO^n\!A \cong A\ot (A/k)^{\otimes n}$
into the desired differential $\hO^n\!A\ra\hO^{n+1}\!A$. It is a standard fact 
that in this way we obtain a differential graded algebra, and
that this algebra has an appropriate universality property in the category of differential graded
algebras. 

In the commutative setting, K{\"a}hler differential forms
 play the role of universal
objects (see \cite[1.3.7., p.26]{l-jl98}). For a 
commutative algebra $A$, the $A$-module of (first order) K{\"a}hler 
differentials is defined as 
$\hO^1_{A|k}:=\hO^1A/(\hO^1A)^2$. Here we
view $(\hO^1A)^2$ as an ideal in $A\ot A$.
\note{
(where, as above, $\ker m$ is the kernel of the multiplication map 
$m:A\ot A\ra A$). The differential is given as composition of the
quotient map with the differential $d:A\ra\hO^1A$,
$d_u(a)=1\ot a-a\ot 1+(\ker m)^2$.
It has the following universal property:
If $\hG$ is a symmetric $A$-bimodule and $\partial:A\ra\hG$ is a derivation, 
then there
exists a unique $A$-bimodule morphism $\phi:\hO^1_{A|k}\ra\hG$ such that
$\phi\ci d_u=\partial$.
}
If $C^\infty(M)$ is the algebra of smooth functions on 
 a  compact smooth
manifold $M$, and $\hO^1_{dR}(M)$ stands for the module of de Rham
 differential 1-forms, one obtains the following
 sequence of surjective $C^\infty(M)$-bimodule maps:
\beq\label{kadr}
\hO^1C^\infty(M)\lra \hO^1_{C^\infty(M)|\R}\lra \hO^1_{dR}(M).
\eeq
 The second
map 
is
surjective but (in general) not injective.
This lack of injectivity follows from the remarkable fact that for a 
transcendental function $f$, the expression
$\mathrm{d}f-\frac{df}{dx}\mathrm{d}x$ is non-zero as a K\"ahler
form despite
being 
obviously zero in the de Rham calculus, see \cite[p.78]{h-j90} and 
\cite[p.42]{vkl86}.
\footnote{We are grateful to M.~Wodzicki for a discussion clarifying this
point.}

To prove the surjectivity, note first 
that every manifold $M$ admits an atlas such that
every coordinate function can be extended to all of $M$. Indeed,
 consider
a chart $(U,\phi)$ with domain homeomorphic to an open unit ball in 
$\R^n$. 
There is a function $g$ on the open
unit ball whose support is contained in the ball of 
radius~$\frac{3}{4}$, and which 
is the constant function 1 on the ball of radius~$\frac{1}{2}$.
 Multiplying the coordinate functions
with $\phi^*g$, we end up with a new chart with smaller domain 
(homeomorphic to a ball with radius~$\frac{1}{2}$). All coordinate 
functions 
of this chart can be extended to $M$ by~0. Obviously, we can cover $M$ 
with these charts, and in the
compact case we can select finitely many of them.

Now, let 
$\{(U_i,\phi_i)\}_{i\in\{1,\ldots,n\}}$ be such an atlas, let $\phi_i^k$ 
be the $k$-th
coordinate function given by $\phi_i$, let $\widetilde{\phi_i^k}$ be its 
extension to $M$, and let $\{\psi_i\}_{i\in\{1,\ldots,n\}}$ be a smooth 
partition of 
unity subordinate to the covering $\{U_i\}_{i\in\{1,\ldots,n\}}$, i.e.,
$\mathrm{supp}\ \psi_i\subset U_i$, $\sum_{i=1}^n\psi_i=1$.
For any $\ha\in\hO^1_{dR}(M)$, its restriction 
$\ha|_{U_i}=\sum_if^k_i\mathrm{d}\phi_i^k$, for certain
$f^k_i\in C^\infty(U_i)$. Let $\widetilde{\psi_if^k_i}$ denote the
extension by zero of $\psi_i|_{U_i}f^k_i$ to an element of $C^\infty(M)$.
Then 
\[
\ha=
\sum_{i=1}^n\sum_{k=1}^{\dim M}
\widetilde{\psi_if^k_i}\mathrm{d}\widetilde{\phi^k_i}
\]
is the desired presentation of $\ha$ as a finite sum of elements
of the form $f\mathrm{d}g$. This proves the surjectivity of 
the second map in (\ref{kadr}).

\subsection{Projectivity and connections}\label{projco}

Recall that a left $B$-module $P$ is said to be
a {\em projective module} provided it is a direct summand of a
free module. Equivalently, 
$P$ is a projective
left $B$-module if and only if the module structure map $B\otimes
P\ra P$ splits as a left $B$-module map. Note that this last
description is possible since we restricted our considerations to
algebras over a field, so that $B\otimes P$ is a free left
$B$-module, and the above statement is equivalent to say that $P$
is a direct summand of a free module.

For a module $P$, being a direct summand of a free module
means that there exists an idempotent automorphism of a free module whose
image is isomorphic with $P$. Indeed, if $P\st{i}{\ra}B\ot V\st{\pi}{\ra}P$,
$\pi\ci i=\id$, then $(i\ci\pi)^2=i\ci\pi$ and $(i\ci\pi)(B\ot V)\cong P$ via
appropriately restricted $i$ and $\pi$. Similarly, if $e$ is an idempotent
automorphism and $f:P\ra e(B\ot V)$ is an isomorphism, then the composite maps
\[
P\st{f}{\lra}e(B\ot V)\st{\subseteq}{\lra}B\ot V\quad\mbox{and}\quad
B\ot V\st{e}{\lra}e(B\ot V)\st{f^{-1}}{\lra}P
\]
realise $P$ as a direct summand of $B\ot V$. 

\bre\label{sil}
At this point let us recall the following simple property of
idempotent operators, which will be used later on. 
If $E$, $F$ are idempotent linear operators on a vector space $V$
with the same kernel, then $EF=E$. Indeed, since
 $v - F(v)\in \ker F$, we have $v-F(v)\in \ker E$.
Therefore, $ E(v) = E(v - F(v) +F(v)) = E(F(v)). $
\ere

Every projective module can be equivalently characterised
by the existence of a {\em connection}  \cite[Proposition~8.2]{cq95},
 i.e.\ a linear map 
 \[
 \nabla: P\lra \Omega^1\!B\underset{B}{\otimes} P
 \quad\mbox{such that}\quad 
\nabla(bp) = \d b\ten{B} p + b\nabla(p),\quad\forall\; b\in B,\, p\in P.
\]
Indeed, the projectivity of $P$ is tantamount to the existence of 
a left $B$-linear splitting $s$ of the product map $m:B\ot P\ra P$. Taking such a splitting
and viewing elements of
$\Omega^1B\otimes_B P$ in $B\otimes P$, we can define a connection
by the  formula:
\begin{equation}\label{canon.con}
\nabla(p) = 1\otimes p - s(p).
\end{equation}
Thus every projective module admits a connection. Vice versa, if $P$ is a left
$B$-module admitting a connection $\nabla$, then the formula
\[\label{conspl}
s(p)=1\ot p-\nabla(p)
\]
defines a left $B$-linear splitting of the multiplication map, so that $P$ is projective.
The formulas (\ref{canon.con}) and (\ref{conspl}) provide mutually inverse bijections
between the space of splittings and the space of connections. 

This correspondence
between  splittings and connections can be applied to show that every connection
is always obtained from the differential and an idempotent.  Indeed, let $\nabla$ 
be a connection and $s_\nabla$ the splitting of the multiplication map 
associated to $\nabla$ via (\ref{conspl}). Then
\[\label{grass}
\nabla(p)=1\ot p-s_\nabla(p)=((\id\ot m)\ci(\d\ot\id)\ci s_\nabla)(p),\quad\fa\;p\in P.
\]
In general, if $P$ is a direct summand of a free module $B\ot V$ via the maps
$P\st{i}{\ra}B\ot V\st{\pi}{\ra}P$, then the formula
\[\label{grass2}
\nabla(p)=((\id\ot \pi)\ci(\d\ot\id)\ci i)(p),\quad\fa\;p\in P,
\]
determines a connection. Such connections are called {\em Gra\ss mann connections}.

The definition of connection can be easily extended from $P$ to $\hO B\ot_{B}P$.
Indeed, there exists a unique extension of $\nabla$ to a degree-one map 
$\hO P\ot_{B}P\ra \hO P\ot_{B}P$  satisfying the graded Leibniz rule:
\[
\nabla(\ha\xi) =\d \ha\ten{B}\xi +(-1)^n\ha\nabla(\xi),\quad\forall\; \ha\in\hO^n\! B\ten{B}P,\, 
\xi\in \hO B\ten{B}P.
\]
One can immediately verify that $\nabla^2$ is linear over $\hO B$. This is a crucial
property allowing us to view $\nabla^2$ as a matrix with entries in $\hO B$. We call $\nabla^2$
the {\em curvature} of a connection $\nabla$. This concept  is fundamental for
the construction of the Chern character. 

In noncommutative geometry projective modules that are also
finitely generated are of particular interest, since they
correspond to vector bundles.  A
left $B$-module $P$ is finitely generated projective if and only
if there exist a positive integer $n$, an $n\times n$-matrix
$e = (e_i^j)_{i,j =1}^n$ with $e_i^j \in B$ such that
$e^2=e$, and  an isomorphism
of left $B$-modules
$
P\cong B^ne
$.
Here $B^n$ is understood as a space of row vectors $(b_1,\ldots,
b_n)$, and the multiplication represents the standard
multiplication of a row vector by a matrix. 
Any two idempotent matrices $e$ and $\tilde{e}$ determine isomorphic 
finitely generated projective modules
if and only if these two idempotent matrices are stably similar --- i.e. 
by enlarging
$e$ and $\tilde{e}$ by zeroes, there then exists an invertible matrix 
which conjugates one enlarged matrix to the other.
 (cf. \cite[Lemma~1.2.1]{r-j94}).
  
Finitely generated projective modules are  equivalently
characterised by the existence of {\em finite dual
bases}\index{finite dual basis}. Thus $P$ is a finitely generated
projective left $B$-module if and only if there exists a finite
set $\{e_i\in P,\; e^i\in \hom_B(P,B)\}_i$ such that, for all
$p\in P$, $p=\sum_i e^i(p)e_i $. The corresponding idempotent
matrix  $e$ has entries $e_i^j: = e^j(e_i)$, while the
splitting $s: B\ra B\otimes P$ of the product map reads
$s(p) = \sum_i e^i(p)\otimes e_i$.
Much as in \eqref{grass}, the
connection $\nabla$ associated to $s$ via \eqref{canon.con} 
can be expressed in terms of the dual basis
$\{e_i\in P,\; e^i\in\hom_B(P,B)\}_i$ and the corresponding
idempotent  as
\begin{equation}\label{gras}
\nabla(e_i) = \sum_j\d e_i^j \ten{B} e_j.
\end{equation}
Indeed, to check \eqref{gras}, note first that the definition of a
dual bases and the corresponding idempotent imply that $e_i
=\sum_j e_i^je_j$. Hence equation~\eqref{canon.con} yields
\[
\nabla(e_i) = 1\otimes e_i - \sum_je^j(e_i)\otimes e_j =
\sum_j1\ot e_i^j\ten{B}e_j - \sum _j e_i^j\ot 1\ten{B}e_j =
\sum_j\d e_i^j \ten{B} e_j.
\]
Finally, let us observe that a straightforward computation shows that
\[\label{curv}
\nabla^2(e_i)= -\sum_{k,l}\d e_i^k\d e_k^l\ten{B}e_l\,.
\]
Applying $\id\ot_Be^j$ to the right hand side of this equation 
gives the matrix of two-forms: $-(\d e)(\d e)e=-e(\d e)(\d e)$.

\section{Differential geometry of principal and vector bundles}
\setcounter{equation}{0}

Thus far we have analysed purely topological aspects of principal
and associated bundles and recalled some relevant general algebraic
 concepts. Herein, we study them in their traditional context of
differential geometry. Throughout this section, spaces are assumed to
be compact Hausdorff, which removes most of the subtleties discussed
in Section~1. On the other hand, this assumption makes the geometry on
these spaces accessible by the aforementioned general (i.e., not bound
to commutative cases) algebraic tools. We use them to recast the standard
differential-geometric concepts of
associated
vector bundles and connections on principal bundles into a 
language that not only lends itself to noncommutative generalisations,
but also functions naturally therein. 
We end the paper by exemplifying
our general considerations on the Dirac-monopole connection on the
Hopf fibration $S^3\ra S^2$.
 
To be in line with the usual setting of noncommutative geometry,
 we always work
with complex-valued functions and complexified versions of vector fields, differential forms
and other geometric objects. Therefore, we adopt a shorthand notation, 
$C^\infty(M,\C)=C^\infty(M),~\hG^\infty(E,\C)=\hG^\infty(E)$, 
$H_{dR}(M,\C)=H_{dR}(M)$, etc.

\subsection{Associated vector bundles}\label{avb}

The goal of this subsection is to give a description of the module 
of continuous (or smooth)
sections of an associated
vector bundle as an appropriate cotensor product. We first  identify
sections of an associated vector bundle with equivariant maps  
(Lemma~\ref{homma}).
Then we show that the latter are in a natural correspondence 
with elements of a certain cotensor product (Lemma~\ref{eqcot}).

Let $G$ be a compact group acting principally
on the right
 on a compact Hausdorff space $X$. Dualising this right action 
gives a $*$-homomorphism
$\Delta_R:C(X)\ra C(X\times G)\cong C(X)\overline{\ot} C(G)$. 
Here, $\overline{\ot}$ denotes the $C^*$-completed tensor product. 
(Recall that all $C^*$-tensor products are equivalent for commutative
$C^*$-algebras, and correspond to the Cartesian product of underlying
spaces.) Similarly, the group multiplication in $G$
gives rise to a $*$-homomorphism $\hD:C(G)\ra C(G\times G)\cong C(G)
\overline{\ot}C(G)$. By analogy with the algebraic situation,
 we say that $\hD$ is a coproduct on $C(G)$ and $\hD_R$ is
a right coaction of $C(G)$ on $C(X)$ (cf.~\cite{t-m85}).

\bre
The principal map $F^G$ (see (\ref{pmap})) has as its dual 
(pull-back) the canonical map
\beq
{\rm can}:C(X)\overline{\ot}C(X)\lra
C(X)\overline{\ot}C(G),\quad \text{can}(a\ot a')=a\hD_R(a').
\eeq
The $C^*$-algebra $C(X/G)$ can be identified with the subalgebra of $C(X)$ of functions that are constant along every fibre,
i.e.,\ the subalgebra of invariants:
\beq
C(X)^{coC(G)}=\{b\in C(X)~|~\hD_R(b)=b\ot 1\}.
\eeq
Note that $\ker ({\rm can})$ coincides with the closure of the linear span of elements
of the form $ab\ot a'-a\ot ba'$, with $a,a'\in C(X)$ and $b\in C(X)^{coC(G)}$.
This can be proved
using the Stone-Weierstra{\ss}
theorem.
\ere

By  Peter-Weyl theory, the matrix elements
of irreducible unitary representations of $G$  generate a  Hopf 
algebra  $\cO(G)$ which is a dense $*$-subalgebra
of the $C^*$-algebra $C(G)$, and any coaction 
${}_{V}\!\hD:V\ra  C(G)\ot V$ on a finite-dimensional vector space $V$
 takes
values in $ \cO(G)\ot V$. The coproduct of the Hopf algebra $\cO(G)$
is the restriction to $\cO(G)$ of the pull-back of the multiplication
map $G\times G\ra G$.
(Due to Woronowicz's generalisation of Peter-Weyl theory
 \cite{w-sl87}, \cite{w-sl98},
analogous statements are true for representations of 
compact quantum groups.)

For any left coaction of $C(G)$ on $V$, we can
define the cotensor product:
\beq
C(X)\coten{C(G)} V:=\{u\in C(X)\ot V~|~(\hD_R\ot\id)(u)=(\id\ot{}_V\hD)(u)\}.
\eeq
If $V$ is finite dimensional with a basis $\{v_1,...,v_n\}$ and
 $\sum_i p_i\ot v_i\in C(X)\Box_{C(G)} V$, then it follows
that
$\hD_R(p_i)\in C(X)\ot \cO(G)$, i.e., the
$p_i$'s are elements of
\beq\label{polfib}
\cC(X):=\{p\in C(X)~|~\hD_R(p)\in C(X)\ot\cO(G)\}.
\eeq
It is immediate  that $\cC(X)$ is a right $\cO(G)$-comodule algebra, and 
we have
\beq
C(X)\coten{C(G)} V=\cC(X)\coten{\cO(G)}V.
\eeq
This is an equality of $C(X/G)$-modules 
(with the module structure by multiplication in the left
leg of the cotensor product). The algebra $\cC(X)$ is fundamental 
in our approach
to associated vector bundles, 
and should be thought of as the algebra of functions
that are continuous along the base space and polynomial along the fibre.

Next, let
\beq
\hom_G(X,V)=\{f\in C(X,V)\ |\ f(xg)=\rho(g^{-1})f(x)\}.
\eeq
 This is a $C(X/G)$-module
by pointwise multiplication. We can now claim (cf.~Lemma~\ref{cote}):
\ble\label{eqcot}
Let $\rho:G\ra GL(V)$ be a finite-dimensional representation, $\{v_i\}_{i=1}^{n}$ a basis
of $V$ with the dual basis  $\{v^i\}_{i=1}^{n}$, and 
$_V\hD:V\ra \cO(G)\ot V$ the left coaction given by
\[
_V\hD(v_i)=\sum_{j=1}^nv^j(\rho(\;\;\;.\!^{-1})(v_i))\ot v_j\,.
\]
Then $\hom_G(X,V)\cong \cC(X)\Box_{\cO(G)}V$ as $C(X/G)$-modules.
\ele 
\bpf
Observe first that the equivariance of $f\in\hom_G(X,V)$ entails 
\[
(v^j\cc f)(xg)=\sum_{k=1}^nv^j(\rho(g^{-1})(v_k))(v^k\cc f)(x).
\]
Therefore,
$\hD_R (v^j\cc f)(x,g)=\llp\sum_{k=1}^n(v^k\cc f)\ot v^j(\rho(\;\;\;.^{-1})(v_k)\lrp(x,g)$, so that 
 $v^j\cc f\in\cC(X)$.
Furthermore, it is immediate that 
$\sum_{j=1}^n\hD_R(v^j\cc f)\ot v_j=\sum_{j=1}^n(v^j\cc f)\ot {}_V\hD(v_j)$. Thus, we obtain a well-defined
 $C(X/G)$-linear map
$\hom_G(X,V)\ni f\to \sum_{j=1}^n(v^j\cc f)\ot v_j\in \cC(X)\Box_{\cO(G)}V$.
The other way around, if  $\sum_{j=1}^nf^j\ot v_j\in \cC(X)\Box_{\cO(G)}V$, then 
\[
\sum_{j=1}^nf^j(xg)v_j
=\sum_{j=1}^n(\hD_R(f^j))(x,g)v_j
=\sum_{j,k=1}^nf^j(x) v^k(\rho(g^{-1})(v_j))v_k
=\rho(g^{-1})\sum_{j=1}^nf^j(x)v_j\,.
\]
Consequently,
$\sum_{j=1}^nf^jv_j\in \hom_G(X,V)$, and 
$\sum_{j=1}^nf^j\ot v_j\mapsto \sum_{j=1}^nf^jv_j$ 
yields the desired  inverse homomorphism. The $C(X/G)$-linearity
is obvious.
\hfill\epf

Combining this with Lemma \ref{homma}, 
we arrive at:
\bth\label{secot}
Let a compact group $G$ act principally on a compact Hausdorff space $X$, $\rho$ be a representation
of $G$ in a finite dimensional vector space $V$,
and $E$ be the vector bundle associated to the principal bundle 
$(X,\pi,M,G)$ via $\rho$.
Then the space $\hG(E)$ of continuous sections of $E$ is isomorphic as a $C(M)$-module to
the cotensor product
$\cC(X){\Box}_{\cO(G)}V$ defined by the coactions given by the formulas
$
\hD_R(f)(x,g)=f(xg)
$
and
$
_V\hD(v)=v^j(\rho(\;\;\;.\!^{-1})(v))\ot v_j\,.
$
\ethe

The foregoing algebraic formulation of associated bundles is also possible in the
smooth setting. To begin with, there is a well-known analogue of Lemma
\ref{homma} (see, e.g., \cite[16.14, Exercise~8]{d-j74}). Then, it is straightforward
to check that the above considerations are valid also when replacing the
algebra of continuous functions on a compact Hausdorff space by the 
algebra of smooth functions on a compact manifold, and the compact group by a compact
Lie group. 
One finally arrives at an isomorphism of the $C^\infty(M)$-module
$\Gamma^{\infty}(E)$ of smooth sections of a smooth associated vector bundle
with the cotensor product $\cC^{\infty}(X)\Box_{\cO(G)}V$, where
\[\label{smas}
\cC^{\infty}(X):=\{p\in C^\infty(X)~|~\Delta_R(p)\in C^\infty(X)\ot \cO
(G)\}.
\]

\subsection{Connections and gauge transformations}

In this subsection, we consider finite-dimensional
 compact smooth manifolds
without boundary, and assume mappings to be  smooth. We denote
by $T$ the tangent functor, and by $\Lambda^k(M)$ 
the space of differential $k$-forms
on a manifold $M$.
Let $(X,\pi,M,G)$ be a principal bundle. 
The tangent spaces to the fibres form the canonical subbundle 
$\mathop{\rm Ver}\subset TX$ coinciding
with the kernel
of $T\pi$. We shall denote by $R:X\times G\ra X$ the right action, and
by
\[
R_g:X\ni x\mapsto xg\in X,\quad 
R_x:G\ni g\mapsto xg\in X,\quad 
R_{\check{x}}:\pi^{-1}(\pi(x))\ni y\mapsto \check{\tau}(x,y)\in G,
\]
the induced maps. (Recall that $\check{\tau}$ stands for 
the translation map (\ref{classtrans}).)
The adjoint action of $G$ on itself will be denoted by
\[
\mathrm{Ad}^G:g\lra \mathrm{Aut}(G),\quad 
\mathrm{Ad}^G_g(g'):= gg'g^{-1}.
\]

A connection on $(X,\pi,M,G)$ is usually defined as 
 a $G$-equivariant horizontal distribution $\mathop{\rm Hor}$
in $TX$ that is complementary to the distribution $\mathop{\rm Ver}$
 of vertical subspaces
 \cite[p.63]{kn63}. Here equivariance means $TR_g(\mathop{\rm Hor}_x)
\inc \mathop{\rm Hor}_{xg}$, where 
 $x$ and $g$ are arbitrary elements of $X$ and $G$, respectively.
Equivalently, this can be rephrased by saying that a connection on a
principal bundle is an equivariant horizontal lift of tangent vectors
(cf.\ \cite[20.2.2]{d-j71}). The idea of a horizontal lift, combined
with the existence of finite dual bases  for the finitely generated
projective modules of smooth sections of tangent and cotangent bundles,
allows us to prove very useful identifications. We take it for granted
that connections always exist for paracompact (in particular compact)
manifolds.

The horizontal-distribution definition of a connection can be easily reformulated in a dual way, referring to differential
forms instead of tangent vectors.
First, recall that a horizontal form is, by definition, a form that vanishes if at least one of its
arguments is a vertical vector (tangent
to a fibre). We denote  the space of horizontal 1-forms on $X$ by $\Lambda_{\rm hor}^1(X)$.
\bde\label{co*}
A connection on a principal bundle is a $C^\infty(X)$-linear 
idempotent homomorhism
$\hP:\Lambda^1(X)\ra\Lambda^1(X)$ 
such that $\ker\hP=\Lambda_{\rm hor}^1(X)$ and 
$\hP\ci R_g^*=R_g^*\ci\hP$, for all $g\in G$.
\ede
Let us denote the space of connections described in the foregoing definition by ${\cal P}(X)$. Since connections are known to form
an affine space, it is natural ask what kind of a vector space is
generated by differences of elements of ${\cal P}(X)$. To this end,
let us define the following subspace of {\em nilpotent} endomorphisms:
\beq
{\cal N}(X):=\{\,N\in \endo_{C^\infty(X)}(\Lambda^1(X))
~|~ \im\, N\inc\Lambda_{\rm hor}^1(X)\inc\ker N,\; 
R_g^*\cc N=N\cc R_g^*,\;
\forall g\in G\,\}.
\eeq
It is straightforward to verify that ${\cal N}(X)$ is a 
$C^\infty(M)$-module. In particular, it is a vector space,
and we have:
\bpr\label{affco*}
The space of connections ${\cal P}(X)$ is an {\em affine space}
over ${\cal N}(X)$.
\epr
\bpf
Let $N:=\Pi-\Pi'$ be the difference  of two connections. Obviously,
it is an equivariant 
$C^\infty(X)$-module endomorphism with
$\Lambda^1_{\rm hor}(X)\subseteq \ker N$. 
If $E$ is any idempotent linear 
operator on 
$\Lambda^1(X)$ with
$\ker E=\Lambda_{\rm hor}^1(X)$,
then Remark~\ref{sil} implies
$EN=E\Pi-E\Pi'=E-E=0$. Hence 
$\im N\subset \ker E=\Lambda_{\rm hor}^1(X)$.
On the other hand, if $\Pi\in {\cal P}(X)$ and $N\in {\cal N}(X)$,
then again it is clear that
 $\Pi+N$ is an equivariant $C^\infty(X)$-module endomorphism. 
Since $\Pi\ci(\Pi+N)=\Pi$, we obtain 
$\ker(\Pi+N)\inc\ker\Pi=\Lambda^1_{\mathrm{hor}}(X)$.
The reverse inclusion is immediate, so that 
$\ker(\Pi+N)=\Lambda^1_{\mathrm{hor}}(X)$.
Finally, as $\im(\id-\Pi)=\ker\Pi=\Lambda^1_{\mathrm{hor}}(X)$,
we have $N\ci(\id-\Pi)=0$. Consequently,
 $(\Pi+N)^2=\Pi+N$.
\hfill\epf

Traditionally, connection forms are defined as appropriate 
differential 1-forms on $X$ with values in the Lie algebra of $G$.
The following is a version of connection forms avoiding the use of 
tangent vectors and Lie algebras. 
\bde\label{coform}
A {\em connection form} is a linear map 
${\omega}:C^\infty(G)\ra\Lambda^1(X)$
with the properties:
\begin{itemize}
\item[(i)]
${\omega}(h_1h_2)={\omega}(h_1)h_2(e)+
h_1(e){\omega}(h_2)$,\quad $\fa h_1,h_2\in C^\infty(G)$,\quad $e\in G$ 
the neutral element;
\item[(ii)]
$\mathrm{ev}_e\ci R_x^*\ci\omega=\mathrm{ev}_e\ci \d$,
\quad $\fa x\in X$,\quad $\mathrm{ev}_e: 
\hG^\infty(T^*G)\ni\ha\mapsto\mathrm{ev}_e(\ha):=\ha_e\in T^*_eG$;
\item[(iii)]
$R_g^*\ci\ho=\ho\ci(\mathrm{Ad}_{g^{-1}}^G)^*$,\quad $\fa g\in G$.
\end{itemize}
\ede
\noindent It is straightforward to show:
\bpr\label{cofot}
Let $\tilde{\omega}$ be a traditionally defined connection 1-form,
$Y$ a vector field on $X$, and $h$ a smooth function on $G$.
Then the formula $({\omega}(h))(Y)=(\tilde{\omega}(Y))(h)$ implements
a bijective correspondence between connection forms and traditionally
defined connection forms.
\epr
We shall denote the space of connection forms 
by $\cF(X)$.
There is an obvious affine structure on  $\cF(X)$. 
The difference of any
two such connection
forms satisfies the defining conditions (i) and (iii),
 and annihilates vertical 
vectors due to the defining condition (ii).

Let us now pass to the third description of the concept of a connection:
covariant differentiation. This time there is no need to reformulate
its traditional definition. What we need to show instead is that it is
an equivalent formulation of the connection understood as an appropriate
equivariant distribution in the cotangent bundle~$T^*X$. Recall that the
equivalence of connections and traditionally defined 
connection forms is considered completely
standard and is proven in many books, e.g., in \cite{kn63}. Therefore,
since the traditional definition of connection forms is equivalent
to ours, we will take  for granted the equivalence of connections and
connection forms, and focus on the equivalence of connections
 and covariant differentiations. The latter are axiomatically defined
in the following manner:
\bde
An (exterior) {\em covariant differentiation} is a linear map
 $D:C^\infty(X)\ra \Lambda_{\rm hor}^1(X)$ with the following properties:
\begin{itemize}
\item[(i)]
$D(fg)=D(f)g+fD(g)$,\quad $\fa f,g\in C^\infty(X)$;
\item[(ii)]
$\d b=Db$,\quad $\fa b\in \pi^*(C^\infty(M))$;
\item[(iii)]
$R_g^*\ci D=D\ci R_g^*$,\quad $\fa g\in G$.
\end{itemize}
\ede

To proceed further, we need the following technical result:
\ble\label{hofo}
The space $\Lambda_{\rm hor}^1(X)$ of horizontal 1-forms coincides with 
$C^\infty(X)\pi^*(\Lambda^1(M))$.
\ele
\bpf
The inclusion 
$C^\infty(X)\pi^*(\Lambda^1(M))\subseteq\Lambda_{\rm hor}^1(X)$ is 
immediate. 
To prove the reverse  inclusion, let us
choose a connection. It defines the horizontal projection,
vertical projection and horizontal lift. We denote these maps 
by $\mathrm{hor}$, 
$\mathrm{ver}$ and $\widehat{\phantom{Y}}$, respectively.   
Let $\{e_\mu,e^\nu\}_{\mu,\nu}$ be  finite dual bases of the finitely
 generated projective $C^\infty(M)$-modules $\Gamma^\infty(TM)$
and $\Gamma^\infty(T^*M)=\Lambda^1(M)$. Then, for any horizontal
form $\ha\in\Lambda_{\rm hor}^1(X)$ and tangent vector $Y_x\in T_xX$,
we obtain:
\begin{align}
\left(\mbox{$\sum_\mu$}\;\ha(\widehat{e_\mu})\;\pi^*(e^\mu)\right)(Y_x)
&=
\mbox{$\sum_\mu$}\;
\ha\left(\left(\widehat{(e_\mu)_{\pi(x)}}\right)_x\right)
\;e^\mu(\pi_*Y_x)
\nonumber\\ &=
\ha\left(\left(\mbox{$\sum_\mu$}\;
(e_\mu)_{\pi(x)}\;e^\mu(T\pi(Y_x))
\right)^{\!\!\!\widehat{\phantom{yy}}}_{\;x}\right)
\nonumber\\ &=
\ha\left(\widehat{T\pi(Y_x)}_{x}\right)
\nonumber\\ &=
\ha\left(\mathrm{hor}(Y_x)\right)
\nonumber\\ &=
\ha\llp\mathrm{hor}(Y_x)+\mathrm{ver}(Y_x)\lrp
\nonumber\\ &=
\ha(Y_x)\,.
\end{align}
Hence $\ha=\sum_\mu\;\ha(\widehat{e_\mu})\;\pi^*(e^\mu)\in 
C^\infty(X)\pi^*(\Lambda^1(M))$.
\hfill\epf

We are now ready for:
\bpr\label{coder}(cf.\ \cite[p.254--255]{ghv73})
The formulas $j_{pd}(\Pi)=(\id-\Pi)\ci \d$ and
$j_{dp}(D)(\d f)=(\d-D)(f)$ define mutually inverse
bijections between the space of connections and the space
of covariant differentiations.
\epr
\bpf
Let  $\Pi$ be a connection in the sense of Definition~\ref{co*}. 
Since it is an idempotent whose kernel is the space of
horizontal 1-forms, we have 
$\im(\id-\hP)=\ker\hP=\Lambda_{\rm hor}^1(X)$. Therefore, 
$D_\Pi:=(\id-\Pi)\ci \d$ is a linear map from $C^\infty(X)$ into
$\Lambda^1_{\rm hor}(X)$. Define $D_\Pi:=(\id-\hP)\ci d$.
The Leibniz rule  for $D_\Pi$ comes from the Leibniz rule
for $\d$ and $C^\infty(X)$-linearity of $\Pi$. The form $\d(f\ci\pi)$
is evidently horizontal, so that $D(f\ci\pi)=\d(f\ci\pi)$.
Finally, the equivariance  of $D_\Pi$ is obvious from the 
equivariance properties of $\Pi$ and
$\d$. Thus $D_\Pi$ enjoys all the necessary properties.

Now, let  $D$ be a covariant differentiation, $Y$ a vector field on 
$X$, and $f$ a smooth function on $X$. Define $D_Y(f):=(D(f))(Y)$.
The Leibniz rule for $D$ implies that $D_Y$ is a derivation of 
$C^\infty(X)$. Therefore, it is a vector field on $X$, and we can
define $\Pi_D(\ha):=\ha(Y-D_Y)$, $\ha\in\Lambda^1(X)$. 
The thus defined $\Pi_D$ is evidently
 $C^\infty(X)$-linear. Hence, as any $\ha\in\Lambda^1(X)$ is a finite
$C^\infty(X)$-linear combination of exact 1-forms (see the end of
Section~\ref{uda}), it is uniquely 
determined by its values on exact 1-forms, where it coincides with
$j_{dp}(D)$.  Thus $j_{dp}(D)=\Pi_D$, and  the formula for $j_{dp}(D)$ defines a
 $C^\infty(X)$-linear endomorphism of $\Lambda^1(X)$. It is immediate
that $\Pi_D$ is equivariant. Furthermore, if  $\Pi_D(\ha)=0$, then
$\ha(Y)=\ha(D_Y)$. Presenting \ha\ as a finite sum $\sum_{k}f_k\d g_k$ 
yields $\ha(Y)=\sum_{k}f_k(Dg_k)(Y)$. 
Since $\im(D)\inc \Lambda^1_{\mathrm{hor}}(X)$, we obtain
$\ha(Y)=0$ if $Y$ is
vertical. Hence $\ker\hP_D\inc\Lambda_{\rm hor}^1(X)$.
The reverse inclusion follows from Lemma~\ref{hofo}, which allows
us to choose the functions $g_k$ from $\pi^*(C^\infty(M))$. Finally,
the just proven inclusion $\ker\hP_D\supseteq\Lambda_{\rm hor}^1(X)$
together with $\Lambda^1_{\mathrm{hor}}(X)\supseteq\im(D)$
implies that $(\hP_D)^2=\hP_D$.

To end the proof, we need to verify that $D_{\Pi_D}=D$ and 
$\Pi_{D_\Pi}=\Pi$. The former is immediate, and the latter is
obvious when remembering that any 1-form is a finite sum
$\sum_{k}f_k\d g_k$.
\hfill\epf

Thus we have shown 
that connections on a principal bundle can be equivalently 
defined as exterior
covariant differentiations.
(Observe that the equation $D_\Pi =(\id-\Pi)\ci\d$ of the foregoing 
proposition is just the usual definition
$D=\d\ci \hor$.)
We denote the space of all exterior covariant differentiations by 
${\cal D}(X)$. It has an obvious
affine structure.
The difference of any
two exterior covariant differentiations is a linear map
$C^\infty(X)\ra\Lambda_{\rm hor}^1(X)$ that
is equivariant, satisfies the Leibniz rule and annihilates 
$\pi^*(C^\infty(M))$.

\note{
More generally,
a connection on a principal bundle defines an exterior covariant differentiation
also on vector valued forms, i.e. elements of $\bigwedge^k(X)\ot V$,
$V$ a finite dimensional vector space. One defines
\beq
D\alpha(Y_1,\ldots,Y_k):=d\alpha(\hor Y_1,\ldots.\hor Y_k).
\eeq
The exterior covariant differentiation $D$ is a derivation: Let
$V_1,V_2,V_3$ be vector spaces, and $\eta:V_1\times V_2\ra V_3$
be a bilinear map. Then an exterior product $\land_\eta:
(\bigwedge^k(X)\ot V_1)\times (\bigwedge^l(X)\ot V_2)\ra
(\bigwedge^{k+l}(X)\ot V_3)$ is defined for every $k,l$ by
\beq
\alpha\land_\eta\beta(Y_1,\ldots,Y_{k+l})=
\frac{1}{k!l!}\sum_{\pi\in{\cal S}_{k+l}}\epsilon(\pi)
\eta(\alpha(Y_{\pi(1)},\ldots,Y_{\pi(k)}),
\beta(Y_{\pi(k+1)},\ldots,Y_{\pi(k+l)})),
\eeq
where ${\cal S}_{k+l}$ is the group of permutations of $k+l$ elements
and $\epsilon$ is the sign of a permutation. Then
\beq\label{der}
D(\alpha\land_\eta\beta)=D\alpha\land_\eta\beta+(-1)^k\alpha
\land_\eta D\beta,
\eeq
as a consequence of the derivation property of $d$.
If $\sigma:G\ra GL(V)$ is a linear representation of $G$ on a finite dimensional
vector space $V$,
a $k$-form $\alpha$ on $X$ with values in $V$ is said to
be $\sigma$-equivariant (or: of type $\sigma$), if
\beq\label{equi}
R_g^*\alpha=\sigma(g^{-1})\alpha.
\eeq
A connection form is a $\frak{g}$-valued 1-form of type $Ad$.
The exterior covariant differentiation defined by a connection preserves
horizontal forms and forms of type $\sigma$. If $\alpha$ is both
horizontal and of type $\sigma$, then
\beq
D\alpha=d\alpha + \tilde{\omega}\land_\sigma\alpha.
\eeq
In the last formula, $d$ is the usual exterior differential (of vector valued forms),
$\tilde{\omega}$ is the connection form, and $\land_\sigma$ is given by
\beq
(\tilde{\omega}\land_\sigma\alpha)_x(Y_1,\ldots,Y_{k+1})=
\sum_{i=1}^{k+1}(-1)^{i+1}\sigma_*(\tilde{\omega}_x(Y_i))
(\alpha_x(Y_0,\ldots,\widehat{Y_i},\ldots,Y_{k+1})),
\eeq
where $\sigma_*:\frak{g}\ra End(V)$ is the derived homomorphism of
the representation $\sigma:G\ra GL(V)$, and $Y_i\in T_xX$.
The curvature form $\Omega$ of a connection is the covariant exterior
differential of the connection
form $\tilde{\omega}$, i.e.,
\beq
\Omega=D\tilde{\omega}.
\eeq
The curvature form is a horizontal 2-form with values in $\frak{g}$
of type $Ad$ (cf.\ (\ref{zuad})).
For a 1-form $\tilde{\omega}$ on $X$ with values in $\frak{g}$, define a 2-form 
$[\tilde{\omega},\tilde{\omega}]$ by
\beq
[\tilde{\omega},\tilde{\omega}](Y_1,Y_2)=[\tilde{\omega}(Y_1),\tilde{\omega}(Y_2)]-
[\tilde{\omega}(Y_2),\tilde{\omega}(Y_1)].
\eeq
Then the curvature form satisfies the structure equation
\beq\label{stru}
\Omega=d\tilde{\omega} + \frac{1}{2}[\tilde{\omega},\tilde{\omega}]
\eeq
and the Bianchi identity
\beq\label{bia}
D\Omega=0.
\eeq
For any horizontal form $\alpha$ of type $\sigma$,
\beq
D^2\alpha=\Omega\land_\sigma\alpha.
\eeq
For more details about covariant differentiations we refer to \cite[20.3]{d-j71},
\cite[Chapter II, \S5]{kn63} and \cite[Kapitel II,\S9]{sw72}.
}

We end this section by considering the behaviour of connections under
gauge transformations. A  gauge 
transformation of a principal bundle $(X,\pi,M,G)$
is a vertical automorphism of the bundle, i.e., a $G$-equivariant
fibre-preserving diffeomorphism
$\gamma:X\ra X$. 
Gauge transformations form a group with group multiplication given by
the composition of maps. We denote this group  by ${\cal G}(X)$,
and define its  left action
on the spaces ${\cal P}(X)$, ${\cal F}(X)$ and
 ${\cal D}(X)$ of connections,
connection forms
and covariant differentiations, respectively,  by
\[\label{pig}
\Pi^\gamma:=(\hg^{-1})^*\ci\Pi\ci\hg^*,\quad
{\ho}^\gamma:=(\hg^{-1})^*\ci\ho,\quad
D^\gamma:=(\hg^{-1})^*\ci D\ci\hg^*\,.
\]
This action commutes with our identifications of connections,
connection forms
and covariant differentiations. These identifications respect also
the affine structure.
More precisely, putting together the results of this section, one
 can verify the following:
\bth
Let ${\cal G}(X)$ be the group of gauge transformations 
acting on the 
 spaces ${\cal P}(X)$, ${\cal F}(X)$ and ${\cal D}(X)$ of connections, 
connection forms and
exterior covariant differentiations according to (\ref{pig}).
The following formulas define ${\cal G}(X)$-equivariant  and affine
 isomorphisms:
\begin{align}
&
j_{dp}: {\cal D}(X)\ni D\longmapsto \Pi^D\in {\cal P}(X),\quad 
\Pi^D(\d f):=(\d-D)f;
\\
&
j_{pf}: {\cal P}(X)\ni \Pi\longmapsto \ho^\Pi\in {\cal F}(X),\quad 
(\ho^\Pi(h))(Y_x):=(\mathrm{ver}^\Pi Y_x) 
(R_{\check{x}}^{\,*}\,h);
\\
&
j_{fd}: {\cal F}(X)\ni \ho\longmapsto D^\omega\in {\cal D}(X),\quad 
(D^\omega(f))(Y_x):=\llp\d f-\ho(R_x^{\,*}\,f)\lrp(Y_x)\,.\label{jfd}
\end{align}
Here $f\in C^\infty(X)$, $h\in C^\infty(G)$, $Y_x\in T_xX$, and
$\mathrm{ver}^\Pi Y_x$ is defined by 
$(\mathrm{ver}^\Pi Y_x)(f):=(\Pi(\d f))(Y_x)$.
\ethe
\note{
It is straightforward that the bijective mappings $j_{dp}:{\cal D}(X)\ra {\cal P}(X)$
and $j_{pd}:{\cal P}(X)\ra{\cal D}(X)$,
given in the proof of Proposition~\ref{coder}, are
compatible with the right actions of ${\cal G}(X)$ and with the affine structures. Let us give also
explicit formulas for bijective mappings $j_{pf}:{\cal P}(X)\ra{\cal F}(X)$ and
$j_{fd}:{\cal F}(X)\ra{\cal D}(X)$ being compatible with all these structures.
For $h\in C^\infty(G)$ and $x\in X$, define $ h_x\in C^\infty(X_{\pi(x)})$
($X_{\pi(x)}$ the
fibre of $X$ through $x$) by $ h_x(y)=f(g)$, for $y=xg$. Then, 
for given $\Pi\in{\cal P}(X)$,
define $j_{pf}(\Pi)\in{\cal F}(X)$ by
\beq\label{jpf}
(j_{pf}(\Pi)(h))_x=\Pi_x(d(h_x)),
\eeq
 where $d$ is the differential along
the fibre. Let now $f\in C^\infty(X)$, and define for $x\in X$
the function $f_x\in C^\infty(G)$ by $f_x(g)=f(xg)$. 
For given ${\ho}\in {\cal F}(X)$, define
$j_{fp}({{\ho}})\in{\cal P}(X)$ by
\beq\label{jfp}
{j_{fp}({\ho})}_x(df_x)=({\ho}(f_x))_x.
\eeq
Composing with $j_{pd}$, we obtain the
desired $j_{fd}:{\cal F}(X)\ra {\cal D}(X)$, sending ${\ho}$ to
$j_{fd}({\ho})$ given by
\beq\label{jfd}
(j_{fd}({\ho})f)_x=(df)_x-({\ho}(f_x))_x.
\eeq
}

\subsection{Covariant derivatives and the Chern character}

Let us 
 recall that the correspondence between equivariant functions on a 
principal bundle and
sections of an associated vector bundle
 extends to a correspondence between
horizontal equivariant forms on the principal bundle and form-valued 
sections
of the associated vector bundle. Since covariant differentiations
on the principal bundle
  preserve equivariant forms,  by this correspondence, they give
rise to  covariant derivatives (connections) 
on the associated vector bundle.
Recall also that a connection on a vector bundle $E$ over the base $M$ 
(cf.\ Section~\ref{projco}) can be defined as a linear map 
\[
\nabla:\Gamma^\infty(E)\ra \Lambda^1(M)
\!\!\!\!\underset{C^\infty(M)}{\ot}\!\!\!\!
\Gamma^\infty(E)
\]
satisfying the Leibniz rule
\beq\label{connec}
\nabla(fs)=f\nabla(s)+df\!\!\!\!\underset{C^\infty(M)}{\ot}\!\!\!\!s,
\quad f\in C^\infty(M),\quad s\in\Gamma^\infty(E)\,.
\eeq
In the first part of this section, we will analyse covariant derivatives 
 from the point of view adopted in Section~\ref{avb},
i.e., identifying the module of smooth sections
of an associated bundle with a suitable cotensor product. 

To identify $\Lambda^1_{\rm hor}(X)$ with
$\Lambda^1(M)\ot_{C^\infty(M)}C^\infty(X)$,
note first that any connection, by providing us with a horizontal
lift, defines a push-down:
\[
\pi_*:\Lambda^1_{\rm hor}(X)\lra
\hom_{C^\infty(M)}(\Gamma^\infty(TM),C^\infty(X)),\quad
(\pi_*(\ha))(Y):=\ha(\widehat{Y}),\quad Y\in \Gamma^\infty(TM)\,.
\]
It follows from Lemma~\ref{hofo} that its inverse is given by the
pull-back $\pi^*$. Next, since $\Gamma^\infty(TM)$ is finitely generated
projective, we have the
canonical isomorphism  
\beq
\hom_{C^\infty(M)}(\Gamma^\infty(TM),C^\infty(X))\lra 
\Lambda^1(M)\!\!\!\!\underset{C^\infty(M)}{\ot}\!\!\!\!C^\infty(X)\,.
\eeq
With the help of a finite dual bases $\{f_\mu, f^\nu\}_{\mu,\nu}$ 
of $\Gamma^\infty(TM)$ and $\Lambda^1(M)$, 
the composition of these two isomorphisms can be explicitly
written as
\beq\label{liot}
\Phi:\Lambda^1_{\rm hor}(X)\lra
\Lambda^1(M)\!\!\!\!\underset{C^\infty(M)}{\ot}\!\!\!\!C^\infty(X),
\quad
\Phi(\ha)= \sum_\nu f^\nu
\!\!\!\!\underset{C^\infty(M)}{\ot}\!\!\!\!
\ha(\widehat{f_\nu}).
\eeq
The inverse of this isomorphism is simply the multiplication. 
If $\ha=f\hb$ for $f\in C^\infty(X)$ and $\hb\in\pi^*(\Lambda^1(M))$,
 then
$\Phi(f\hb)=\hb\ot_{C^\infty(M)}f$.


Furthermore, observe that the right action of $G$ on $X$ yields a 
right coaction $\hD_R$ of $C^\infty(G)$ on $C^\infty(X)$, and
$\hD^1_R:=\id\ot_{C^\infty(M)}\hD_R$ defines a 
right coaction on $\Lambda^1(M)\ot_{C^\infty(M)}C^\infty(X)$. 
 As in \eqref{smas}, let 
$\cC^\infty(X):=
\{p\in C^\infty(X)~|~\hD_R(p)\in C^\infty(X)\ot \cO(G)\}$.
At this point, we need the following technical lemma:
\ble
If $D:C^\infty(X)\ra \Lambda^1_{\mathrm{hor}}(X)$  is a covariant
 differentiation, then
\[
(\Phi\ci D)(\cC^\infty(X))\subset \Lambda^1(M)
\!\!\!\!\underset{C^\infty(M)}{\ot}\!\!\!\!
\cC^\infty(X)\,.
\]
\ele
\bpf
For brevity, we put $\breve{D}:=\Phi\ci D$. First we need to prove 
\beq\label{pds}
(\hD_R^1\ci\breve{D})(f)\in \Lambda^1(M)
\!\!\!\!\underset{C^\infty(M)}{\ot}\!\!\!\!
C^\infty(X)\ot\cO(G),\quad\fa f\in\cC^\infty(X)\,.
\eeq
Writing  $\hD_R(f)$ as a finite sum
$\sum_kf_k\ot h^k\in C^\infty(X)\ot\cO(G)$, 
 for a fixed $g\in G$, we can compute:
\begin{align}
\left(\hD_R^1(\breve{D}(f))\right)(g)
&=
\sum_\nu f^\nu
\!\!\!\underset{C^\infty(M)}{\ot}\!\!\!
\left(\hD_R\left((Df)(\widehat{f_\nu})\right)\right)(g)
\nonumber\\ &=
\sum_\nu f^\nu
\!\!\!\underset{C^\infty(M)}{\ot}\!\!\!
R_g^*\left((Df)(\widehat{f_\nu})\right)
\nonumber\\ &=
R_g^*\llp\breve{D}(f)\lrp
\nonumber\\ &=
\breve{D}\left(R_g^*f\right)
\nonumber\\ &=
\breve{D}\llp(\hD_Rf)(g)\lrp
\nonumber\\ &=
\left(\mbox{$\sum_k$}\breve{D}(f_k)\ot h^k\right)(g)\,.
\end{align}
Here we used the equivariance of $D$ 
and the fact that $\Phi$  commutes with $R_g^*$. Since this 
computation works for any
$g\in G$, the inclusion (\ref{pds}) is proven.

Next, consider the following identities:
\begin{align}\label{dutrick}
&
\Lambda^1(M)
\!\!\!\underset{C^\infty(M)}{\ot}\!\!\!
C^\infty(X)
\ni\sum_\nu f^\nu
\!\!\!\underset{C^\infty(M)}{\ot}\!\!\!
F_\nu
=
\sum_{\mu,\nu}f^\mu f^\nu(f_\mu)
\!\!\!\underset{C^\infty(M)}{\ot}\!\!\!
F_\nu
=
\sum_{\mu}f^\mu 
\!\!\!\underset{C^\infty(M)}{\ot}\!\!\!
\mbox{$\sum_\nu$}f^\nu(f_\mu)F_\nu\,,
\\ &
\hD_R\left(\mbox{$\sum_{\nu}$}f^\nu(f_\mu)F_\nu\right)
=\sum_\nu f^\nu(f_\mu) \hD_R(F_\nu)=
\llp(m\ot\id)\ci(f_\mu
\!\!\!\underset{C^\infty(M)}{\ot}\!\!\!
\id\ot\id)\ci(\id
\!\!\!\underset{C^\infty(M)}{\ot}\!\!\!
\hD_R)\lrp(\mbox{$\sum_\nu$}f^\nu
\!\!\!\underset{C^\infty(M)}{\ot}\!\!\!
F_\nu)\,,
\nonumber
\end{align}
where $m$ is the multiplication map. Now it follows that
\[
\hD_R^1\left(\mbox{$\sum_\nu$}f^\nu
\!\!\!\underset{C^\infty(M)}{\ot}\!\!\!
F_\nu\right)\in \Lambda^1(M)
\!\!\!\!\underset{C^\infty(M)}{\ot}\!\!\!\!
C^\infty(X)\ot\cO(G)
\nonumber\]\[
\Downarrow
\]\[
\sum_{\nu}f^\nu(f_\mu)F_\nu
\in\cC^\infty(X),\quad\fa \mu\,.
\nonumber\]
 Combining this implication with (\ref{pds}), we obtain the assertion
of this lemma.
\hfill\epf 

We are now ready to associate a covariant derivative to a covariant
differentiation. Restricting $\breve{D}$ to $\cC^\infty(X)$,
we obtain a linear map
\beq\label{nablaco}
\nabla:=(\breve{D}|_{\cC^\infty(X)}\ot\id):\;\cC^\infty(X)\coten{\cO(G)}V \lra
\Lambda^1(M)
\!\!\!\!\underset{C^\infty(M)}{\ot}\!\!\!\!
\cC^\infty(X)\coten{\cO(G)}V.
\eeq
It is straightforward to check that the thus defined
$\nabla$ satisfies the Leibniz rule. 

\note{
Our next aim is to relate the connection $\nabla$ to a Gra{\ss}mann 
connection.
Any $s\in \cC^\infty(X)\Box_{\cO(G)}V$ is of the form
$s=\sum_A s^A\ot v_A$, where the $\{v_A\}_A$ is a basis of the 
finite-dimensional vector space $V$.
Let  $\{e_k,e^l\}_{k,l}$ be finite dual bases of 
$\cC^\infty(X)\Box_{\cO(G)}V$ and its dual. Put
$e_i=:\sum_A s^A_i\ot v_A$.
Using the horizontal lift with respect to the connection 
corresponding to $D$
in the definition of the isomorphism $\Phi$, we obtain 
\beq
\nabla(e_i)=\sum_{\nu,A}f^\nu
\!\!\!\!\underset{C^\infty(M)}{\ot}\!\!\!\!
 \widehat{f_\nu}(s^A_i)\ot v_A\,.
\eeq
For any $k$, since 
$ \sum_A\hat{f}_\nu(s^A_k)\ot v_A\in \cC^\infty(X)\Box_{\cO(G)}V$,
there exist functions $T^j_{\nu k}\in C^\infty(M)$ such that
\beq\label{tjnk}
\sum_A\hat{f}_\nu(s^A_k)\ot v_A=
\sum_jT^j_{\nu k}e_j\,.
\eeq
Consequently, we can write
\beq\label{nabt}
\nabla(e_k)=\sum_{\nu,j}f^\nu
\!\!\!\underset{C^\infty(M)}{\ot}\!\!\! 
T^j_{\nu k}e_j=
\sum_{\nu,j}f^\nu T^j_{\nu k}
\!\!\!\underset{C^\infty(M)}{\ot}\!\!\!
e_j.
\eeq
On the other hand, the  Gra{\ss}mann connection 
corresponding to the idempotent $e:=(e^k_l):=(e^k(e_l))$
is given by
\beq\label{grasp}
\nabla^e(e_k):=\sum_j\d e^j_k
\!\!\!\underset{C^\infty(M)}{\ot}\!\!\!
 e_j.
\eeq
\bpr
The connection $\nabla$ on $\cC^\infty(X)\Box_{\cO(G)}V$ that is
associated to an exterior covariant differentiation
$D$ coincides with the Gra{\ss}mann connection $\nabla^e$
corresponding to the idempotent $e$
{\em if and only if} the functions $T^j_{\nu k}$ from
equation (\ref{tjnk}) satisfy
\beq\label{tgra}
\sum_{j,\nu}f^\nu T^j_{\nu i}e^l_j=\sum_kde^k_ie^l_k.
\eeq
\epr
\bpf
Assuming $\nabla=\nabla^e$, acting with $\id\ot e^l$ on the right hand sides
of (\ref{nabt}) and (\ref{grasp}) and multiplying out gives (\ref{tgra}).
If (\ref{tgra}) is satisfied, one uses $\sum_le^l_ke_l=e_k$ and the fact that
one can move elements of $C^\infty(M)$ over the tensor product to obtain
equality of the right hand sides of (\ref{nabt}) and (\ref{grasp}).\hfill
\epf
}

Chern classes and the Chern character are fundamental topological
invariants of vector bundles. Our next aim is to recall some  of
these constructions.
Let $(E,\pi,M)$ be a complex vector
bundle. Then, for each $i\geq 0$, the $i$-th Chern class
$c_i(E)\in H^{2i}_{dR}(M)$ 
can be defined by a certain set of axioms
(\cite[pp.305--312]{kn69}, \cite[Chapter 17]{h-d94}). It is
crucial for computations that these classes can be realized in
terms of the curvatures of connections. This is the celebrated
Chern-Weil construction. For a treatment in an algebraic setting 
we refer to
\cite[Section~8.1]{l-jl98}, for the classical setting to \cite{g-pb95}
and \cite{ms74}.
\note{
Let $\nabla$ be a connection on the vector bundle $E$.
If $(e_1,\ldots,e_k)$ is a local
frame of $E$ defined over some open $U\subset M$, then 
\beq 
\nabla e_i=\ho_i^j\underset{C^\infty(U)}{\ot} e_j 
\eeq 
defines a
$k\times k$-matrix $(\ho_i^j)$ of (locally defined) differential
1-forms, called the connection form corresponding to the local
frame. 
 Any such locally defined matrix
of 1-forms locally defines a connection. Under a change
${e'}_k=h_k^le_l$ of the local frame the connection form
transforms as (written in matrix notation) 
\beq 
{\ho'} =h\ho
h^{-1}+dh h^{-1}. 
\eeq
}

Many of the considerations below follow verbatim their
counterparts in the universal
differential algebra setting of Section~\ref{projco}.
 The graded Leibniz rule 
 allows us to extended a
 covariant derivative $\nabla$   to a linear map
\[
\nabla:\Lambda^*(M)
\!\!\!\underset{C^\infty(M)}{\ot}\!\!\!
\hG^\infty(E)\lra\Lambda^{*+1}(M)
\!\!\!\underset{C^\infty(M)}{\ot}\!\!\!
\hG^\infty(E)\,,\quad
\nabla^1(\ha
\!\!\!\underset{C^\infty(M)}{\ot}\!\!\!
\xi):=\d\ha
\!\!\!\underset{C^\infty(M)}{\ot}\!\!\!
\xi +(-1)^{\mathrm{deg}\alpha} \ha\wedge\nabla(\xi)
\,.
\]
The composition $\nabla^2:=\nabla\ci\nabla:\hG^\infty(E)\ra
\Lambda^2(M)\ot_{C^\infty(M)}\hG^\infty(E)$ is called the {\em curvature}
 of a
connection $\nabla$. One easily checks that 
 $\nabla^2$ is an endomorphism of $\Lambda(M)
\ot_{C^\infty(M)}\hG^\infty(E)$ over the differential graded algebra
$\Lambda(M)$ of de Rham differential forms.
In particular, the curvature can  be considered 
as a 2-form on $M$ with values
in the endomorphisms of $\hG^\infty(E)$.  
\note{
and can be represented with respect
to a local frame by a matrix of 2-forms, 
\beq 
\hO e_i=\hO_i^j\underset{C^\infty(U)}{\ot}
e_j, 
\eeq 
transforming under a change of the local frame as 
\beq
\hO'=h\hO h^{-1}. 
\eeq 
}

Since $\hG^\infty(E)$ is a finitely generated projective 
$C^\infty(M)$-module, we can choose finite dual bases
$\{e_k,e^l\}_{k,l}$. The  Gra{\ss}mann connection 
corresponding to the idempotent $e:=(e^k_l):=(e^k(e_l))$
is given by
\beq\label{grasp}
\nabla^e(e_k):=\sum_j\d e^j_k
\!\!\underset{C^\infty(M)}{\ot}\!\!
 e_j\,.
\eeq
For a general $\xi=\sum_k\xi^ke_k\in \hG^\infty(E)$, it follows from the Leibniz rule
that
$ 
\nabla^e(\xi)=\sum_{j,k}\d \xi^k e^j_k{\ot}_{C^\infty(M)}e_j 
$.
This is sometimes written in the 
matrix notation as $\nabla^e(\xi)=(\d \xi) e$.
One finds that the
corresponding curvature is given by
\beq\label{ogras}
(\nabla^e)^2\left(\mbox{$\sum_k$}\xi^ke_k\right)
=-\sum_{j,k,l,m}\xi^ke^j_k\d e^l_j\d e^m_l
\!\!\underset{C^\infty(M)}{\ot}\!\!
e_m\,.
\eeq
In the matrix notation one would write it 
$(\nabla^e)^2(\xi)=-\xi e\,\d e\,\d e$.
Note that the minus sign arises from the convention to consider $\hG^\infty(E)$ as a left
$C^\infty(M)$-module.
\note{
Let now $P$ be an invariant $\C$-valued polynomial on
$\frak{g}:=\frak{gl}(\C,k)=End(\C^k)$, i.e., $P(hAh^{-1})=P(A)$
for $A\in\frak{g}$. For any complex vector bundle $E$, such a polynomial
$P$ defines a map $\Lambda^*(M)\ot_{C^\infty(M)}\hom(E,E)\ra \Lambda^*(M\ot\C)$.
In particular,
}

Now let us view the curvature  of a connection $\nabla$
 as an element of 
$\Lambda^2(M)\ot_{C^\infty(M)}\mathrm{End}_{C^\infty(M)}(\hG^\infty(E))$.
%
Since the trace of an endomorphism of a finitely generated projective
 module 
over a commutative ring is well defined and linear over this ring,
we can form
\bea
ch(\nabla^2)&=&\Tr(e^{i\nabla^2/2\pi})=k+ch_1(\nabla^2)+\ldots+ch_k
(\nabla^2)+\ldots,\\
ch_k(\nabla^2)&=&(i/2\pi)^k\Tr(\nabla^{2k})/k!\,. 
\eea 
One can show that every $ch_k(\nabla^2)$ 
is a closed form and that its de Rham 
cohomology class does not depend on the choice of a connection $\nabla$.
The corresponding classes are denoted by $ch_k(E)$ and lie in $H^{2k}_{dR}(M)$.  
Their sum is
the Chern
character $ch(E)$. The class $ch_1(E)$ is also known as the first Chern class
and is usually denoted by $c_1(E)$. 

The Chern character has the following
properties for two complex vector bundles $E_1$, $E_2$: 
\bea
ch(E_1\oplus E_2)&=&ch(E_1)+ch(E_2),\\
ch(E_1\ot E_2)&=&ch(E_1)ch(E_2). 
\eea 
(In the second formula
the multiplication on the right is the cup product in cohomology.)
\note{
Moreover (\cite{ah62}, \cite{k-m67}), the
associated map $ch:K^0(M)\ot\Q\ra H^{\mbox{even}}(M,\Q)$ is an
isomorphism that can be extended to a ring isomorphism
$K^*(M)\ot\Q\ra H^*(M,\Q)$ mapping $K^1(M)\ot\Q$ onto
$H^{\mbox{odd}}(M,\Q)$.
}
For the first Chern class, one has 
\beq\label{c1}
c_1(\nabla)=-\frac{1}{2\pi i}\Tr\nabla^2. 
\eeq 
\note{
One can show that the
integral of the $2j$-form $c_j(\nabla)$ over some $2j$-cycle 
has to be integer 
(sometimes called Chern number).
}
If   $L_1,\ldots,L_n$ are complex line bundles over a compact
space, then (see \cite[V.3.10]{k-m78})
\note{
$c_0(L)=1$, $c_1(L)=\chi(L)$ and $c_i(L)=0$ for
$i\neq 0,1$ (see \cite[V.3.15,3)]{k-m78}). Here, $\chi(L)$ is the
Euler class of the oriented vector bundle $L$. Since the Euler
class fulfills $\chi(L_1\ot L_2)=\chi(L_1)+\chi(L_2)$ for two line
bundles $L_1,~L_2$ , we obtain
}
\beq\label{c1ten} 
c_1(L_1\ot\cdots\ot L_n)=\sum_{k=1}^nc_1(L_k).
\eeq 
 For the Chern character of
a line bundle $L$, we have 
$ 
ch(L)=\exp(c_1(L))
$ (see \cite[V.3.23]{k-m78}).
If the bundle lives over a
2-manifold, all higher than linear terms in $ch(L)$ vanish, so that
in this case $ch(L)=1+c_1(L)$.

\subsection{Hopf fibration and Dirac monopole}

We begin this section with the description of a general topological
construction of
a principal bundle starting from a compact group $G$. When specialising
 to $G=U(1)$, we shall recover the Hopf fibration $S^3\ra S^2$.

Let $I=[0,1]$ be the closed unit interval and let $G$ be a compact group.
On the space $I\times G\times G$, consider the 
 equivalence relation where $R$ given by
\beq
\ba{ccc}
(0,g,h)\sim (0,g',h')&\Longleftrightarrow&h=h'\\
(t,g,h)\sim (t',g',h')&\Longleftrightarrow&t=t',~g=g',~h=h',~t\neq 0,1,\\
(1,g,h)\sim (1,g',h')&\Longleftrightarrow&g=g'.
\ea
\eeq
The quotient space 
$X:=(I\times G\times G)/R$ is like a partial suspension of 
$G\times G$, with the first copy of $G$ collapsed at the left end, and 
the second
copy collapsed at the right end. The space $X$ is compact
Hausdorff.
\note{
\bre
Note that we can realise the space $X$ in a slightly different manner:
The map $[0,1]\times G\times G\ra [0,1]\times G\times G$ given by
$(t,g,h)\mapsto (t,gh^{-1},h)$ is an equivariant homeomorphism, if we supply
the second copy of $[0,1]\times G\times G$ with the right $G$-action
$((t,g,h),k)\mapsto(t,g,hk)$. This action is compatible with the equivalence
relation $R'$ on $[0,1]\times G\times G$ given by 
$(0,g,h)\sim (0,g',h') \Longleftrightarrow h=h'$,
$(t,g,h)\sim(t',g',h')\Longleftrightarrow t=t', g=g', h=h',~t,t'\neq 0,1$,
$(1,g,h)\sim (1,g',h')\Longleftrightarrow gh=g'h'$. Thus, the space $X$ is
$G$-equivariantly homeomorphic to $X':=([0,1]\times G\times G)/R'$.
\ere
}
\vspace*{-0mm}\begin{figure}[h]
\[
\includegraphics[width=90mm]{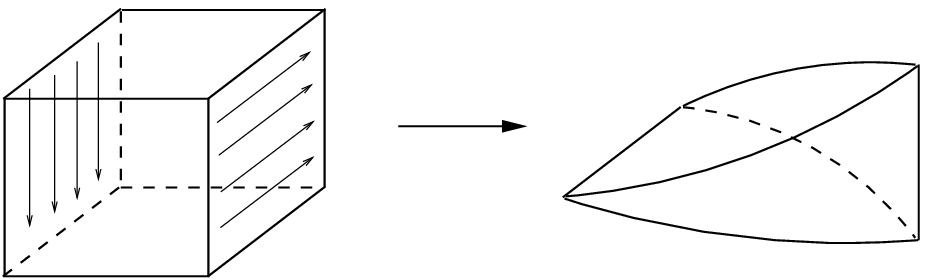}
\nonumber
\]
\end{figure}\vspace*{-0mm}

The diagonal action of $G$ on $I\times G\times G$ factorizes to the quotient, so that the formula
$
([(t,g,h)],k)\mapsto [(t,gk,hk)]
$
makes $X$ a right $G$-space. It is immediate that this action is
 principal. The quotient $X/G$ is homeomorphic to the suspension
$SG:=(I\times G)/(\{0\}\times G\cup \{1\}\times G)$ of $G$.
(Note that this notion of suspension differs from the one used in 
$K$-theory.)
Indeed, one can verify that the map $X/G\ra SG$ given by
\beq
[[(t,g,h)]]\longmapsto \left\{\ba{ccc} 
x_0:=[(0,k)]& \text{for}& t=0,\quad k\in G\\ 
\left[(t,gh^{-1})\right] & \text{for}& t\neq 0,1\\
x_1:=[(1,k)]& \text{for}& t=1,\quad k\in G
\ea\right.
\eeq
is a homeomorphism.
The composition of the canonical quotient map with this homeomorphism
defines a principal $G$-bundle $\pi:X\ra SG$. Explicitly, $\pi$ is given 
by
\beq
\pi([(t,g,h)])=\left\{\ba{ccc} 
x_0& \text{for}& t=0\\ 
\left[(t,gh^{-1})\right]& \text{for}& t\neq 0,1\\
x_1& \text{for}&t=1
\ea\right..
\eeq

 
\bpr
The principal bundle $(X,\pi,SG,G)$ is trivial if and only if $G$ is contractible.
\epr
\bpf
Suppose that the bundle is trivial. Then there exists a continuous
global section $\sigma:SG\ra X$.
This section is necessarily of the form:
\begin{align}
&\sigma(x_0)=[(0,g_0,h)],\\
&\sigma[(t,g)]=[(t,s(t,g),g^{-1}s(t,g))]\quad\text{for }t\neq 0,1,\\
&\sigma(x_1)=[(1,h,g_1)].
\end{align}
Here $s:\,]0,1[\times G\ra G$ is
a continuous map. The continuity of \hs\ implies that
$\lim_{[(t,g)]\ra x_0}s(t,g)=g_0$ and 
$\lim_{[(t,g)]\ra x_1}g^{-1}s(t,g)=g_1$. 
On the other hand, consider the following map:
\beq
H(t,g):=\left\{\ba{ccc} 
g_0g_1^{-1}& \text{for}& t=0\\ 
s(t,g)g_1^{-1}& \text{for}& t\in\, ]0,1[\\
g& \text{for}&t=1
\ea\right..
\eeq
It is evidently continuous at any point in $]0,1[\times G$. 
The continuity at any point $(0,k)$ follows from the first of the
above limits. To conclude the continuity at any point $(1,k)$, note
that $H$ is continuous if and only if $H'$ defined by 
$H'(t,g):=g^{-1}H(t,g)$ is continuous, and use the second limit.
Hence $H$ is a homotopy between $\id_G$ and the constant map
$g\mapsto g_0g_1^{-1}$.
The implication in the other direction is clear.
\epf

We now proceed to exhibit the local triviality of our bundle. 
Define $U_0:=\{[(t,g)]\in SG~|~t\in[0,1[\,\}$ and
$U_1:=\{[(t,g)]\in SG~|~t\in\,]0,1]\}$. They form an open cover
admitting local sections $\sigma_0:U_0\ra X$ and
$\sigma_1:U_1\ra X$ are defined by
\beq
\sigma_0(x_0)=[(0,e,h)],~~~\sigma_0([(t,g)])=[(t,e,g^{-1})],
\eeq
\beq
\sigma_1(x_1)=[(1,h,e)],~~~\sigma_1([(t,g)])=[(t,g,e)].
\eeq
These two local sections determine a transition function 
$\theta:U_0\cap U_1\cong\, ]0,1[\times G\ra G$ via (see the end
of Section~\ref{topo})
\beq\label{idtrans}
\theta([(t,g)])=
\check{\tau}\llp\sigma_0([(t,g)]),\hs_1([(t,g)])\lrp=g.
\eeq
Restricting to the cones $C_0:=\{[(t,g)]\in SG~|~t\in[0,\frac{1}{2}]\}$
and $C_1:=\{[(t,g)]\in SG~|~t\in[\frac{1}{2},1]\}$, we get piecewise 
trivialisations. The restriction of the transition function 
(\ref{idtrans}) to 
$C_0\cap C_1\cong G$ can be considered as the identity map of $G$.
Thus we arrive at the following geometrical picture:
The base space $SG$ is glued from the two cones $C_0$ and $C_1$,
and the total space $X$  is glued from the two trivial pieces 
$C_0\times G$ and $C_1\times G$ via the above transition function.
\note{
The Dirac monopole on the Hopf fibration is among the most popular
examples of connections on principal bundles. We recommend
\cite{n-gl97} and \cite{n-gl00} for a detailed study of this case.
}

Let us now specialise to the case $G=U(1)$, and see that we obtain
the Hopf fibration. To this end, first recall the usual description
of this bundle (see \cite{n-gl97}). 
Consider the 3-sphere $S^3$ as $\{(a,c)\in
\C^2~|~|a|^2+|c|^2=1\}$. The Hopf fibration is
obtained by the right action $((a,c),e^{i\phi})\mapsto
(ae^{i\phi},ce^{i\phi})$ of $U(1)$ on $S^3$. Identifying $S^3$
with $SU(2)$ by identifying $(a,c)$ with the matrix \beq
\begin{pmatrix}a&-\bar{c}\\c&\bar{a}\end{pmatrix},
\eeq we obtain the right action of $U(1)$ as the right
multiplication with ${\rm diag}(e^{i\phi},e^{-i\phi})$. The
quotient $S^3/U(1)=SU(2)/U(1)$ is diffeomorphic to the two-sphere
$S^2$.
The bundle projection $\pi:S^3\ra S^2$ is  explicitly given by
\beq\label{prs3s2} 
\pi(a,c)=(2a\bar{c},|a|^2-|c|^2). 
\eeq
 By Theorem \ref{lotri}, this bundle is a locally trivial principal bundle.

It is straigtforward to check that the following formulas define a 
$U(1)$-equivariant
homeomorphism $f:X=([0,1]\times U(1)\times U(1))/R\ra S^3$ and its inverse:
\bea
f[(t,g,h)]&=&\left(\sqrt{t}\ g,\sqrt{1-t}\ h\right)\,,\\
f^{-1}(a,c)&=&
\left[\left(|a|^2,\frac{a}{|a|},\frac{c}{|c|}\right)\right].
\eea
The suspension of $U(1)$ is identified with $S^2$ via the following homeomorphism
$\phi$ and its inverse:
\bea
\phi[(t,g)]&=&\left(2g\sqrt{t-t^2},2t-1\right),\\
\phi^{-1}(z,s)&=&\left[\left(\frac{s+1}{2},\frac{z}{|z|}\right)\right].
\eea
(Note that the endpoints $x_0$ and $x_1$ of the suspension of $U(1)$ are 
identified 
with the poles $(0,-1)$ and $(0,1)$
of~$S^2$, respectively.)

Next, consider $V_0=S^2\setminus\{(0,1)\}$ and $V_1=S^2\setminus \{(0,-1)\}$.
Then $\pi^{-1}(V_0)=\{(a,c)\in S^3~|~c\neq 0\}$ and $\pi^{-1}(V_1)=\{(a,c)\in
S^3~|~a\neq 0\}$,
and we get locally trivialising maps $\chi_0:\pi^{-1}(V_0)\ra V_0\times U(1)$ and
$\chi_1:\pi^{-1}(V_1)\ra
V_1\times U(1)$ by
\beq
\chi_0(a,c)=(\pi(a,c),\frac{c}{|c|}),~~\chi_1(a,c)=(\pi(a,c),\frac{a}{|a|}).
\eeq 
The corresponding transition function $\theta_{01}:V_0\cap
V_1\ra U(1)$ is given by 
\beq\label{transu1}
\theta_{01}(\pi(a,c))=\frac{a}{|a|}\frac{|c|}{c}. 
\eeq 

It is  straightforward to verify that in our concrete case $G=U(1)$ 
the open sets $V_0$ and $V_1$ are precisely
the images under $\phi$ of the open sets $U_0$ and $U_1$ defined above, that
the local trivialisations $\chi_0$, $\chi_1$ correspond to the sections
$\hs_0$, $\hs_1$, and that the transition functions (\ref{transu1}) and
(\ref{idtrans}) coincide up to the homeomorphism $\phi$.
Note that the restriction of $\theta_{01}$ to the equator $|a|=|c|$ is
the identity function, which is in accordance with the above remark on 
piecewise
triviality.

As is well-known, the 3-sphere can be considered as a gluing of
two solid tori along their boundaries. Indeed, $S^3$ can be split
into the union of $T_1=\{(a,c)\in S^3~|~|a|^2\in[0,\frac{1}{2}]\}$
and $T_2=\{(a,c)\in S^3~|~|a|^2\in[\frac{1}{2},1]\}$ (Heegaard
splitting, see Picture~\ref{hee}).
To prove that $T_1$ and $T_2$ are solid tori, we map $S^3$ 
onto 
\beq
 X:=\{(z_1,z_2)\in \C^2~|~|z_i|\leq 1,
(1-|z_1|^2)(1-|z_2|^2)=0\} 
\eeq
by a homeomorphism $g$ given by
\beq
g(a,c)=\frac{\sqrt{2}(a,c)}{\sqrt{1+||a|^2-|c|^2|}}. 
\eeq 
Its
inverse is defined by 
\beq
g^{-1}(z_1,z_2)=\frac{(z_1,z_2)}{\sqrt{1+|z_1|^2|z_2|^2}}. 
\eeq
See \cite{hms,bhms}  for proofs and details. It is
immediate to check that 
\beq g(T_1\cap T_2)=\{(z_1,z_2)\in\C^2~|~
|z_1|=|z_2|=1\}, 
\eeq \beq
g(T_1)=X_1:=\{(z_1,z_2)\in\C^2~|~|z_1|=1,|z_2|\leq 1\}=S^1\times
D^2, \eeq \beq g(T_2)=X_2:=\{(z_1,z_2)\in\C^2~|~|z_1|\leq 1,|z_2|=
1\}=D^2\times S^1. 
\eeq 
Thus we have shown that $S^3\cong
D^2\times S^1\bigcup S^1\times D^2$.
\begin{equation}
\label{hee}
\includegraphics[width=70mm]{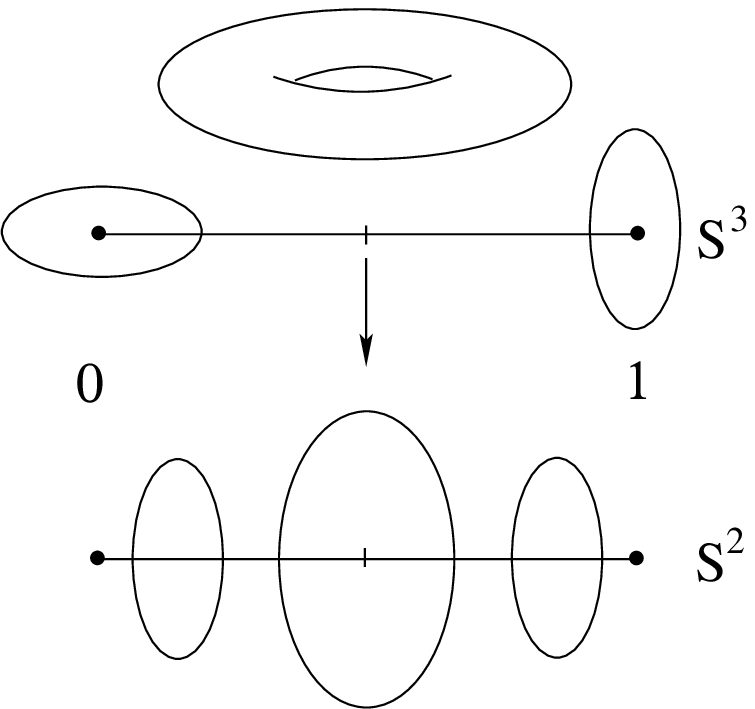}
\end{equation}

All the above can be written in a dual algebraic way, using
the coordinate functions $\ha$ and $\hg$, which assign to an
$SU(2)$-matrix its components $a$ and $c$.
Let us recall that the right action of $U(1)$ on $SU(2)$ gives
rise to a right coaction $\hD_R:C(SU(2))\ra C(SU(2))\ot C(U(1))$,
and that $C(S^2)$ can be identified with the subalgebra of
invariants, $C(S^2)=\{f\in C(SU(2))~|~\hD_R(f)=f\ot 1\}$. Note
also that the $*$-subalgebras $\cO(SU(2))\subset C(SU(2))$
(respectively $\cO(U(1))\subset C(U(1))$) can be identified with
the free $*$-algebras generated by $\ha,\hg$ modulo commutativity
and the relation $\ha^*\ha+\hg^*\hg=1$ (respectively by one element
$u$ modulo commutativity and the relation $u^*u=1$). The
comultiplication of $C(SU(2))$ is given on the generators by 
\beq
\hD(\ha)=\ha\ot\ha-\hg^*\ot\hg,~\hD(\hg)=\hg\ot\ha + \ha^*\ot\hg,
\eeq 
and the antipode and the counit by 
\beq
S(\ha)=\ha^*,~S(\hg)=-\hg,~\epsilon(\ha)=1,~\epsilon(\hg)=0. 
\eeq
The right coaction  on the generators is
\beq 
\hD_R(\ha)=\ha\ot
u,~\hD_R(\hg)=\hg\ot u. 
\eeq
Note that here all the maps (also the antipode) are $*$-homomorphisms, 
so that it suffices
to define them on generators.
Using the Poincar\'e-Birkhoff-Witt basis $\ha^k{\hg^*}^l\hg^m$,
${\hg^*}^p\hg^r{\ha^*}^s$, $k,l,m,p,r,s\in\N$, $k>0$, of
$\cO(SU(2))$, one can easily derive that the subalgebra
$\cO(S^2)\subset C(S^2)$ of invariant polynomials is spanned by
the elements $\ha^k{\hg^*}^l\hg^m$ with $k+m-l=0$ and the elements
${\hg^*}^p\hg^r{\ha^*}^s$ with $r-p-s=0$.

Our goal is now to describe vector bundles associated to the
$U(1)$-Hopf principal bundle, and to compute Chern numbers for
these bundles.  The
irreducible representations (all one-dimensional) of $U(1)$ are
parametrized by $\Z$ and given by $e^{i\phi}\mapsto e^{i\mu\phi}$,
$\mu\in\Z$. For any of these representations we have a line bundle
$L_\mu$ associated to the Hopf fibration.
Since the integer $\mu$ functions as the winding number of a 
representation defining $L_\mu$, we call it the {\em winding number} of the bundle $L_\mu$.
 Let $P_\mu$ denote the space of continuous
sections of $L_\mu$, i.e., $P_\mu:=\hG(L_\mu)$. Then we have the
following isomorphisms, coming from Lemma \ref{homma} and Theorem
\ref{secot}: 
\beq 
P_\mu\cong
\hom_{U(1)}(SU(2),\C_\mu)\cong\cC(SU(2))\!\!\coten{\cO(U(1))}\!\!
^\mu\C\,. 
\eeq
 Here $\C_\mu$ is the symbol for the space $\C$
viewed as a right $U(1)$-space with action $(g,z)\mapsto
g^{-\mu}z$, and $^\mu\C$ signifies the space $\C$ as a left
$\cO(U(1))$-comodule
with left coaction $_\C\hD(1)=u^{-\mu}\ot 1$.
By the first isomorphism, $P_\mu$ is identified with the
$C(S^2)$-submodule $\{f\in C(SU(2))~|~f(xg)=g^{-\mu}f(x)\}$ of
$C(SU(2))$, and by the second one with the $C(S^2)$-submodule $\{f\in
\cC(SU(2))~|~\hD_R(f)=f\ot u^{-\mu}\}$ of 
$\cC(SU(2))=\bigoplus_{\mu\in\Z} P_\mu$ (algebraic direct
sum). The last equality is immediate from the fact that the elements
$u^\mu,~\mu\in\Z$, form a linear basis of $\cO(U(1))$.

Our next aim is to obtain
explicit  idempotents determining the line bundles.
 We will follow the lines of
\cite{hm99}, where this was done for the quantum group $SU_q(2)$ 
(cf.\ \cite{l-g00}).
With the help of the Poincar\'e-Birkhoff-Witt basis
$\ha^k{\hg^*}^l\hg^m,\;{\hg^*}^p\hg^r{\ha^*}^s,\;
k,l,m,p,r,s\in\N,\; k>0$,
of $\cO(SU(2))$,
one can show that $f\in \cO(SU(2))\cap P_\mu$ is an element of
$\sum_{k=0}^{-\mu}\cO(S^2)\ \ha^{-\mu-k}\hg^{k}$
for $\mu\leq 0$, and of
$\sum_{k=0}^{\mu}\cO(S^2)\ {\hg^*}^{k}{\ha^*}^{\mu-k}$
for $\mu\geq 0$.
This suggests that we should have
\begin{eqnarray}\label{pmuc}
&& P_\mu=\left\{
\begin{array}{ll}
\sum_{k=0}^{-\mu}C(S^2)\ \ha^{-\mu-k}\hg^{k}
& \mbox{for $\mu\leq 0$}\\
\sum_{k=0}^{\mu}C(S^2)\ {\hg^*}^{k}{\ha^*}^{\mu-k}
& \mbox{for $\mu\geq 0$}.
\end{array}
\right.
\end{eqnarray}

It is obvious that the right hand side of (\ref{pmuc}) is
contained in $P_\mu$. In order to prove the reverse inclusion, let
$f\in P_\mu$, i.e., $\hD_R(f)=f\ot u^{-\mu}$. Let us consider the
case $\mu\leq 0$. Then the functions
$f_k:=f{\ha^*}^{-\mu-k}{\hg^*}^k$, $k=0,\ldots,-\mu$, are in
$C(S^2)$.
On the other hand, 
\beq
\sum_{k=0}^{-\mu}f_k\begin{pmatrix}-\mu\\k
\end{pmatrix}\ha^{-\mu-k}\hg^k=
f\sum_{k=0}^{-\mu}{\ha^*}^{-\mu-k}{\hg^*}^k\begin{pmatrix}-\mu\\k
\end{pmatrix}\ha^{-\mu-k}\hg^k=f,
\eeq 
because the last sum coincides with
$(m\ci(S\ot\id)\ci\hD)(\ha^{-\mu})=\epsilon(\ha^{-\mu})=1$. This
proves the first equation (\ref{pmuc}). The second equation is
derived similarly.

The foregoing argument gives us a key for finding idempotents.
Indeed (again for the case $\mu\leq 0$), the two column vectors
$\tilde{v}_\mu,\tilde{w}_\mu$ defined by 
\begin{align}
&
\tilde{v}_\mu^\top=\left({\ha^*}^{-\mu},\ldots,
{\scriptsize\begin{pmatrix}-\mu\\k\end{pmatrix}}
{\ha^*}^{-\mu-k}{\hg^*}^k,\ldots,{\hg^*}^{-\mu}\right), 
\\ &
\tilde{w}_\mu^\top=\llp{\ha}^{-\mu},\ldots,{\ha}^{-\mu-k}{\hg}^k,\ldots,
{\hg}^{-\mu}\lrp,
\end{align}
satisfy $\tilde{v}_\mu^\top \tilde{w}_\mu=1$. Therefore,
$(\tilde{w}_\mu \tilde{v}_\mu^\top)^2= \tilde{w}_\mu
\tilde{v}_\mu^\top \tilde{w}_\mu \tilde{v}_\mu^\top =\tilde{w}_\mu
\tilde{v}_\mu^\top$. 
For the case $\mu\geq 0$ we similarly use
\bea
&
\tilde{v}_\mu^\top=\left({\ha}^\mu,\ldots,
{\scriptsize\begin{pmatrix}\mu\\k\end{pmatrix}}{\ha}^{\mu-k}{\hg}^k,
\ldots, {\hg}^\mu\right),
\\ &
\tilde{w}_\mu^\top=\llp{\ha^*}^\mu,\ldots,{\ha^*}^{\mu-k}{\hg^*}^k,
\ldots,{\hg^*}^\mu\lrp.
\eea 
Let us put $\tilde{E}_\mu=\tilde{w}_\mu \tilde{v}_\mu^\top$.
The entries of this $(|\mu|+1)\times(|\mu|+1)$-matrix are elements
of $C(S^2)$.
\bpr\label{pn+1} 
The left $C(S^2)$-modules $P_\mu$ and
$C(S^2)^{|\mu|+1}\tilde{E}_\mu$ are isomorphic. 
\epr 
\bpf 
We
define a map $C(S^2)^{|\mu|+1}\tilde{E}_\mu\ra P_\mu$ by
$x^\top\tilde{E}_\mu\mapsto x^\top \tilde{w}_\mu$. This map is
well defined if and only if $x^\top(1-\tilde{E}_\mu)\mapsto 0$. 
To prove the latter, we observe that,  
due to $\tilde{v}_\mu^\top \tilde{w}_\mu=1$,
\[
x^\top(1-\tilde{w}_\mu \tilde{v}_\mu^\top) \longmapsto x^\top
\tilde{w}_\mu- x^\top \tilde{w}_\mu \tilde{v}_\mu^\top
\tilde{w}_\mu=0\,.
\]
 The surjectivity of this
map follows from (\ref{pmuc}). Moreover, from $x^\top
\tilde{w}_\mu=0$ we infer that $x^\top \tilde{w}_\mu \tilde{v}_\mu^\top
=x^\top\tilde{E}_\mu=0$, so that the map is injective. It is
obviously a homomorphism of left $C(S^2)$-modules. 
\hfill\epf
\bre\label{herm}
The idempotent $\tilde{E}_\mu$ is not hermitian, which is due to
the fact that the binomial coefficients appear only in $\tilde{v}_\mu$.
However, hermiticity can be achieved by a
similarity transformation. Let 
\beq
A_\mu:={\rm diag}\left(1,\ldots,\sqrt{\scriptsize\begin{pmatrix}|\mu|\\k
\end{pmatrix}}
,\ldots,1\right). 
\eeq 
Then the vectors
$w_\mu=A_\mu\tilde{w}_\mu$, $v_\mu^\top=\tilde{v}_\mu^\top
A_\mu^{-1}$
give rise to $E_\mu:=w_\mu v_\mu^\top 
=A_\mu\tilde{E}_\mu A_\mu^{-1}$. Hence, 
 from $v_\mu^\top=w_\mu^*\,$, we conclude that $E_\mu^*=E_\mu$.
\ere

Next, let us denote the real coordinates in $\C\times\R\cong \R^3$ by
$x_1,x_2,x_3$, and consider $x_1,x_2,x_3$ as continuous functions
on $S^2\subset \C\times\R$. Then the formula (\ref{prs3s2}) means
\[
\pi^*(1+x_3)=2\ha^*\ha,\quad
\pi^*(1-x_3)=2\hg^*\hg,\quad\pi^*(x_1+ix_2)=2\ha^*\hg,\quad
\pi^*(x_1-ix_2)=2\ha\hg^*.
\]
 Omitting the pull-back $\pi^*$ in these formulas we can write the following
explicit expressions for the projectors $E_{\pm 1}$:
\beq\label{bott-1}
E_{-1}=\tilde{E}_{-1}=\left(\ba{cc}\ha^*\ha & {\ha}{\hg^*}\\
{\ha}^*{\hg}& \hg^*\hg \ea\right)=
\frac{1}{2}\left(\ba{cc}
1+x_3&x_1-ix_2\\
x_1+ix_2&1-x_3 \ea\right)
\eeq
and
\beq\label{bott1}
E_1=\tilde{E}_1=\left(\ba{cc}\ha^*\ha & {\ha^*}{\hg}\\
{\ha}{\hg^*}& \hg^*\hg \ea\right)=
\frac{1}{2}\left(\ba{cc}
1+x_3&x_1+ix_2\\
x_1-ix_2&1-x_3 \ea\right).
\eeq

Let us now see which of these idempotents corresponds to the
tautological line bundle. Consider $S^2$ as $\C P^1$ with homogeneous 
coordinates $[z_1:z_2]$. Then the fibre of the tautological line bundle 
at $[z_1:z_2]$ is the 
complex line $\{(\lambda z_1,\lambda z_2)\,|\, \lambda\in\C\}$, 
and
\beq
\frac{1}{|z_1|^2+|z_2|^2}\left(\ba{c}z_1^*\\z_2^*\ea\right)
\left(\ba{cc}z_1 &z_2\ea\right)
\eeq
is an idempotent defining the left module of sections of this
bundle. Restricting to $S^3\subset \C^2$, we have 
$(|z_1|^2+|z_2|^2)^{-1/2}z_1=\ha$, $(|z_1|^2+|z_2|^2)^{-1/2}z_2=\hg$,
and the above idempotent coincides with $E_1$. This means that the tautological
line bundle coincides with the associated line bundle $L_1$ coming from the
identity representation.

\wegdamit{ Now, formula (\ref{i}) gives us a splitting
$s:\slq\ra\asq\ot\slq$, and we can claim: \bpr\label{ml} Put
\begin{eqnarray*}\label{ekl}
&& (e_n)_{kl}=\left\{
\begin{array}{ll}
\ha^{-n-k}\hg^{k}{\scriptsize\left(\ba{c}-n\\
l\ea\right)}{\hg^*}^l{\ha^*}^{-n-l}
& \mbox{for $n\leq 0$}\\
\hb^{k}\hd^{n-k}{\scriptsize\left(\ba{c}n\\ l\ea\right)}
(-1)^{-l}\ha^{n-l}\hg^l & \mbox{for $n\geq 0$}.
\end{array}
\right.
\end{eqnarray*}
Then, for any $n\in\IZ$, $e_n\in M_{|n|+1}(\asq)$, $e_n^2=e_n$,
and $\asq^{|n|+1}e_n$ is isomorphic to $P_n$ as a left
\asq-module. \epr\bpf
Observe that for $n\geq 0$ we can write $e_n=uv^T$, where
$u^T=({\ha^*}^n,...,(-{\hg^*})^k{\ha^*}^{n-k},...,(-\hg^*)^n)$ and
$v^T=(S({\ha^*}^n),...,
{\scriptsize\left(\ba{c}n\\
k\ea\right)}S(\hg^k{\ha^*}^{n-k}),...,S(\hg^n))$. Since
\[
v^Tu =\sum_{k=0}^{n}
{\scriptsize\left(\ba{c}n\\
k\ea\right)}S(\hg^k{\ha^*}^{n-k})(-\hg^*)^k{\ha^*}^{n-k}
=S(({\ha^*}^n)\1)({\ha^*}^n)\2=\he({\ha^*}^n)=1,
\]
we can directly see that $e_n^2=e_n$. The case $n\leq 0$ is
similar.
An isomorphism $P_n\ra \asq^{|n|+1}e_n$ of left \asq-modules can
be given as follows. }

From the proof of the Serre-Swan Theorem, together with the
fact that $S^2$ is covered by two trivialising neighbourhoods, we
know that there exist idempotents in $M_2(C(S^2))$ defining the
modules $P_\mu$. (See
 \cite[p.34--35]{k-m78} for a general formula for an idempotent
in terms of a partition of unity and transition functions subordinate to 
local trivialisations, and \cite[p.77]{fgv01} for concrete formulas.) 
Let us   give now some heuristic arguments
leading to such idempotents. From (\ref{pmuc}) we know that
$\ha^{-\mu},\ldots,\ha^{-\mu-k}\hg^k,\ldots,\hg^{-\mu}$ are
generators of the $C(S^2)$-module $P_\mu$ (for $\mu\leq 0$). If we
want to choose only two of these as generators, we immediately see that
the only possiblity is to take $\ha^{-\mu}$ and $\hg^{-\mu}$,
because any other pair has common zeros. Now, for sure we have
$C(S^2)\ha^{-\mu} + C(S^2)\hg^{-\mu}\subset P_\mu$. To prove the
reverse inclusion, it is sufficient to find $g_1,~g_2\in
P_{-\mu}$ 
such that 
$g_1\ha^{-\mu}+g_2\hg^{-\mu}=1$. Indeed, then 
$f\in P_\mu$ can be written as $f=fg_1\ha^{-\mu}+fg_2\hg^{-\mu}$, with $fg_1,fg_2\in C(S^2)$.
It is
immediate that 
\beq
g_1=\frac{{\ha^*}^{-\mu}}{|\ha|^{-2\mu}+|\hg|^{-2\mu}}
\quad\text{and}\quad
g_2=\frac{{\hg^*}^{-\mu}}{|\ha|^{-2\mu}+|\hg|^{-2\mu}}
\eeq
satisfy these requirements.
Thus we have the equality 
$C(S^2)\ha^{-\mu} + C(S^2)\hg^{-\mu}= P_\mu\,$.

We can also prove an analogue of Proposition~\ref{pn+1}
for $2\times 2$-idempotents coming from the above considerations.
Furthermore, we can argue as in Remark~\ref{herm} to
 arrive at a hermitian
$2\times 2$-idempotent. It can be constructed as follows.
%
%
%
Define column
vectors $r_\mu$ by 
\beq\label{p-1}
r_\mu^*=
\frac{1}{\sqrt{|\ha|^{2|\mu|}+|\hg|^{2|\mu|}}}
\left\{\ba{ccc}
\left({\ha^*}^{-\mu},
~{\hg^*}^{-\mu}\right)&
\mbox{for}& \mu< 0\\
\left({\ha}^{\mu},
~{\hg}^{\mu}\right)&
\mbox{for}&\mu> 0
\ea\right. .
\eeq 
They satisfy $r_\mu^*r_\mu=1$
and give a hermitian idempotent $p_\mu:=r_\mu r_\mu^*$
such that
$P_\mu\cong C(S^2)^2p_\mu\,$.
%
%
We have the following explicit formulas for the projectors $p_\mu$
(the pull-back $\pi^*$ omitted):
\hspace*{-5mm}
\bea\label{p-1'}
p_\mu&=&\frac{1}{(\ha^*\ha)^{|\mu|}+(\hg^*\hg)^{|\mu|}}
\left(\ba{cc}(\ha^*\ha)^{|\mu|}
 & 
({\ha}{\hg^*})^{|\mu|}\\
({\ha^*}{\hg})^{|\mu|} &
(\hg^*\hg)^{|\mu|} \ea\right)\\
& =&
\frac{1}{(1+x_3)^{|\mu|}+(1-x_3)^{|\mu|}}
\left(\ba{cc}
(1+x_3)^{|\mu|}&(x_1-ix_2)^{|\mu|}\\
(x_1+ix_2)^{|\mu|}&(1-x_3)^{|\mu|}
\ea\right),~\mu< 0,~~~~ \\[.3cm]
p_\mu&=&\frac{1}{(\ha^*\ha)^{\mu}+(\hg^*\hg)^{\mu}}\left(\ba{cc}(\ha^*\ha)^{\mu}
 & 
({\ha^*}{\hg})^{\mu}\\
({\ha}{\hg^*})^{\mu} &
(\hg^*\hg)^{\mu} \ea\right)\\
& =&
\frac{1}{(1+x_3)^{\mu}+(1-x_3)^{\mu}}
\left(\ba{cc}
(1+x_3)^{\mu}&(x_1+ix_2)^{\mu}\\
(x_1-ix_2)^{\mu}&(1-x_3)^{\mu}
\ea\right),\mu> 0. 
\eea
%
In particular we obtain 
$
p_{\pm 1}=E_{\pm 1}.
$
Finally, one can see directly that the projectors $E_\mu$ and
$p_\mu$ are Murray-von Neumann equivalent for any $\mu$. Indeed, the
$(|\mu|+1)\times 2$-matrix $F_\mu:=r_\mu w_\mu^*$ has entries from
$C(S^2)$ and satisfies $F_\mu F_\mu^*=p_\mu$ and $F_\mu^*
F_\mu=E_\mu$. 

Our next and main goal  is to compute the Chern class $c_1(L_\mu)$
for any integer~$\mu$.
First, let us make some general remarks about computations
with $2\times 2$-idempotents. Any $2\times 2$-matrix $e$ with entries from an algebra $A$ 
can
be written uniquely as $e=\sum_{k=0}^3s_k\hs_k$, where $\hs_k$ are the
Pauli matrices
\beq
\hs_0=\left(\ba{cc}1&0\\0&1\ea\right),\quad
\hs_1= \left(\ba{cc}1&0\\0&-1\ea\right),\quad
\hs_2= \left(\ba{cc}0&1\\1&0\ea\right),\quad
\hs_3= \left(\ba{cc}0&-i\\i&0\ea\right),
\eeq
and $s_k\in A$. The matrices $\hs_k$, $k=1,2,3$, are traceless and satisfy
$\hs_k^2=$I and
$
\hs_j\hs_k=i\varepsilon_{jkl}\hs_l$. Here $(j,k,l)$ is
 any cyclic permutation of $(1,2,3)$, and 
 $\varepsilon_{jkl}$ is the completely antisymmetric tensor.
Assume now that $A:=C(M)$ is the algebra
 of continuous functions on a compact Hausdorff
space~$M$. It is immediate from these formulas
and the linear independence of the $\hs_k$'s that the condition $e^2=e$
is equivalent to the equations
\beq
s_0=s_0^2+s_1^2+s_2^2+s_3^2,~~s_1=2s_0s_1,~~s_2=2s_0s_2,~~s_3=2s_0s_3.
\eeq
Moreover, it follows that $s_0=\frac{1}{2}$ at every point of $M$ 
where at least one of
the $s_k$, $k=1,2,3$, does not vanish. (We disregard the trivial
case $e=0$.) Also, it follows that $s_0=1$ at every point where
all the other $s_k$ vanish. Thus, if $M$ is connected, either $s_0=1$ 
and $s_k=0$, $k=1,2,3$, everywhere (i.e., $e=I$), or $s_0=\frac{1}{2}$ and
$s_1^2+s_2^2+s_3^3=\frac{1}{4}$. If  $M$ is a 
manifold and the entries
of $e$ are differentiable, then
\beq
\d e=\d s_1\hs_1+\d s_2\hs_2+\d s_3\hs_3\,.
\eeq
An easy calculation making use of the above properties of the Pauli matrices
(in particular regarding the trace), leads to
\beq\label{trps}
\Tr (e\,\d e\,\d e)=4i\ (s_1\d s_2\land \d s_3+s_2\d 
s_3\land \d s_1+s_3\d s_1\land \d s_2).
\eeq
 
On the other hand,
 for the line bundles over the 2-sphere,
let us notice that
$L_\mu$ is the tensor product of $|\mu|$ copies of $L_1$ (for
$\mu> 0$) or $L_{-1}$ (for $\mu< 0$).
Also, $L_\mu\ot L_{-\mu}\cong L_0$ for any $\mu\in\Z$ because
the tensor product of a
representation and its conjugate is the trivial representation. 
 Hence, it follows from  (\ref{c1ten}) that
$c_1(L_\mu)=\mu c_1(L_{1})$. 
 In view of (\ref{c1}) and (\ref{ogras}), the right hand side is
determined by
$c_1(L_{1})=\frac{1}{2\pi i}\Tr(E_{1}\d E_{1}\d E_{1})\,$. Furthermore,
from (\ref{bott1}) we have
\beq
E_{1}=\frac{1}{2}(I+x_3\hs_1+x_1\hs_2-x_2\hs_3),
\eeq
so that
$
s_1=\frac{1}{2}x_3,\; s_2=\frac{1}{2}x_1,\; s_3=-\frac{1}{2}x_2\,
$.
 Putting this into (\ref{trps}), we immediately
obtain
\beq
c_1(L_{1})=-\frac{1}{4\pi}(x_1\d x_2\land \d x_3+x_2 \d x_3\land \d x_1
+ x_3 \d x_1\land \d x_2). 
\eeq 
Integrating $c_1(L_{1})$
 in the spherical coordinates (for the orientation with 
outward-pointing normal vector) gives~$-1$. 
Consequently,
$
\int_{S^2}c_1(L_{\mu})=-\mu 
$.

In the framework of noncommutative geometry,  the integration
 of $c_1(L_\mu)$ can be understood as
a pairing between the cyclic 2-cocycle given by
$\phi_{S^2}(f_0,f_1,f_2):=\frac{1}{2\pi i}
\int_{S^2}f_0\d f_1\wedge\d f_2$
and the $K_0$-class of the projector
$E_\mu$. 
Therefore, we can summarise  our computation as:
\bth
The pairing of the $K_0$-class of the projector $E_\mu$ 
and the cyclic cohomology class of the cyclic
2-cocycle $\phi_{S^2}$ (fundamental class) is equal to minus the 
winding number of the line bundle $L_\mu\,$:
\beq
\langle [E_\mu],[\phi_{S^2}]\rangle =
\frac{1}{2\pi i}\int_{S^2}\Tr (E_\mu \d E_\mu\land \d E_\mu)= 
\int_{S^2}c_1(L_\mu)=-\mu.
\eeq
\ethe
\noindent

We proceed now to the discussion of
connections. In the general situation of a connected Lie group $G$
with a closed subgroup $H$ giving rise to a prinicipal bundle
$\pi:G\ra G/H$ with the structure group $H$, it is natural to look for
connections that are invariant under the left action of $G$ on
itself. Such connections are characterized by 
$\mathrm{ad}(H)$-invariant
complements of the Lie algebra $\frak{h}$ of $H$ in the Lie
algebra $\frak{g}$ of $G$ (see \cite[II, Theorem~11.1]{kn63}).
Equivalently, a $G$-invariant connection is given by a projection
$\frak{g}\ra\frak{h}$ that commutes with $\mathrm{ad}(H)$. In the
semisimple case, a complement can be defined by orthogonality with
respect to the Killing metric. The corresponding connection is
then called canonical.

In the dual formulation in terms of function algebras, this can be
rephrased as an appropriate splitting of the pull-back of the
imbedding of a closed subgroup into a compact group. To be more precise, 
let $G$ be a compact  group with
closed subgroup $H$, and let $p:\cO(G)\ra\cO(H)$ be the pull-back of
the inclusion map $H\ra G$. Denote by $\hD^G$ and $\hD^H$ the coproducts
 of the
Hopf algebras $\cO(G)$ and $\cO(H)$, respectively.
 If $i:\cO(H)\ra\cO(G)$ is a unital linear
splitting of $p$, then the desired invariance properties of $i$
can be given as follows:
\beq
 (i\ot\id)\ci\hD^H=(\id\ot p)\ci\hD^G\ci i,
\quad(\id\ot i)\ci\hD^H=(p\ot\id)\ci\hD^G\ci i\,. 
\eeq 
We want to show that the connection defining the
Dirac monopole can be obtained along these lines.

To this end, let us imbed $U(1)$ in $SU(2)$ as the subgroup of
diagonal matrices, so that its pull-back reads $p(\ha)=u$, $p(\hg)=0$.
Now, observe that the unital linear map 
$i:\cO(U(1))\ra\cO(SU(2))$
defined by 
\beq\label{csplit}
 i(u^n)=\alpha^n,\quad i({u^*}^{n})={\ha^*}^n,\quad n\in\N, 
\eeq
enjoys all the aforementioned properties (see~\cite[2.58]{dgh01}). 
Starting from $i$, we can
define  $\ho:\cO(U(1))\ra \Lambda^1(SU(2))$ by
$\ho=m\ci(S\ot\d)\ci\hD\ci i$.  
Here  $\d$ denotes the de Rham differential on
$SU(2)$ and $m$ is
the multiplication. In particular, we have
\beq\label{dirpol} 
\ho(u)=\ha^*\d\ha + \hg^*\d\hg. 
\eeq 
We want to prove now that \ho\ is a connection form in the sense
of Definition~\ref{coform}.
First, \ho\ is constructed as the canonical projection on
the de Rham differential forms of a universal-calculus connection
form (see \cite[Example~2.14]{dgh01}). It follows from the general
theory of such connection forms that their projections on the de Rham 
calculus
always satisfy the Leibniz rule
\cite[(2.62)]{dgh01} that is  the first defining property
 of a connection form.
Therefore, $\ho$ is determined by its values on
$u$. With this in mind, a verification of the remaining
defining properties is straightforward.

On the other hand, it is known that the traditionally defined
 Dirac monopole
connection form comes from the canonical invariant connection
(see \cite[p.110]{kn63}), and
 is given by the formula (cf.\ \cite[pp.38-40]{n-gl00})
\beq\label{tdirpol} 
\tilde{\ho}=\ha^*\d\ha+\hg^*\d\hg\,. 
\eeq 
The very form of equations (\ref{dirpol})
and (\ref{tdirpol}) already strongly suggests that $\ho$
and $\widetilde{\ho}$ correspond to each other in the sense of 
Proposition~\ref{cofot}. 
Taking into account that both \ho\ and $\widetilde{\ho}$
satisfy the Leibniz rule, it suffices to
verify the equation of Proposition~\ref{cofot}
 for $h=u$, which is straightforward. 

Finally, let us  see that the \ho-defined
covariant derivative $\nabla$ on $L_{-1}$ 
coincides with the Gra{\ss}mann connection given by the projector 
$p_{-1}$. More precisely, on one hand side the connection form
\ho\ determines a covariant differentiation $D$ by 
Proposition~\ref{cofot}, and the latter induces a covariant derivative
$\nabla$ by the formula (\ref{nablaco}). On the other hand,
the dual bases $e_1:=\ha$, $e_2:=\hg$,  $e^1:=\ha^*$, $e^2:=\hg^*$, 
of $P_{-1}$ define the idempotent $p_{-1}$ by $e^j(e_i)=(p_{-1})^j_i$ 
(see \eqref{p-1} and \eqref{p-1'}). 

On any $f\in\cO(SU(2))$ we have the right coaction of the Hopf
algebra $\cO(U(1))$, so that  (\ref{jfd}) can be
 written as $Df=\d f-f_{(0)}\ho(f_{(1)})$. Hence, 
for the values of $D$ on the generators , we obtain
\beq
D\ha=\d\ha-\ha\ha^* \d\ha-\ha\hg^*\d\hg=\hg^*\hg \d\ha-\ha\hg^*\d\hg=
\hg \d(\ha\hg^*)-\ha \d(\hg^*\hg).
\eeq
Similarly, $
D\hg=\ha \d(\ha^*\hg)-\hg \d(\ha^*\ha)$.
The right hand sides of these two 
formulas are already written as elements
of $C^\infty(S^3)\pi^*(\Lambda^1(S^2))$, 
so that the isomorphism
$\Phi$ defined in (\ref{liot}) takes very simple form.
Plugging it into the formula (\ref{nablaco}) yields:
\bea
\nabla(\ha)&=&
\d(\ha\ha^*)\!\!\underset{C^\infty(S^2)}{\ot}\!\!\ha + 
\d(\ha\hg^*)\!\!\underset{C^\infty(S^2)}{\ot}\!\!\hg\,,\\
\nabla(\hg)&=&
\d(\hg\ha^*)\!\!\underset{C^\infty(S^2)}{\ot}\!\!\ha + 
\d(\hg\hg^*)\!\!\underset{C^\infty(S^2)}{\ot}\!\!\hg\,.
\eea
This coincides with the Gra{\ss}mann connection 
$\nabla^e(e_k)=\d (p_{-1})^j_k\ot_{C^\infty(S^2)}e_j$.
An analogous reasoning can be carried out for the tautological
line bundle $L_1$.

The considerations of this section manifest a well-known fact that
the $K_0$-invariants of vector bundles associated to a principal
bundle can be computed using constructions originating from the
principal bundle. In the setting of noncommutative geometry, 
 strong connections on a principal extension are used
to obtain explicit idempotents representing  finitely
generated projective modules associated to the principal extension.
Further along this line, the Chern-Galois character is used to
produce a cyclic homology class out of a quantum-group representation
defining the associated module~\cite{bh04}.
 Pairing this class with a cyclic
cocycle replaces the integration of characteristic classes given
in terms of the de Rham cohomology.

\section*{Acknowledgements}

The authors are very
grateful to   Tomasz Brzezi\'nski, Ryszard Engelking, Bogna Janisz, Max Karoubi,
Pawe\l\ \L.\ Kasprzak, Ulrich Kr\"ahmer,
Tomasz Maszczyk, Ryszard Nest, Yorck Sommerh\"auser,
Joseph C.\ V\'arilly,  Pawe\l\ Witkowski, Mariusz Wodzicki,
Stanis\l aw L.\ Woronowicz, and Bartosz Zieli\'nski
 for very helpful discussions or technical support.
This work was partially supported by the
European Commission
  grants MERG-CT-2004-513604 (PMH), RITA-CT-2004-505493 (PMH),
 MKTD-CT-2004-509794 (PFB, RM, WS),
  the KBN  grants 1 P03A 036 26 (PMH, RM, WS),
115/E-343/SPB/6.PR UE/DIE 50/2005-2008 (PMH),
and the Mid-Career Academic Grant from the Faculty of
Science and IT, the University of Newcastle (WS).

\fontsize{11pt}{12.5pt}\selectfont

\end{document}